\documentclass{amsart}
\usepackage{hyperref}

\usepackage[all,2cell]{xy}
\usepackage{amssymb,amsmath,stmaryrd,mathrsfs,amsthm}
\usepackage{graphics,color,epsfig}
\UseTwocells

\newcommand{\C}{\mathbb{C}}

\newcommand{\scrO}{\ensuremath{\mathscr{O}}}

\newcommand{\scrS}{\ensuremath{\mathscr{S}}}

\newcommand{\V}{\ensuremath{V}}

\newcommand{\bbM}{\ensuremath{\mathbb{M}}}

\newcommand{\bbR}{\ensuremath{\mathbb{R}}}
\newcommand{\bbS}{\ensuremath{\mathbb{S}}}
\newcommand{\bbT}{\ensuremath{\mathbb{T}}}
\newcommand{\bbZ}{\ensuremath{\mathbb{Z}}}

\newcommand{\ftil}{\ensuremath{\widetilde{f}}}

\newcommand{\fhat}{\ensuremath{\widehat{f}}}

\newcommand{\ten}{\otimes}

\newcommand{\inv}{^{-1}}

\newcommand{\maps}{\colon}
\newcommand{\meet}{\wedge}
\newcommand{\join}{\vee}

\newcommand{\ph}{\varphi}

\newcommand{\Cat}{\ensuremath{\mathbf{Cat}}}
\newcommand{\Gpd}{\ensuremath{\mathbf{Gpd}}}
\newcommand{\Grp}{\ensuremath{\mathbf{Grp}}}

\newcommand{\Set}{\ensuremath{\mathbf{Set}}}

\newcommand{\VCat}{\ensuremath{\V\text{-}\mathbf{Cat}}}

\newcommand{\im}{\ensuremath{\operatorname{im}}}

\newcommand{\Aut}{\ensuremath{\operatorname{Aut}}}

\newcommand{\Gal}{\ensuremath{\operatorname{Gal}}}
\newcommand{\AUT}{\ensuremath{\operatorname{AUT}}}

\newcommand{\Sh}{\ensuremath{\operatorname{Sh}}}

\newcommand{\op}{\ensuremath{^{\mathit{op}}}}
\newcommand{\adj}{\dashv}

\newdir{ >}{{}*!/-10pt/@{>}}
\newcommand{\iso}{\cong}
\newcommand{\eqv}{\simeq}
\newcommand{\too}[1][]{\ensuremath{\overset{#1}{\longrightarrow}}}
\renewcommand{\to}{\ensuremath{\rightarrow}}
\def \stackto #1{ \, {\stackrel{#1}{\longrightarrow}}\, }
\newcommand{\Impl}{\ensuremath{\Rightarrow}}

\newtheorem{thm}{Theorem}

\newtheorem{hyp}[thm]{Hypothesis}
\newtheorem{fact}[thm]{Fact}
\theoremstyle{definition}

\newtheorem{defn}[thm]{Definition}

\newtheorem{eg}[thm]{Example}

\newcommand{\cxymatrix}[1]{\vcenter{\xymatrix{#1}}}

\input xy \xyoption{all} \xyoption{poly} 
\labelmargin-{1.2pt}  

\newcommand{\thetagraph}[3]  
{ \xymatrix{ *{\bullet} \ar@{-} @/^1.5pc/ [r] ^{#1} \ar@{-}
@/_1.5pc/ [r] _{#3} \ar@{-} [r]^{#2} & *{\bullet}
\\}
}

\newcommand{\fourtheta}[4]  
{ \cxymatrix{ *{\bullet} \ar@{-} @/^1.5pc/ [r] ^{#1} \ar@{-}
@/_1.5pc/ [r] _{#4} \ar@{-} @/^/ [r]^{#2} \ar@{-} @/_/ [r]_{#3} &
*{\bullet}
\\}
}

\newcommand{\TetJ}[6]{
\def\lab{\ifcase\xypolynode\or #1 \or #2 \or #3 \fi}
\begin{xy}
\xygraph{!{<3.2pc,0pc>:}
 *{\bullet}
 !P3"A"{~><{@{{}{-}*{\bullet}}} ~>>{_{\lab}}}
 "A0" -@-_{#4} "A1"
 "A0" -@-_{#5} "A2"
 "A0" -@-^{#6} "A3"
}
\end{xy}
}

\newcommand{\TenJ}{
\def\lab{\ifcase\xypolynode\or 12 \or 23 \or 34 \or 45 \or 15 \fi}
\begin{xy}
\xygraph{!{<4pc,0pc>:}
 !P5"A"{~><{@{{}{-}*{\bullet}}} ~>>{_{j_{\lab}}}}
 "A1" -@-_{j_{13}} "A3"
 "A2" -@-_{j_{24}} "A4"
 "A3" -@-_{j_{35}} "A5"
 "A4" -@-_{j_{14}} "A1"
 "A5" -@-_{j_{25}} "A2"
}
\end{xy}
}

\newcommand{\TenA}{
\def\lab{\ifcase\xypolynode\or 12 \or 23 \or 34 \or 45 \or 15 \fi}
\begin{xy}
\xygraph{!{<4pc,0pc>:}
 !P5"A"{~><{@{{}{-}*{\bullet}}} ~>>{_{a_{\lab}}}}
 "A1" -@-_{a_{13}} "A3"
 "A2" -@-_{a_{24}} "A4"
 "A3" -@-_{a_{35}} "A5"
 "A4" -@-_{a_{14}} "A1"
 "A5" -@-_{a_{25}} "A2"
}
\end{xy}
}

\newcommand{\Ten}{%
\def\lab{\ifcase\xypolynode\or 1,2 \or 2,3 \or 3,4 \or 4,5 \or 5,1 \fi}
\begin{xy}
\xygraph{!{<2pc,0pc>:}
 !P5"A"{~><{@{{}{-}*{\bullet}}} ~>>{_{}}}
 "A1" -@-_{} "A3"
 "A2" -@-_{} "A4"
 "A3" -@-_{} "A5"
 "A4" -@-_{} "A1"
 "A5" -@-_{} "A2"
}
\end{xy}
}

\newcommand{\TenL}{
\def\lab{\ifcase\xypolynode\or 1,2 \or 2,3 \or 3,4 \or 4,5 \or 5,1 \fi}
\begin{xy}
\xygraph{!{<4pc,0pc>:}
 !P5"A"{~><{@{{}{-}*{\bullet}}} ~>>{_{\smult}}}
 "A1" -@-_{\smult} "A3"
 "A2" -@-_{\smult} "A4"
 "A3" -@-_{\smult} "A5"
 "A4" -@-_{\smult} "A1"
 "A5" -@-_{\smult} "A2"
}
\end{xy}
}

\newcommand{\tenl}{
\def\lab{\ifcase\xypolynode\or 1,2 \or 2,3 \or 3,4 \or 4,5 \or 5,1 \fi}
\begin{xy}
\xygraph{!{<2pc,0pc>:}
 !P5"A"{~><{@{{}{-}*{\bullet}}} ~>>{_{\smult}}}
 "A1" -@-_{\smult} "A3"
 "A2" -@-_{\smult} "A4"
 "A3" -@-_{\smult} "A5"
 "A4" -@-_{\smult} "A1"
 "A5" -@-_{\smult} "A2"
}
\end{xy}
}

\newcommand{\periodictable}{
\begin{center}
 { \textbf{
\begin{tabular}{|c|c|c|c|}  \hline
        & $\mathbf{\mathit n = 0}$ & $\mathbf{\mathit n = 1}$ &
$\mathbf{\mathit n = 2}$\\ \hline $\mathbf{\mathit k = 0}$ & sets
& categories & 2-categories     \\     \hline
$\mathbf{\mathit k = 1}$  & monoids   & monoidal   & monoidal         \\
        &           & categories & 2-categories     \\     \hline
$\mathbf{\mathit k = 2}$  &commutative& braided    & braided          \\
        & monoids   & monoidal   & monoidal         \\
        &           & categories & 2-categories     \\     \hline
$\mathbf{\mathit k = 3}$  &`'         & symmetric  & sylleptic \\
        &           & monoidal   & monoidal         \\
        &           & categories & 2-categories     \\     \hline
$\mathbf{\mathit k = 4}$  &`'         & `'         & symmetric \\
        &           &            & monoidal         \\
        &           &            & 2-categories     \\     \hline
$\mathbf{\mathit k = 5}$  &`'         &`'          & `'               \\
        &           &            &                  \\
        &           &            &                  \\     \hline
$\mathbf{\mathit k = 6}$  &`'         &`'          & `'               \\
        &           &            &                  \\
        &           &            &                  \\     \hline
\end{tabular}}}
\vskip 0.5em
\end{center}}

\newcommand{\extendedperiodictable}{
\begin{center}
 { \textbf{
\begin{tabular}{|c|c|c|c|c|c|}  \hline
        & $\mathbf{\mathit n = -2}$ & $\mathbf{\mathit n = -1}$ &
         $\mathbf{\mathit n = 0}$ & $\mathbf{\mathit n = 1}$ &
$\mathbf{\mathit n = 2}$\\ \hline 
$\mathbf{\mathit k = 0}$ &?  & ?
& sets & categories & 2-categories     \\     \hline
$\mathbf{\mathit k = 1}$  
& `' & ?
& monoids   & monoidal   & monoidal         \\
& &       &           & categories & 2-categories     \\     \hline
$\mathbf{\mathit k = 2}$  
& `'  &  `'
&commutative& braided    & braided          \\
& &        & monoids   & monoidal   & monoidal         \\
& &        &           & categories & 2-categories     \\     \hline
$\mathbf{\mathit k = 3}$  
& `'  &  `'
&`'         & symmetric  & sylleptic \\
& &        &           & monoidal   & monoidal         \\
& &        &           & categories & 2-categories     \\     \hline
$\mathbf{\mathit k = 4}$  
& `'  &  `'
&`'         & `'         & symmetric \\
& &        &           &            & monoidal         \\
& &        &           &            & 2-categories     \\     \hline
\end{tabular}}}
\vskip 0.5em
\end{center}}

\newcommand{\periodictableII}{
{
\[ %
\xy (-75,60)*{k=0}; (-75,30)*{k=1}; (-75,0)*{k=2};
(-75,-30)*{k=3}; (-75,-60)*{k=4}; (-40,75)*{n=0}; (0,75)*{n=1};
(40,75)*{n=2};
(0,0)*{ 
  \xy 0;/r.20pc/: 
  (25,25)*{};
(-2,12)*{\bullet}="1"+(-2,1)*{ \scriptstyle x};
(1,16)*{\bullet}="2"+(-2,2)*{\scriptstyle x^{\ast}};
(10,14)*{\bullet}="3"+(1,3)*{\scriptstyle x};
(1,-7)*{\bullet}="4"+(-1,3)*{\scriptstyle x}; "3";"4" **\crv{}
\POS?(.3)*{\hole}="J"; ?(0)*\dir{>}; "1";"J" **\crv{(15,-5)};
?(.15)*\dir{>}; "2";"J" **\crv{} ?(.28)*\dir{<};
(-14,10)*{}="TL"; (14,10)*{}="TR"; (14,-10)*{}="BR";
(-14,-10)*{}="BL"; (-6,20)*{}="xTL"; (22,20)*{}="xTR";
(22,0)*{}="xBR"; (-6,0)*{}="xBL";
    "TL";"TR" **\dir{-};
    "TR";"BR" **\dir{-};
    "BR";"BL" **\dir{-};
    "BL";"TL" **\dir{-};
    "xTL";"xTR" **\dir{-};
    "xTR";"xBR" **\dir{-};
    "xBR";"xBL" **\dir{.};
    "TL";"xTL" **\dir{-};
    "TR";"xTR" **\dir{-};
    "BL";"xBL" **\dir{.};
    "BR";"xBR" **\dir{-};
    "xTL";"xBL" **\dir{.};
\endxy
   };
   (0,30)*{
    \xy 0;/r.20pc/: 
(-7,10)*{\bullet}="1"+(0,3)*{x};
(-1,10)*{\bullet}="2"+(0,3)*{x^{\ast}};
(6,10)*{\bullet}="3"+(0,3)*{x}; (1,-10)*{\bullet}="4"+(0,-3)*{x};
 "1"; "2" **\crv{(-7,2) & (-1,2)}; ?(.2)*\dir{>}; ?(.9)*\dir{>};
 "3"; "4" **\crv{(9,4) & (-5,-1)}; ?(.5)*\dir{>};
(-10,10)*{}; (10,10)*{} **\dir{-}; (10,-10)*{}; (10,10)*{}
**\dir{-}; (-10,-10)*{}; (-10,10)*{} **\dir{-}; (-10,-10)*{};
(10,-10)*{} **\dir{-};
\endxy};
(40,0)*{
\xy 0;/r.20pc/:
(8.25,-1.25)*\ellipse(2,.65){-}; (6,18)*{}="b";
  (1,16)*{}="a";
 \vunder~{(1,17.5)}{(6,18)}{(3.5,15.5)}{(6,16)};
 (1,17.5)*{}; (6,16)*{} **\crv{(-5,15)& (6,13) }; \POS?(.75)*{\hole}="J1";
 "J1";(6,18)*{} **\crv{(10,14)&(11,18)};
  (6,-2)*{}="b";
  (1,-4)*{}="a";
 \vunder~{(1,-2.5)}{(6,-2)}{(3.5,-4.5)}{(6,-4)};
 (1,-2.5)*{}; (6,-4)*{} **\crv{(-5,-5)& (6,-7) }; \POS?(.75)*{\hole}="J1";
 "J1";(6,-2)*{} **\crv{(10,-6)&(11,-2)};
 (-1,16)*{}="TL";
 (9,16)*{}="TR";
 (-1,-4)*{}="BL";
 (9.25,-3.5)*{}="BR";
 (18.5,-2.5)*{}="BRR";
 (14.5,-2.5)*{}="BRR2";
 "TL";"BL" **\crv{(-2,13) & (1,6)};
 "TR";"BRR" **\crv{(9,12) & (10,6)};
 "BR";"BRR2" **\crv{(7,12) & (12.5,0)};
(-14,10)*{}="TL"; (14,10)*{}="TR"; (14,-10)*{}="BR";
(-14,-10)*{}="BL"; (-6,20)*{}="xTL"; (22,20)*{}="xTR";
(22,0)*{}="xBR"; (-6,0)*{}="xBL";
    "TL";"TR" **\dir{-};
    "TR";"BR" **\dir{-};
    "BR";"BL" **\dir{-};
    "BL";"TL" **\dir{-};
    "xTL";"xTR" **\dir{-};
    "xTR";"xBR" **\dir{-};
    "xBR";"xBL" **\dir{.};
    "TL";"xTL" **\dir{-};
    "TR";"xTR" **\dir{-};
    "BL";"xBL" **\dir{.};
    "BR";"xBR" **\dir{-};
    "xTL";"xBL" **\dir{.};
(5,22)*{\textbf{4d}}; (25,25)*{};
\endxy};
(-40,0)*{ 
\xy  0;/r.20pc/:
(-5,5)*{\bullet}+(1,3)*{x}; (5,4)*{\bullet}+(1,3)*{x^{\ast}};
(1,-7)*{\bullet}+(1,3)*{x}; (-10,10)*{}; (10,10)*{} **\dir{-};
(10,-10)*{}; (10,10)*{} **\dir{-}; (-10,-10)*{}; (-10,10)*{}
**\dir{-}; (-10,-10)*{}; (10,-10)*{} **\dir{-};
\endxy
};(0,-30)*{ 
\xy 0;/r.20pc/: 
(-4,14)*{\bullet}="1"+(,2)*{ \scriptstyle x};
(0,14)*{\bullet}="2"+(,2)*{\scriptstyle x^{\ast}};
(10,14)*{\bullet}="3"+(1,3)*{\scriptstyle x};
(1,-7)*{\bullet}="4"+(-1,3)*{\scriptstyle x}; "3";"4"
**\crv{(10,4) & (-1,-6)}; ?(.6)*\dir{>};
(-14,10)*{}="TL"; (14,10)*{}="TR"; (14,-10)*{}="BR";
(-14,-10)*{}="BL"; (-6,20)*{}="xTL"; (22,20)*{}="xTR";
(22,0)*{}="xBR"; (-6,0)*{}="xBL";
    "TL";"TR" **\dir{-}; \POS?(.35)*{\hole}="J"; \POS?(.5)*{\hole}="J1";
    "TR";"BR" **\dir{-};
    "BR";"BL" **\dir{-};
    "BL";"TL" **\dir{-};
    "xTL";"xTR" **\dir{-};
    "xTR";"xBR" **\dir{-};
    "xBR";"xBL" **\dir{.};
    "TL";"xTL" **\dir{-};
    "TR";"xTR" **\dir{-};
    "BL";"xBL" **\dir{.};
    "BR";"xBR" **\dir{-};
    "xTL";"xBL" **\dir{.};
 "1";"J" **\crv{}?(.9)*\dir{>};
 "2";"J1" **\crv{} ?(.3)*\dir{<};
 "J";"J1" **\crv{(-4,5) & (0,5)};
 (5,22)*{\textbf{4d}};
 (25,25)*{};
\endxy
};
 (0,-60)*{
   \xy 0;/r.20pc/: 
(-4,14)*{\bullet}="1"+(,2)*{ \scriptstyle x};
(0,14)*{\bullet}="2"+(,2)*{\scriptstyle x^{\ast}};
(10,14)*{\bullet}="3"+(1,3)*{\scriptstyle x};
(1,-7)*{\bullet}="4"+(-1,3)*{\scriptstyle x}; "3";"4"
**\crv{(10,4) & (-1,-6)}; ?(.6)*\dir{>};
(-14,10)*{}="TL"; (14,10)*{}="TR"; (14,-10)*{}="BR";
(-14,-10)*{}="BL"; (-6,20)*{}="xTL"; (22,20)*{}="xTR";
(22,0)*{}="xBR"; (-6,0)*{}="xBL";
    "TL";"TR" **\dir{-}; \POS?(.35)*{\hole}="J"; \POS?(.5)*{\hole}="J1";
    "TR";"BR" **\dir{-};
    "BR";"BL" **\dir{-};
    "BL";"TL" **\dir{-};
    "xTL";"xTR" **\dir{-};
    "xTR";"xBR" **\dir{-};
    "xBR";"xBL" **\dir{.};
    "TL";"xTL" **\dir{-};
    "TR";"xTR" **\dir{-};
    "BL";"xBL" **\dir{.};
    "BR";"xBR" **\dir{-};
    "xTL";"xBL" **\dir{.};
 "1";"J" **\crv{}?(.9)*\dir{>};
 "2";"J1" **\crv{} ?(.3)*\dir{<};
 "J";"J1" **\crv{(-4,5) & (0,5)};
 (5,22)*{\textbf{5d}};
 (25,25)*{};
\endxy
};
(-40,30)*{ 
\xy   
(0,0)*{\bullet}+(1,3)*{x^{\ast}}; (5,0)*{\bullet}+(1,3)*{x};
(-5,0)*{\bullet}+(1,3)*{x}; (-10,0)*{}; (10,0)*{} **\dir{-};
\endxy
};
(-40,60)*{ 
\xy (0,0)*{\bullet}+(2,3)*{x^{\ast}};
\endxy
};
(-40,-30)*{ 
  \xy 0;/r.20pc/: 
(-2.5,4)*{\bullet}+(1,2)*{\scriptstyle x}; (5,7)*{
\bullet}+(3,1)*{\scriptstyle x^{\ast}};
(6,0)*{\bullet}+(1,2)*{\scriptstyle x};
(-14,10)*{}="TL"; (14,10)*{}="TR"; (14,-10)*{}="BR";
(-14,-10)*{}="BL"; (-6,20)*{}="xTL"; (22,20)*{}="xTR";
(22,0)*{}="xBR"; (-6,0)*{}="xBL";
    "TL";"TR" **\dir{-};
    "TR";"BR" **\dir{-};
    "BR";"BL" **\dir{-};
    "BL";"TL" **\dir{-};
    "xTL";"xTR" **\dir{-};
    "xTR";"xBR" **\dir{-};
    "xBR";"xBL" **\dir{.};
    "TL";"xTL" **\dir{-};
    "TR";"xTR" **\dir{-};
    "BL";"xBL" **\dir{.};
    "BR";"xBR" **\dir{-};
    "xTL";"xBL" **\dir{.};
(25,25)*{};
\endxy
}; (-40,-60)*{
  \xy 0;/r.20pc/: 
(-2.5,4)*{\bullet}+(1,2)*{\scriptstyle x}; (5,7)*{
\bullet}+(3,1)*{\scriptstyle x^{\ast}};
(6,0)*{\bullet}+(1,2)*{\scriptstyle x};
(-14,10)*{}="TL"; (14,10)*{}="TR"; (14,-10)*{}="BR";
(-14,-10)*{}="BL"; (-6,20)*{}="xTL"; (22,20)*{}="xTR";
(22,0)*{}="xBR"; (-6,0)*{}="xBL";
    "TL";"TR" **\dir{-};
    "TR";"BR" **\dir{-};
    "BR";"BL" **\dir{-};
    "BL";"TL" **\dir{-};
    "xTL";"xTR" **\dir{-};
    "xTR";"xBR" **\dir{-};
    "xBR";"xBL" **\dir{.};
    "TL";"xTL" **\dir{-};
    "TR";"xTR" **\dir{-};
    "BL";"xBL" **\dir{.};
    "BR";"xBR" **\dir{-};
    "xTL";"xBL" **\dir{.};
(25,25)*{}; (5,22)*{\textbf{4d}};
\endxy
}; (0,60)*{
 \xy  0;/r.20pc/:  
(0,10)*{\bullet}="a"+(2.5,1)*{x};
(0,-10)*{\bullet}="b"+(2.5,1)*{x}; "a";"b" **\dir{-};
?(.45)*\dir{>};
\endxy
}; (40,60)*{
 \xy 0;/r.20pc/: 
(-10,10)*{\bullet}="TL"+(-1,3)*{x}; (10,10)*{}="2"="TR"+(1,3)*{x};
(10,-10)*{\bullet}="BR"+(1,-3)*{x};
(-10,-10)*{\bullet}="BL"+(-1,-3)*{x}; "TL";"TR" **\dir{-};
?(.5)*\dir{>}; "BL";"BR" **\dir{-}; ?(.5)*\dir{>}; "TL";"BL"
**\dir{-}; "TR";"BR" **\dir{-}; (0,3)*{}="a"; (0,-3)*{}="a'";
{\ar@{=>} "a";"a'"};
\endxy
}; (40,30)*{
 \xy 0;/r.20pc/: 
(-9,10)*{\bullet}="1"+(-1,-3)*{x^{\ast}};
(-1,10)*{\bullet}="2"+(-1,-3)*{x};
 (-9,-10)*{\bullet}="b1"+(-1,-3)*{x^{\ast}};
(-1,-10)*{\bullet}="b2"+(-1,-3)*{x};
(13,20)*{\bullet}="4"+(1,3)*{x};
(7,20)*{\bullet}="3"+(1,3)*{x^{\ast}}; (13,0)*{\bullet}="b4";
(7,0)*{\bullet}="b3"+(1,3); "1";"2" **\crv{(-4,15) & (3,15)};
?(.15)*\dir{<}; "3";"4" **\crv{(3,15) & (10,15)}; ?(.85)*\dir{>};
 (-1,-5.25)*{}="M";
"1";"b1"  **\dir{-}; "2";"b2"  **\dir{-}; "b1";"M"  **\dir{-};
?(.76)*\dir{>}; "b3";"M"  **\dir{.}; "b2";"b4"  **\dir{-};
?(.4)*\dir{<}; "4";"b4"  **\dir{-}; "3";"b3"  **\dir{.};
(0,13.25)*{}="z1"; (7,16.15)*{}="z2"; "z1";"z2" **\crv{(2,-3) &
(7,3)};
(-14,10)*{}="TL"; (14,10)*{}="TR"; (14,-10)*{}="BR";
(-14,-10)*{}="BL"; (-6,20)*{}="xTL"; (22,20)*{}="xTR";
(22,0)*{}="xBR"; (-6,0)*{}="xBL";
    "TL";"TR" **\dir{-};
    "TR";"BR" **\dir{-};
    "BR";"BL" **\dir{-};
    "BL";"TL" **\dir{-};
    "xTL";"xTR" **\dir{-};
    "xTR";"xBR" **\dir{-};
    "xBR";"xBL" **\dir{.};
    "TL";"xTL" **\dir{-};
    "TR";"xTR" **\dir{-};
    "BL";"xBL" **\dir{.};
    "BR";"xBR" **\dir{-};
    "xTL";"xBL" **\dir{.};
(25,25)*{};
\endxy
}; (40,-30)*{
 \xy 0;/r.20pc/:
(-1.5,-2)*\ellipse(3,1){-}; (3.5,-2)*\ellipse(3,1){-};
(6,18)*{}="b";
  (1,16)*{}="a";
 \vunder~{(1,17.5)}{(6,18)}{(3.5,15.5)}{(6,16)};
 (1,17.5)*{}; (6,16)*{} **\crv{(-5,15)& (6,13) }; \POS?(.75)*{\hole}="J1";
 "J1";(6,18)*{} **\crv{(10,14)&(11,18)};
 (-1,16)*{}="TL";
 (9,16)*{}="TR";
 (-6,-4)*{}="BLL";
 (0,-4)*{}="BL";
 (10,-4)*{}="BRR";
 (4,-4)*{}="BR";
 (3,5)*{}="C";
 (4.25,6.25)*{}="C2";
 (4.25,15)*{}="C3";
   "TL";"BLL" **\crv{(-2,13) & (-1,6)};
   "TR";"BRR" **\crv{(9,13) & (10,6)};
   "C";"BL" **\crv{};
   "C";"BR" **\crv{(6,10)};
    "C3";"C2" **\dir{.};
(-14,10)*{}="TL"; (14,10)*{}="TR"; (14,-10)*{}="BR";
(-14,-10)*{}="BL"; (-6,20)*{}="xTL"; (22,20)*{}="xTR";
(22,0)*{}="xBR"; (-6,0)*{}="xBL";
    "TL";"TR" **\dir{-};
    "TR";"BR" **\dir{-};
    "BR";"BL" **\dir{-};
    "BL";"TL" **\dir{-};
    "xTL";"xTR" **\dir{-};
    "xTR";"xBR" **\dir{-};
    "xBR";"xBL" **\dir{.};
    "TL";"xTL" **\dir{-};
    "TR";"xTR" **\dir{-};
    "BL";"xBL" **\dir{.};
    "BR";"xBR" **\dir{-};
    "xTL";"xBL" **\dir{.};
(5,22)*{\textbf{5d}}; (25,25)*{};
\endxy
}; (40,-60)*{
\xy 0;/r.20pc/:
(-1.5,-2)*\ellipse(3,1){-}; (3.5,-2)*\ellipse(3,1){-};
(6,18)*{}="b";
  (1,16)*{}="a";
 \vunder~{(1,17.5)}{(6,18)}{(3.5,15.5)}{(6,16)};
 (1,17.5)*{}; (6,16)*{} **\crv{(-5,15)& (6,13) }; \POS?(.75)*{\hole}="J1";
 "J1";(6,18)*{} **\crv{(10,14)&(11,18)};
 (-1,16)*{}="TL";
 (9,16)*{}="TR";
 (-6,-4)*{}="BLL";
 (0,-4)*{}="BL";
 (10,-4)*{}="BRR";
 (4,-4)*{}="BR";
 (3,5)*{}="C";
 (4.25,6.25)*{}="C2";
 (4.25,15)*{}="C3";
   "TL";"BLL" **\crv{(-2,13) & (-1,6)};
   "TR";"BRR" **\crv{(9,13) & (10,6)};
   "C";"BL" **\crv{};
   "C";"BR" **\crv{(6,10)};
    "C3";"C2" **\dir{.};
(-14,10)*{}="TL"; (14,10)*{}="TR"; (14,-10)*{}="BR";
(-14,-10)*{}="BL"; (-6,20)*{}="xTL"; (22,20)*{}="xTR";
(22,0)*{}="xBR"; (-6,0)*{}="xBL";
    "TL";"TR" **\dir{-};
    "TR";"BR" **\dir{-};
    "BR";"BL" **\dir{-};
    "BL";"TL" **\dir{-};
    "xTL";"xTR" **\dir{-};
    "xTR";"xBR" **\dir{-};
    "xBR";"xBL" **\dir{.};
    "TL";"xTL" **\dir{-};
    "TR";"xTR" **\dir{-};
    "BL";"xBL" **\dir{.};
    "BR";"xBR" **\dir{-};
    "xTL";"xBL" **\dir{.};
(5,22)*{\textbf{6d}}; (25,25)*{};
\endxy
}:
\endxy 
\]
} }

\newcommand{\feynmandiagram}{
  \xy 0;/r.22pc/:      
 (-3,10)*{}="TL"; (3,10)*{}="TR";
 (0,-2)*{}="B";
 (0,-12)*{}="BB";
    "TL";"B" **\dir{-}?(.5)*\dir{>};
    "TR";"B" **\dir{~};
    "B";"BB" **\dir{-}?(.5)*\dir{>};
 \endxy
                    \; \; \; + \; \;
\xy 0;/r.22pc/:
 (-3,10)*{}="TL"; (3,10)*{}="TR";
 (0,3)*{}="B";
 (0,-12)*{}="BB";
    "TL";"B" **\dir{-}?(.5)*\dir{>};
    "TR";"B"+(0,-.6) **\dir{~};
    "B";"BB" **\dir{-}?(.5)*\dir{>} ?(.24)*\dir{}="1" ?(.84)*\dir{}="2";
    {\ar@/^.35pc/@{~} "1"+(1,0);"2"+(1,0)};
 \endxy
                    \; \; + \; \;
\xy
 (-3,10)*{}="TL"; (4,10)*{}="TR";
 (0,-2)*{}="B";
 (0,-12)*{}="BB";
    "TL";"B" **\dir{-}?(.5)*\dir{>}?(.17)*\dir{}="1";
    "TR";"B" **\dir{~} ?(.5)*\dir{}="mid";
    "B";"BB" **\dir{-}?(.5)*\dir{>} ?(.7)*\dir{}="2";
    {\ar@/^1pc/@{~}|<<<<{\hole} "1"+(1,0);"2"+(1,0)};
 \endxy
                     \; \; + \; \cdots \; + \;
 \xy
 (-3,10)*{}="TL"; (3,10)*{}="TR";
 (0,3)*{}="B";
 (0,-12)*{}="BB";
    "TL";"B" **\dir{-}?(.35)*\dir{>} ?(.75)*\dir{>};
    "TR";"B"+(0,-.6) **\dir{~};
    "B";"BB" **\dir{-}?(.4)*\dir{>}?(.7)*\dir{>} ?(.24)*\dir{}="1" ?(.85)*\dir{}="2";
    {\ar@/^.35pc/@{~} "1"+(.9,0);"2"+(.9,0)};
     {\ar@{~} (-5.2,-2.7);(0,-8.1)};
   {\ar@{~} (-1.9,8);(-5,2.1)};
     (-5,0)*\xycircle(2,2.8){-};
     {\ar@{~} (-3,0);(0,-5)};
 \endxy
    \; \; + \; \cdots
}

\newcommand{\bigpentagon}{
\xy 0;/r.20pc/:
    (-24.73,8.03)*+{\big(k(hg) \big)f}="l";
    (0,26)*+{\big((kh)g \big)f}="t";
    (24.73,8.03)*+{(kh)(gf)}="r";
    (15.28,-21.03)*+{k\big((h(gf)\big)}="br";
    (-15.28,-21.03)*+{k\big((hg)f\big)}="bl";
     {\ar^{a} "t";"l"};
     {\ar_{a} "l";"bl"};
     {\ar_{a} "br";"bl"};
     {\ar_{a} "r";"br"};
     {\ar^{a} "t";"r"};
(-51.35,16.68)*+{
    \xy 0;/r.16pc/:
     (-10.6,-14.63)*{}="t1";
     (-17.1,5.56)*{}="t2";
     (0,18)*{}="t3";
     (17.1,5.56)*{}="t4";
     (10.6,-14.63)*{}="t5";
    {\ar@{-}^{f} "t1";"t2"};
    {\ar@{-}^{g} "t2";"t3"};
    {\ar@{-}^{h} "t3";"t4"};
    {\ar@{-}^{k} "t4";"t5"};
    {\ar@{-} "t1";"t5"};
    {\ar@{-} "t2";"t4"};
    {\ar@{-} "t2";"t5"};
   \endxy};
(0,54)*+{
    \xy 0;/r.16pc/:
     (-10.6,-14.63)*{}="t1";
     (-17.1,5.56)*{}="t2";
     (0,18)*{}="t3";
     (17.1,5.56)*{}="t4";
     (10.6,-14.63)*{}="t5";
    {\ar@{-}^{f} "t1";"t2"};
    {\ar@{-}^{g} "t2";"t3"};
    {\ar@{-}@{-}^{h} "t3";"t4"};
    {\ar@{-}^{k} "t4";"t5"};
    {\ar@{-} "t1";"t5"};{\ar@{-} "t2";"t5"};
    {\ar@{-} "t3";"t5"};
   \endxy};
(51.36,16.68)*+{
    \xy 0;/r.16pc/:
     (-10.6,-14.63)*{}="t1";
     (-17.1,5.56)*{}="t2";
     (0,18)*{}="t3";
     (17.1,5.56)*{}="t4";
     (10.6,-14.63)*{}="t5";
   {\ar@{-}^{f} "t1";"t2"};
   {\ar@{-}^{g} "t2";"t3"};
   {\ar@{-}^{h} "t3";"t4"};
   {\ar@{-}^{k} "t4";"t5"};
   {\ar@{-} "t1";"t5"};
   {\ar@{-} "t1";"t3"};
   {\ar@{-} "t3";"t5"};
  \endxy};
(31.73,-43.68)*+{
    \xy 0;/r.16pc/:
     (-10.6,-14.63)*{}="t1";
     (-17.1,5.56)*{}="t2";
     (0,18)*{}="t3";
     (17.1,5.56)*{}="t4";
     (10.6,-14.63)*{}="t5";
    {\ar@{-}^{f} "t1";"t2"};
    {\ar@{-}^{g} "t2";"t3"};
    {\ar@{-}^{h} "t3";"t4"};
    {\ar@{-}^{k} "t4";"t5"};
    {\ar@{-} "t1";"t5"};
    {\ar@{-} "t1";"t3"};
    {\ar@{-} "t1";"t4"};
  \endxy};
(-31.73,-43.68)*+{
    \xy 0;/r.16pc/:
     (-10.6,-14.63)*{}="t1";
     (-17.1,5.56)*{}="t2";
     (0,18)*{}="t3";
     (17.1,5.56)*{}="t4";
     (10.6,-14.63)*{}="t5";
   {\ar@{-}^{f} "t1";"t2"};
   {\ar@{-}^{g} "t2";"t3"};
   {\ar@{-}^{h} "t3";"t4"};
   {\ar@{-}^{k} "t4";"t5"};
   {\ar@{-} "t1";"t5"};
   {\ar@{-} "t2";"t4"};
   {\ar@{-} "t1";"t4"};
  \endxy};
\endxy
}

\newcommand{\tinytwothreemove}{
 \xy 0;/r.09pc/:
 (6.18,19)*{}="t1"; 
 (-16.18,11.74)*{}="t2";
 (-16.18,-11.74)*{}="t3";
(6.18,-19)*{}="t4";
 (20,0)*{}="t5";
   {\ar@{-}"t1";"t2"};
   {\ar@{-} "t2";"t3"};
   {\ar@{-} "t3";"t4"};
   {\ar@{-} "t4";"t5"};
   {\ar@{-} "t1";"t5"};
   {\ar@{-} "t1";"t3"}; {\ar@{-} "t3";"t5"};
   {\ar@{-}|>>>>>>{ \hole \; \hole} "t1";"t4"};
   {\ar@{-}|<<<<<<{\hole}|<<<<<<<<<<{ \hole} "t2";"t4"};
\endxy
\quad   = \quad
 \xy 0;/r.09pc/:
 (6.18,19)*{}="t1"; 
 (-16.18,11.74)*{}="t2";
 (-16.18,-11.74)*{}="t3";
(6.18,-19)*{}="t4";
 (20,0)*{}="t5";
   {\ar@{-}"t1";"t2"};
   {\ar@{-} "t2";"t3"};
   {\ar@{-} "t3";"t4"};
   {\ar@{-} "t4";"t5"};
   {\ar@{-} "t1";"t5"};
   {\ar@{-} "t1";"t3"}; {\ar@{-} "t3";"t5"};
   {\ar@{-}|<<<<<<{ \hole \; \hole}|>>>>>>{ \hole \; \hole} "t1";"t4"};
   {\ar@{-}|<<<<<<{ \hole } "t2";"t5"};
   {\ar@{-}|<<<<<<{\hole}|<<<<<<<<<<{ \hole} "t2";"t4"};
\endxy
}

\newcommand{\bigtwothreemove}{
\xy 
 (6.18,19)*{}="t1"; 
 (-16.18,11.74)*{}="t2";
 (-16.18,-11.74)*{}="t3";
(6.18,-19)*{}="t4";
 (20,0)*{}="t5";
   {\ar@{-}"t1";"t2"};
   {\ar@{-} "t2";"t3"};
   {\ar@{-}|<<<<<<<<<<<<<<<{ \hole \; \hole}|>>>>>>>>>>>>>>>{ \hole \; \hole} "t2";"t4"};
   {\ar@{-} "t3";"t4"};
   {\ar@{-} "t4";"t5"};
   {\ar@{-} "t1";"t5"};
   {\ar@{-} "t1";"t3"}; {\ar@{-} "t3";"t5"};
   {\ar@{-}|>>>>>>>>>>>>>>>{ \hole \; \hole}  "t1";"t4"};
\endxy
\qquad \qquad  = \qquad  \qquad
 \xy 
 (6.18,19)*{}="t1"; 
 (-16.18,11.74)*{}="t2";
 (-16.18,-11.74)*{}="t3";
(6.18,-19)*{}="t4";
 (20,0)*{}="t5";
   {\ar@{-}"t1";"t2"};
   {\ar@{-} "t2";"t3"};
   {\ar@{-} "t3";"t4"};
   {\ar@{-} "t4";"t5"};
   {\ar@{-} "t1";"t5"};
   {\ar@{-} "t1";"t3"}; {\ar@{-} "t3";"t5"};
   {\ar@{-}|<<<<<<<<<<<<<<<{ \hole \; \hole}|>>>>>>>>>>>>>>>{ \hole \; \hole} "t1";"t4"};
   {\ar@{-}|<<<<<<<<<<<<<<{ \hole \; \hole} "t2";"t5"};
   {\ar@{-}|<<<<<<<<<<<<<<<{ \hole \; \hole}|>>>>>>>>>>>>>>>{ \hole \; \hole} "t2";"t4"};
\endxy
}

\newcommand{\mediumtwothreemove}{
 \xy 0;/r.17pc/:
 (6.18,19)*{}="t1"; 
 (-16.18,11.74)*{}="t2";
 (-16.18,-11.74)*{}="t3";
(6.18,-19)*{}="t4";
 (20,0)*{}="t5";
   {\ar@{-}"t1";"t2"};
   {\ar@{-} "t2";"t3"};
   {\ar@{-} "t3";"t4"};
   {\ar@{-} "t4";"t5"};
   {\ar@{-} "t1";"t5"};
   {\ar@{-} "t1";"t3"}; {\ar@{-} "t3";"t5"};
   {\ar@{-}|>>>>>>>>>>>{ \hole \; \hole} "t1";"t4"};
   {\ar@{-}|<<<<<<<<<<<{ \hole \; \hole}|>>>>>>>>>>>{ \hole \; \hole} "t2";"t4"};
\endxy
\qquad   =   \qquad
 \xy 0;/r.18pc/:
 (6.18,19)*{}="t1"; 
 (-16.18,11.74)*{}="t2";
 (-16.18,-11.74)*{}="t3";
(6.18,-19)*{}="t4";
 (20,0)*{}="t5";
   {\ar@{-}"t1";"t2"};
   {\ar@{-} "t2";"t3"};
   {\ar@{-} "t3";"t4"};
   {\ar@{-} "t4";"t5"};
   {\ar@{-} "t1";"t5"};
   {\ar@{-} "t1";"t3"}; {\ar@{-} "t3";"t5"};
   {\ar@{-}|<<<<<<<<<<<{ \hole \; \hole}|>>>>>>>>>>>{ \hole \; \hole} "t1";"t4"};
   {\ar@{-}|<<<<<<<<<<<<{ \hole \; \hole} "t2";"t5"};
   {\ar@{-}|<<<<<<<<<<<{ \hole \; \hole}|>>>>>>>>>>>{ \hole \; \hole} "t2";"t4"};
\endxy
}

\newcommand{\smalltwothreemove}{
 \xy 0;/r.10pc/:
 (6.18,19)*{}="t1"; 
 (-16.18,11.74)*{}="t2";
 (-16.18,-11.74)*{}="t3";
(6.18,-19)*{}="t4";
 (20,0)*{}="t5";
   {\ar@{-}"t1";"t2"};
   {\ar@{-} "t2";"t3"};
   {\ar@{-} "t3";"t4"};
   {\ar@{-} "t4";"t5"};
   {\ar@{-} "t1";"t5"};
   {\ar@{-} "t1";"t3"}; {\ar@{-} "t3";"t5"};
   {\ar@{-}|>>>>>>>{ \hole \; \hole} "t1";"t4"};
   {\ar@{-}|<<<<<<<{  \hole}|>>>>>>{  \hole} "t2";"t4"};
\endxy
\qquad   = \qquad
 \xy 0;/r.10pc/:
 (6.18,19)*{}="t1"; 
 (-16.18,11.74)*{}="t2";
 (-16.18,-11.74)*{}="t3";
(6.18,-19)*{}="t4";
 (20,0)*{}="t5";
   {\ar@{-}"t1";"t2"};
   {\ar@{-} "t2";"t3"};
   {\ar@{-} "t3";"t4"};
   {\ar@{-} "t4";"t5"};
   {\ar@{-} "t1";"t5"};
   {\ar@{-} "t1";"t3"}; {\ar@{-} "t3";"t5"};
   {\ar@{-}|<<<<<<<{ \hole \; \hole}|>>>>>>>{ \hole \; \hole} "t1";"t4"};
   {\ar@{-}|<<<<<<<{ \hole } "t2";"t5"};
   {\ar@{-}|<<<<<<<{  \hole}|>>>>>>{  \hole} "t2";"t4"};
\endxy
}

\newcommand{\smallfouronemove}{
\xy 0;/r.15pc/:
 (-10,-5 )*{}="1";
 (8,-10)*{}="2";
 (15,0)*{}="3";
 (1,12)*{}="4";
    {\ar@{-} "1";"2" };
    {\ar@{-}"2";"3" };
    {\ar@{-} "4";"3" };
    {\ar@{-} "1";"4" };
    {\ar@{-} "4";"2" };
    {\ar@{.}|>>>>>>>>>>{\hole \hole} "1";"3"};
 \endxy
\qquad =\qquad
 \xy 0;/r.15pc/:
 (-10,-5 )*{}="1";
 (8,-10)*{}="2";
 (15,0)*{}="3";
 (1,12)*{}="4";
 (0,3)*{}="m";
    {\ar@{-} "1";"m" };
    {\ar@{-} "2";"m" };
    {\ar@{-} "3";"m" };
    {\ar@{-} "4";"m" };
    {\ar@{-} "1";"2" };
    {\ar@{-}"2";"3" };
    {\ar@{-} "4";"3" };
    {\ar@{-} "1";"4" };
    {\ar@{-} "4";"2" };
    {\ar@{.}|>>>>>>>>>>{\hole \hole} "1";"3"};
 \endxy
 }

\newcommand{\mediumfouronemove}{
 \xy
 (-10,-5 )*{}="1";
 (8,-10)*{}="2";
 (15,0)*{}="3";
 (1,12)*{}="4";
    {\ar@{-} "1";"2" };
    {\ar@{-}"2";"3" };
    {\ar@{-} "4";"3" };
    {\ar@{-} "1";"4" };
    {\ar@{-} "4";"2" };
    {\ar@{.}|>>>>>>>>>>{\hole \hole} "1";"3"};
 \endxy
\qquad =\qquad
 \xy
 (-10,-5 )*{}="1";
 (8,-10)*{}="2";
 (15,0)*{}="3";
 (1,12)*{}="4";
 (0,3)*{}="m";
    {\ar@{-} "1";"m" };
    {\ar@{-} "2";"m" };
    {\ar@{-} "3";"m" };
    {\ar@{-} "4";"m" };
    {\ar@{-} "1";"2" };
    {\ar@{-}"2";"3" };
    {\ar@{-} "4";"3" };
    {\ar@{-} "1";"4" };
    {\ar@{-} "4";"2" };
    {\ar@{.}|>>>>>>>>>>{\hole \hole} "1";"3"};
 \endxy
}

\newcommand{\stickassociator}{
\xy 0;/r.13pc/:
 (-30,0)*{
  \xy 0;/r.15pc/:
    (-8,10)*{}="TL";
    (8,10)*{}="TR";
    (-2,10)*{}="X'";
    (-5,7)*{}="XM'";
    (2,10)*{}="X";
    (5,7)*{}="XM";
    (0,2)*{}="M";
    (0,-10)*{}="B";
    "TL";"M" **\dir{-};
    "X";"XM" **\dir{-};
    "X'";"XM'" **\dir{-};
    "TR";"M" **\dir{-};
    "B";"M" **\dir{-};
    \endxy
    }="1";
 (30,0)*{
    \xy 0;/r.15pc/:
    (-8,10)*{}="TL";
    (8,10)*{}="TR";
    (-2,10)*{}="C";
    (2,10)*{}="X";
    (0,8)*{}="XM";
    (3,5)*{}="CM";
    (0,2)*{}="M";
    (0,-10)*{}="B";
    "TL";"M" **\dir{-};
    "X";"XM" **\dir{-};
    "C";"CM" **\dir{-};
    "TR";"M" **\dir{-};
    "B";"M" **\dir{-};
    \endxy
    }="5";
 (15,-30)*{
    \xy 0;/r.15pc/:
    (-8,10)*{}="TL";
    (8,10)*{}="TR";
    (2,10)*{}="C";
    (-2,10)*{}="X";
    (0,8)*{}="XM";
    (-3,5)*{}="CM";
    (0,2)*{}="M";
    (0,-10)*{}="B";
    "TL";"M" **\dir{-};
    "X";"XM" **\dir{-};
    "C";"CM" **\dir{-};
    "TR";"M" **\dir{-};
    "B";"M" **\dir{-};
    \endxy
    }="4";
 (0,20)*{
    \xy 0;/r.15pc/:
    (-8,10)*{}="TL";
    (8,10)*{}="TR";
    (2,10)*{}="C";
    (-2,10)*{}="X";
    (-5,7)*{}="XM";
    (-3,5)*{}="CM";
    (0,2)*{}="M";
    (0,-10)*{}="B";
    "TL";"M" **\dir{-};
    "X";"XM" **\dir{-};
    "C";"CM" **\dir{-};
    "TR";"M" **\dir{-};
    "B";"M" **\dir{-};
\endxy
    }="2";
 (-15,-30)*{
 \xy 0;/r.15pc/: 
    (8,10)*{}="TL";
    (-8,10)*{}="TR";
    (-2,10)*{}="C";
    (2,10)*{}="X";
    (5,7)*{}="XM";
    (3,5)*{}="CM";
    (0,2)*{}="M";
    (0,-10)*{}="B";
    "TL";"M" **\dir{-};
    "X";"XM" **\dir{-};
    "C";"CM" **\dir{-};
    "TR";"M" **\dir{-};
    "B";"M" **\dir{-};
\endxy
    }="3";
    {\ar@{=>} "2";"1"};
    {\ar@{=>} "2";"5"};
    {\ar@{=>} "1";"3"};
    {\ar@{=>} "4";"3"};
    {\ar@{=>} "5";"4"};
\endxy}

\newcommand{\LTQTtable}{
\begin{tabular}{|c|c|}
  \hline
  \textbf{2D Lattice Field Theory} & \textbf{3D Lattice Field Theory} \\
  \hline
  $\vcenter{\xy (0,0)*{\LARGE \bullet}; \endxy}$  unlabeled
&
  $\vcenter{\xy (0,0)*{\LARGE \bullet}; \endxy}$  unlabeled
\\
  $\vcenter{\xy {\ar (-6,0)*+{\bullet};
  (6,0)*+{\bullet}}; (0,5)*{};(0,-5)*{};\endxy}$  $A \in \Vect$
 & $\vcenter{\xy {\ar (-6,0)*+{\bullet};
  (6,0)*+{\bullet}}; (0,5)*{};(0,-5)*{};\endxy}$  $A \in$ 2-$\Vect$ \\
  $ \vcenter{  \xy 0;/r.16pc/:
  (0,6)*{\bullet};
 (-10,0)*{}="L";
 (10,0)*{}="R";
 (0,16)*{}="T";
 (0,6)*{}="M";
 (0,-4)*{}="B";
 (-10,12)*{}="TL";
 (10,12)*{}="TR";
    "T";"L" **\dir{.};
    "R";"T" **\dir{.};
    "L";"R" **\dir{.};
    "TL";"M" **\dir{-}?(.5)*\dir{>};
    "TR";"M" **\dir{-}?(.5)*\dir{>};
    "M";"B" **\dir{-}?(.6)*\dir{>};
 \endxy}$ $m \maps A \tensor A \to A$ &   $ \vcenter{  \xy 0;/r.16pc/:
  (0,6)*{\bullet};
 (-10,0)*{}="L";
 (10,0)*{}="R";
 (0,16)*{}="T";
 (0,6)*{}="M";
 (0,-4)*{}="B";
 (-10,12)*{}="TL";
 (10,12)*{}="TR";
    "T";"L" **\dir{.};
    "R";"T" **\dir{.};
    "L";"R" **\dir{.};
    "TL";"M" **\dir{-}?(.5)*\dir{>};
    "TR";"M" **\dir{-}?(.5)*\dir{>};
    "M";"B" **\dir{-}?(.6)*\dir{>}; (0,-12)*{};
 \endxy}$ $m \maps A \tensor A \to A$ \footnote{Care must be taken in the definition
 of this tensor product but using a basis this can be done without to much difficulty.} \\
  $\xy 0;/r.16pc/:
 (0,10)*{}="mt";
 (0,-10)*{}="mb";
 (-16,0)*{}="l";
 (16,0)*{}="r";
  (6,0)*{}="xr";
  (-6,0)*{}="xl";
  (-12,10)*{}="xlt";
  (-12,-10)*{}="xlb";
    (12,10)*{}="xrt";
  (12,-10)*{}="xrb";
  (6,0)*{\bullet};
  (-6,0)*{\bullet};
    "mt";"mb" **\dir{.};
    "mb";"l" **\dir{.};
    "mt";"l" **\dir{.};
    "mb";"r" **\dir{.};
    "mt";"r" **\dir{.};
    "xl";"xr" **\dir{-}?(.5)*\dir{<};
    "xl";"xlt" **\dir{-}?(.5)*\dir{<};
    "xl";"xlb" **\dir{-}?(.65)*\dir{>};
    "xr";"xrt" **\dir{-}?(.5)*\dir{<};
    "xr";"xrb" **\dir{-}?(.5)*\dir{<};
 \endxy
 \quad = \quad
 \xy 0;/r.16pc/:
 (0,10)*{}="mt";
 (0,-10)*{}="mb";
 (-16,0)*{}="l";
 (16,0)*{}="r";
  (0,-4)*{}="xr";
  (0,4)*{}="xl";
  (10,10)*{}="xlt";
  (-10,10)*{}="xlb";
    (10,-10)*{}="xrt";
  (-10,-10)*{}="xrb";
  (0,4)*{\bullet};
  (0,-4)*{\bullet};
    "l";"r" **\dir{.};
    "mb";"l" **\dir{.};
    "mt";"l" **\dir{.};
    "mb";"r" **\dir{.};
    "mt";"r" **\dir{.};
    "xl";"xr" **\dir{-}?(.58)*\dir{>};
    "xl";"xlt" **\dir{-}?(.5)*\dir{<};
    "xl";"xlb" **\dir{-}?(.5)*\dir{<};
    "xr";"xrt" **\dir{-}?(.5)*\dir{<};
    "xr";"xrb" **\dir{-}?(.65)*\dir{>};
 \endxy$ &  $\xy 0;/r.16pc/:
 (0,10)*{}="mt";
 (0,-10)*{}="mb";
 (-16,0)*{}="l";
 (16,0)*{}="r";
  (6,0)*{}="xr";
  (-6,0)*{}="xl";
  (-12,10)*{}="xlt";
  (-12,-10)*{}="xlb";
    (12,10)*{}="xrt";
  (12,-10)*{}="xrb";
  (6,0)*{\bullet};
  (-6,0)*{\bullet};
    "mt";"mb" **\dir{.};
    "mb";"l" **\dir{.};
    "mt";"l" **\dir{.};
    "mb";"r" **\dir{.};
    "mt";"r" **\dir{.};
    "xl";"xr" **\dir{-}?(.5)*\dir{<};
    "xl";"xlt" **\dir{-}?(.5)*\dir{<};
    "xl";"xlb" **\dir{-}?(.65)*\dir{>};
    "xr";"xrt" **\dir{-}?(.5)*\dir{<};
    "xr";"xrb" **\dir{-}?(.5)*\dir{<};
 \endxy
 \quad \xy {\ar@{=>}^{\scs \alpha} (-3,0);(3,0)}; \endxy \quad
 \xy 0;/r.16pc/:
 (0,10)*{}="mt";
 (0,-10)*{}="mb";
 (-16,0)*{}="l";
 (16,0)*{}="r";
  (0,-4)*{}="xr";
  (0,4)*{}="xl";
  (10,10)*{}="xlt";
  (-10,10)*{}="xlb";
    (10,-10)*{}="xrt";
  (-10,-10)*{}="xrb";
  (0,4)*{\bullet};
  (0,-4)*{\bullet};
    "l";"r" **\dir{.};
    "mb";"l" **\dir{.};
    "mt";"l" **\dir{.};
    "mb";"r" **\dir{.};
    "mt";"r" **\dir{.};
    "xl";"xr" **\dir{-}?(.58)*\dir{>};
    "xl";"xlt" **\dir{-}?(.5)*\dir{<};
    "xl";"xlb" **\dir{-}?(.5)*\dir{<};
    "xr";"xrt" **\dir{-}?(.5)*\dir{<};
    "xr";"xrb" **\dir{-}?(.65)*\dir{>};
 \endxy$\\
  $\vcenter{\xy 0;/r.16pc/:
  (0,6)*{\bullet};
(0,24)*{}; 
 (-10,0)*{}="L";
 (10,0)*{}="R";
 (0,16)*{}="T";
 (0,6)*{}="M";
 (0,-4)*{}="B";
 (-10,12)*{}="TL";
 (10,12)*{}="TR";
    "T";"L" **\dir{.};
    "R";"T" **\dir{.};
    "L";"R" **\dir{.};
    "TL";"M" **\dir{-};
    "TR";"M" **\dir{-};
    "M";"B" **\dir{-};
 \endxy}
\quad = \quad
 \vcenter{\xy 0;/r.16pc/:
 (0,24)*{}; 
(-10,0)*{}="L";
 (10,0)*{}="R";
 (0,16)*{}="T";
 (0,6)*{}="M";
    "L";"T" **\dir{.};
    "R";"T" **\dir{.};
    "L";"R" **\dir{.};
    "T";"M" **\dir{.};
    "R";"M" **\dir{.};
    "L";"M" **\dir{.};
 (0,-4)*{}="B";
 (-10,12)*{}="TL";
 (10,12)*{}="TR";
 (-3.5,8)*{}="tl";
 (3.5,8)*{}="tr";
 (0,2.5)*{}="b";
    "TL";"tl" **\dir{-};
    "TR";"tr" **\dir{-};
    "b";"B" **\dir{-};
    "tl";"tr" **\dir{-};
    "tr";"b" **\dir{-};
    "tl";"b" **\dir{-}?(.6)*\dir{};
 \endxy}$ &  $\vcenter{\xy 0;/r.16pc/:
  (0,6)*{\bullet};
(0,24)*{}; 
 (-10,0)*{}="L";
 (10,0)*{}="R";
 (0,16)*{}="T";
 (0,6)*{}="M";
 (0,-4)*{}="B";
 (-10,12)*{}="TL";
 (10,12)*{}="TR";
    "T";"L" **\dir{.};
    "R";"T" **\dir{.};
    "L";"R" **\dir{.};
    "TL";"M" **\dir{-};
    "TR";"M" **\dir{-};
    "M";"B" **\dir{-};
 \endxy}
\quad \xy {\ar@{=>}^{} (-3,1);(3,1)};
          {\ar@{=>}^{} (3,-4);(-3,-4)};\endxy \quad
 \vcenter{\xy 0;/r.16pc/:
 (0,24)*{}; 
(-10,0)*{}="L";
 (10,0)*{}="R";
 (0,16)*{}="T";
 (0,6)*{}="M";
    "L";"T" **\dir{.};
    "R";"T" **\dir{.};
    "L";"R" **\dir{.};
    "T";"M" **\dir{.};
    "R";"M" **\dir{.};
    "L";"M" **\dir{.};
 (0,-4)*{}="B";
 (-10,12)*{}="TL";
 (10,12)*{}="TR";
 (-3.5,8)*{}="tl";
 (3.5,8)*{}="tr";
 (0,2.5)*{}="b";
    "TL";"tl" **\dir{-};
    "TR";"tr" **\dir{-};
    "b";"B" **\dir{-};
    "tl";"tr" **\dir{-};
    "tr";"b" **\dir{-};
    "tl";"b" **\dir{-}?(.6)*\dir{};
 \endxy}$ \\
   & $\stickassociator$ \\
   &  \\
  \hline
\end{tabular}
}

\newcommand{\LTQTtableII}{
\begin{tabular}{|c|c|}
  \hline
  \textbf{2D Lattice Field Theory} & \textbf{3D Lattice Field Theory} \\
  \hline
  $\vcenter{\xy (0,0)*{\LARGE \bullet}; \endxy}$  unlabeled
&
  $\vcenter{\xy (0,0)*{\LARGE \bullet}; \endxy}$  unlabeled
\\
  $\vcenter{\xy {\ar (-6,0)*+{\bullet};
  (6,0)*+{\bullet}}; (0,5)*{};(0,-5)*{};\endxy}$  $A \in \Vect$
 & $\vcenter{\xy {\ar (-6,0)*+{\bullet};
  (6,0)*+{\bullet}}; (0,5)*{};(0,-5)*{};\endxy}$  $A \in$ 2-$\Vect$ \\
  $ \vcenter{  \xy 0;/r.16pc/:
  (0,6)*{\bullet};
 (-10,0)*{}="L";
 (10,0)*{}="R";
 (0,16)*{}="T";
 (0,6)*{}="M";
 (0,-4)*{}="B";
 (-10,12)*{}="TL";
 (10,12)*{}="TR";
    "T";"L" **\dir{.};
    "R";"T" **\dir{.};
    "L";"R" **\dir{.};
    "TL";"M" **\dir{-}?(.5)*\dir{>};
    "TR";"M" **\dir{-}?(.5)*\dir{>};
    "M";"B" **\dir{-}?(.6)*\dir{>};
 \endxy}$ $m \maps A \tensor A \to A$ &   $ \vcenter{  \xy 0;/r.16pc/:
  (0,6)*{\bullet};
 (-10,0)*{}="L";
 (10,0)*{}="R";
 (0,16)*{}="T";
 (0,6)*{}="M";
 (0,-4)*{}="B";
 (-10,12)*{}="TL";
 (10,12)*{}="TR";
    "T";"L" **\dir{.};
    "R";"T" **\dir{.};
    "L";"R" **\dir{.};
    "TL";"M" **\dir{-}?(.5)*\dir{>};
    "TR";"M" **\dir{-}?(.5)*\dir{>};
    "M";"B" **\dir{-}?(.6)*\dir{>}; (0,-12)*{};
 \endxy}$ $m \maps A \tensor A \to A$ \footnote{Care must be taken in the definition
 of this tensor product but using a basis this can be done without to much difficulty.} \\
  $\xy 0;/r.16pc/:
 (0,10)*{}="mt";
 (0,-10)*{}="mb";
 (-16,0)*{}="l";
 (16,0)*{}="r";
  (6,0)*{}="xr";
  (-6,0)*{}="xl";
  (-12,10)*{}="xlt";
  (-12,-10)*{}="xlb";
    (12,10)*{}="xrt";
  (12,-10)*{}="xrb";
  (6,0)*{\bullet};
  (-6,0)*{\bullet};
    "mt";"mb" **\dir{.};
    "mb";"l" **\dir{.};
    "mt";"l" **\dir{.};
    "mb";"r" **\dir{.};
    "mt";"r" **\dir{.};
    "xl";"xr" **\dir{-}?(.5)*\dir{<};
    "xl";"xlt" **\dir{-}?(.5)*\dir{<};
    "xl";"xlb" **\dir{-}?(.65)*\dir{>};
    "xr";"xrt" **\dir{-}?(.5)*\dir{<};
    "xr";"xrb" **\dir{-}?(.5)*\dir{<};
 \endxy
 \quad = \quad
 \xy 0;/r.16pc/:
 (0,10)*{}="mt";
 (0,-10)*{}="mb";
 (-16,0)*{}="l";
 (16,0)*{}="r";
  (0,-4)*{}="xr";
  (0,4)*{}="xl";
  (10,10)*{}="xlt";
  (-10,10)*{}="xlb";
    (10,-10)*{}="xrt";
  (-10,-10)*{}="xrb";
  (0,4)*{\bullet};
  (0,-4)*{\bullet};
    "l";"r" **\dir{.};
    "mb";"l" **\dir{.};
    "mt";"l" **\dir{.};
    "mb";"r" **\dir{.};
    "mt";"r" **\dir{.};
    "xl";"xr" **\dir{-}?(.58)*\dir{>};
    "xl";"xlt" **\dir{-}?(.5)*\dir{<};
    "xl";"xlb" **\dir{-}?(.5)*\dir{<};
    "xr";"xrt" **\dir{-}?(.5)*\dir{<};
    "xr";"xrb" **\dir{-}?(.65)*\dir{>};
 \endxy$ &  $\xy 0;/r.16pc/:
 (0,10)*{}="mt";
 (0,-10)*{}="mb";
 (-16,0)*{}="l";
 (16,0)*{}="r";
  (6,0)*{}="xr";
  (-6,0)*{}="xl";
  (-12,10)*{}="xlt";
  (-12,-10)*{}="xlb";
    (12,10)*{}="xrt";
  (12,-10)*{}="xrb";
  (6,0)*{\bullet};
  (-6,0)*{\bullet};
    "mt";"mb" **\dir{.};
    "mb";"l" **\dir{.};
    "mt";"l" **\dir{.};
    "mb";"r" **\dir{.};
    "mt";"r" **\dir{.};
    "xl";"xr" **\dir{-}?(.5)*\dir{<};
    "xl";"xlt" **\dir{-}?(.5)*\dir{<};
    "xl";"xlb" **\dir{-}?(.65)*\dir{>};
    "xr";"xrt" **\dir{-}?(.5)*\dir{<};
    "xr";"xrb" **\dir{-}?(.5)*\dir{<};
 \endxy
 \quad \xy {\ar@{=>}^{\scs \alpha} (-3,0);(3,0)}; \endxy \quad
 \xy 0;/r.16pc/:
 (0,10)*{}="mt";
 (0,-10)*{}="mb";
 (-16,0)*{}="l";
 (16,0)*{}="r";
  (0,-4)*{}="xr";
  (0,4)*{}="xl";
  (10,10)*{}="xlt";
  (-10,10)*{}="xlb";
    (10,-10)*{}="xrt";
  (-10,-10)*{}="xrb";
  (0,4)*{\bullet};
  (0,-4)*{\bullet};
    "l";"r" **\dir{.};
    "mb";"l" **\dir{.};
    "mt";"l" **\dir{.};
    "mb";"r" **\dir{.};
    "mt";"r" **\dir{.};
    "xl";"xr" **\dir{-}?(.58)*\dir{>};
    "xl";"xlt" **\dir{-}?(.5)*\dir{<};
    "xl";"xlb" **\dir{-}?(.5)*\dir{<};
    "xr";"xrt" **\dir{-}?(.5)*\dir{<};
    "xr";"xrb" **\dir{-}?(.65)*\dir{>};
 \endxy$\\
  $\vcenter{\xy 0;/r.16pc/:
  (0,6)*{\bullet};
(0,24)*{}; 
 (-10,0)*{}="L";
 (10,0)*{}="R";
 (0,16)*{}="T";
 (0,6)*{}="M";
 (0,-4)*{}="B";
 (-10,12)*{}="TL";
 (10,12)*{}="TR";
    "T";"L" **\dir{.};
    "R";"T" **\dir{.};
    "L";"R" **\dir{.};
    "TL";"M" **\dir{-};
    "TR";"M" **\dir{-};
    "M";"B" **\dir{-};
 \endxy}
\quad = \quad
 \vcenter{\xy 0;/r.16pc/:
 (0,24)*{}; 
(-10,0)*{}="L";
 (10,0)*{}="R";
 (0,16)*{}="T";
 (0,6)*{}="M";
    "L";"T" **\dir{.};
    "R";"T" **\dir{.};
    "L";"R" **\dir{.};
    "T";"M" **\dir{.};
    "R";"M" **\dir{.};
    "L";"M" **\dir{.};
 (0,-4)*{}="B";
 (-10,12)*{}="TL";
 (10,12)*{}="TR";
 (-3.5,8)*{}="tl";
 (3.5,8)*{}="tr";
 (0,2.5)*{}="b";
    "TL";"tl" **\dir{-};
    "TR";"tr" **\dir{-};
    "b";"B" **\dir{-};
    "tl";"tr" **\dir{-};
    "tr";"b" **\dir{-};
    "tl";"b" **\dir{-}?(.6)*\dir{};
 \endxy}$ &  $\vcenter{\xy 0;/r.16pc/:
  (0,6)*{\bullet};
(0,24)*{}; 
 (-10,0)*{}="L";
 (10,0)*{}="R";
 (0,16)*{}="T";
 (0,6)*{}="M";
 (0,-4)*{}="B";
 (-10,12)*{}="TL";
 (10,12)*{}="TR";
    "T";"L" **\dir{.};
    "R";"T" **\dir{.};
    "L";"R" **\dir{.};
    "TL";"M" **\dir{-};
    "TR";"M" **\dir{-};
    "M";"B" **\dir{-};
 \endxy}
\quad \xy {\ar@{=>}^{} (-3,1);(3,1)};
          {\ar@{=>}^{} (3,-4);(-3,-4)};\endxy \quad
 \vcenter{\xy 0;/r.16pc/:
 (0,24)*{}; 
(-10,0)*{}="L";
 (10,0)*{}="R";
 (0,16)*{}="T";
 (0,6)*{}="M";
    "L";"T" **\dir{.};
    "R";"T" **\dir{.};
    "L";"R" **\dir{.};
    "T";"M" **\dir{.};
    "R";"M" **\dir{.};
    "L";"M" **\dir{.};
 (0,-4)*{}="B";
 (-10,12)*{}="TL";
 (10,12)*{}="TR";
 (-3.5,8)*{}="tl";
 (3.5,8)*{}="tr";
 (0,2.5)*{}="b";
    "TL";"tl" **\dir{-};
    "TR";"tr" **\dir{-};
    "b";"B" **\dir{-};
    "tl";"tr" **\dir{-};
    "tr";"b" **\dir{-};
    "tl";"b" **\dir{-}?(.6)*\dir{};
 \endxy}$ \\
   $\xy (0,15)*{};(0,-15)*{}; \endxy$ & $\smalltwothreemove$ \\
   $\xy (0,-10)*{}; \endxy$& $\smallfouronemove$ \\
  \hline
\end{tabular}
}

\newcommand{\bubblenattrans}{
  \vcenter{
 \xy 
    (-6,4)*{};(6,4)*{};
      **\crv{(5,-10) & (-5,-10)};
\endxy}
\qquad  \xy {\ar@{=>} (-3,0);(3,0)} \endxy  \qquad
 \vcenter{\xy 0;/r.14pc/:
  (0,6)*{};
 (0,6)*{}="M";
 (-6,-4)*{}="B";
 (-9.5,-10)*{}="b";
 (-1.5,-16)*{}="b'";
 (-4,14)*{}="TL";
 (4,12)*{}="TR";
 (12,6)*{}="tr";
    "TL";"M" **\dir{-} ;
    "TR";"M" **\dir{-};
    "M";"B" **\dir{-};
(-16,14)*{}="a";
    "a";"B" **\dir{-};
   "b";"B" **\dir{-};
   "b'";"tr" **\dir{-} ;
  "TR";"tr" **\crv{(8,18) & (16,12)};
  "b";"b'" **\crv{ (-14,-16)&(-6,-22) };
 \endxy}
 }

\newcommand{\dualbubblenattrans}{
  \vcenter{
 \xy 
    (-6,4)*{};(6,4)*{};
      **\crv{(5,-10) & (-5,-10)};
\endxy}
\qquad  \xy {\ar@{=>} (3,0);(-3,0)} \endxy  \qquad
 \vcenter{\xy 0;/r.14pc/:
  (0,6)*{};
 (0,6)*{}="M";
 (-6,-4)*{}="B";
 (-9.5,-10)*{}="b";
 (-1.5,-16)*{}="b'";
 (-4,14)*{}="TL";
 (4,12)*{}="TR";
 (12,6)*{}="tr";
    "TL";"M" **\dir{-} ;
    "TR";"M" **\dir{-};
    "M";"B" **\dir{-};
(-16,14)*{}="a";
    "a";"B" **\dir{-};
   "b";"B" **\dir{-};
   "b'";"tr" **\dir{-} ;
  "TR";"tr" **\crv{(8,18) & (16,12)};
  "b";"b'" **\crv{ (-14,-16)&(-6,-22) };
 \endxy}
 }

\newcommand{\YYYQ}{
\psset{unit=0.5cm}
\begin{pspicture}[.6](4.5,5)
 \psline[linestyle=dashed](0,3)(4,3.5)
 \psline(0,3)(3,2)
 \psline(4,3.5)(3,2)
 \psline(0,3)(2,5)
 \psline(2,5)(3,2)
 \psline(2,5)(4,3.5)
\end{pspicture}
\quad \xy {\ar@{<->} (0,0);(6,0)}; \endxy \quad
\begin{pspicture}[.6](4.5,5)
\pscustom[fillstyle=gradient,
    gradbegin=white, gradend=lightgray,
    gradmidpoint=.5,gradangle=140]{
  \psbezier(0,3)(1,3.8)(3,4.3)(4,3.5)
   \psline(3,2)
   \psline(0,3)
  }
  \pscustom[fillstyle=gradient,
    gradbegin=white, gradend=gray,
    gradmidpoint=.5,gradangle=140]{
  \psbezier(0,3)(1.3,.8)(3.7,.8)(4,3.5)
   \psline(3,2)
   \psline(0,3)
  }
 \psline[linestyle=dashed](0,3)(4,3.5)
 \psline(0,3)(2,5)
 \psline(2,5)(3,2)
 \psline(2,5)(4,3.5)
 \psline[linestyle=dashed](2.25,2.75)(0,3)
 \psline[linestyle=dashed](2.25,2.75)(4,3.5)
 \psline[linestyle=dashed](2.25,2.75)(3,2)
 \psdots(2.25,2.75)
\end{pspicture}
\quad \xy {\ar@{<->} (0,0);(6,0)}; \endxy \quad
\begin{pspicture}[.6](4.5,5)
 \psline[linestyle=dashed](0,3)(4,3.5)
 \psline(0,3)(3,2)
 \psline(4,3.5)(3,2)
 \psline(0,3)(2,5)
 \psline(2,5)(3,2)
 \psline(2,5)(4,3.5)
  \psline(1.9,3.6)(0,3)
  \psline(1.9,3.6)(2,5)
 \psline(1.9,3.6)(4,3.5)
 \psline(1.9,3.6)(3,2)
 \psdots(1.9,3.6)
\end{pspicture}
}

\newcommand{\DIMofH}{
 \vcenter{\xy
    (-5,8)*{}="x1";
    (5,8)*{}="x2";
    (-5,5)*{}="m1";
    (5,5)*{}="m2";
    (-5,-5)*{}="k1";
    (5,-5)*{}="k2";
    (-5,-8)*{}="y1";
    (5,-8)*{}="y2";
 \vtwist~{"m1"}{"m2"}{"k1"}{"k2"};
 "x1";"m1" **\dir{-} ?(.25)*\dir{>};
 "x2";"m2" **\dir{-} ?(0)*\dir{<};
 "k1";"y1" **\dir{-} ?(0)*\dir{<};
 "k2";"y2" **\dir{-} ?(.27)*\dir{>};
    "x1";"x2" **\crv{(-5,16) & (5,16)};
    "y1";"y2" **\crv{(-5,-14) & (5,-14)};
        (-7.5,8)*{H};
\endxy}
 \qquad \qquad
 \vcenter{ \xy
    (0,14)*+{1}="1";
    (0,5)*+{e^i \tensor e_i}="2";
    (0,-5)*+{e_i \tensor e^i}="3";
    (0,-14)*+{\delta^i_i = \dim(H)}="4";
        {\ar@{|->} "1";"2"};
        {\ar@{|->} "2";"3"};
        {\ar@{|->} "3";"4"};
\endxy }
}

\newcommand{\spinNET}{
\xy
    (-4,0)*{}="b";
    (4,8)*{}="t";
    (9,-6)*{}="v";
    (14,6)*{}="u";
    (-5,8)*{}="e";
    (-20,1)*{}="q";
    (-10,-10)*{}="c";
    (10,-15)*{}="d";
    (22,4)*{}="z";
    (24,-9)*{}="w";
        {\ar@/^1pc/@{-} "q";"e"};
        {\ar@/_.1pc/@{-} "q";"b"};
        {\ar@{-} "b";"e"};
        {\ar@/^.2pc/@{-} "e";"t"};
        {\ar@/_.5pc/@{-} "u";"t"};
        {\ar@/^.5pc/@{-} "u";"t"};
        {\ar@/_.2pc/@{-} "b";"v"};
        {\ar@/_.2pc/@{-} "z";"u"};
        {\ar@/_1pc/@{-} "q";"c"};
        {\ar@/_.8pc/@{-} "c";"d"};
        {\ar@/^.6pc/@{-} "c";"d"};
        {\ar@/^.3pc/@{-} "z";"v"};
        {\ar@/^.2pc/@{-} "z";"w"};
        {\ar@/^.4pc/@{-} "w";"d"};
        {\ar@{-} "v";"w"};
 \endxy
 }

\newcommand{\spinNETII}{
\def\objectstyle{\scriptscriptstyle}
\xy
    (-4,0)*{\bullet}="b";
    (4,8)*{\bullet}="t";
    (9,-6)*{\bullet}="v";
    (14,6)*{\bullet}="u";
    (-5,8)*{\bullet}="e";
    (-20,1)*{\bullet}="q";
    (-10,-10)*{\bullet}="c";
    (10,-15)*{\bullet}="d";
    (22,4)*{\bullet}="z";
    (24,-9)*{\bullet}="w";
        {\ar@/^1pc/@{-}^1 "q";"e"};
        {\ar@/_.1pc/@{-}^1 "q";"b"};
        {\ar@{-}^1 "b";"e"};
        {\ar@/^.2pc/@{-}_1 "e";"t"};
        {\ar@/_.5pc/@{-}_1 "u";"t"};
        {\ar@/^.5pc/@{-}^1 "u";"t"};
        {\ar@/_.2pc/@{-}^1 "b";"v"};
        {\ar@/_.2pc/@{-}^1 "z";"u"};
        {\ar@/_1pc/@{-}_1 "q";"c"};
        {\ar@/_.8pc/@{-}_1 "c";"d"};
        {\ar@/^.6pc/@{-}^1 "c";"d"};
        {\ar@/^.3pc/@{-}^1 "z";"v"};
        {\ar@/^.2pc/@{-}^1 "z";"w"};
        {\ar@/^.4pc/@{-}^1 "w";"d"};
        {\ar@{-}_1 "v";"w"};
 \endxy
 }

\newcommand{\PLUGzigzag}{ \xy
    (-10,-12)*{}="1E";
    (-10,0)*{}="1";
    (0,0)*{}="2";
    (10,0)*{}="3"+(3,0)*{H};
    (10,12)*{}="3B";
 "2";"1" **\crv{(0,10)& (-10,10)}
     ?(.03)*\dir{>}  ?(1)*\dir{>};
 "3";"2" **\crv{(10,-10)& (0,-10)}
     ?(.03)*\dir{>}  ;
 "1";"1E" **\dir{-};
 "3B";"3" **\dir{-};
     (-5,8.5)*{\scriptstyle e_H};
     (5,-9)*{\scriptstyle i_H};
\endxy
 \qquad \qquad
\xy
    (0,11)*+{\psi}="1";
    (0,0)*+{e_i \tensor e^i \tensor \psi}="2";
    (0,-11)*+{e_i \tensor \psi^i = \psi}="3";
        {\ar@{|->} "1";"2"};
        {\ar@{|->} "2";"3"};
\endxy
}

\newcommand{\opetopicCHART}{
\begin{center}\makebox[0pt]{
\begin{tabular}{|c|c|c|c|c|}
  \hline
   \textbf{Objects} & \textbf{Morphisms} & \textbf{2-morphisms} & \textbf{3-morphisms} & $\cdots$ \\
  \hline \hline
   $\bullet$ & $\xy
  (-6,0)*+{\bullet}="1";
  (0,8)*+{}; 
  (0,-8)*+{}; 
  (6,0)*+{\bullet}="2";
  {\ar "1";"2"};
 \endxy$
    & $
  \def\objectstyle{\scriptstyle}
\xy 
 (-10,-5)*+{\bullet}="x";
 (10,-5)*+{\bullet}="y";
 (-10,1)*{\bullet}="1";
 (-7,6)*{\bullet}="2";
 (0,6)*+{\dots}="3";
 (7,6)*{\bullet}="4";
 (10,1)*{\bullet}="5";
  {\ar "x"; "y"};
  {\ar "x"; "1"};
  {\ar "1"; "2"};
  {\ar "2"; "3"};
  {\ar "3"; "4"};
  {\ar "4"; "5"};
  {\ar "5"; "y"};
  {\ar@{=>} (0,3); (0,-3)};
\endxy$
 & $\def\objectstyle{\scriptstyle} \xy 0;/r.14pc/:
 (-10,-5)*+{\bullet}="1";
 (10,-5)*+{\bullet}="2";
 (-6,6)*+{\bullet}="3";
 (6,6)*+{\bullet}="3'";
  {\ar "1"; "2"};
  {\ar "3'"; "2"};
  {\ar "1"; "3"};
  {\ar "3"; "3'"};
  {\ar "1"; "3'"};
\endxy
\;\; \xy {\ar@3{->}^{} (-2,0);(2,0)}; \endxy \; \xy 0;/r.14pc/:
 (-10,-5)*+{\bullet}="1";
 (10,-5)*+{\bullet}="2";
 (-6,6)*+{\bullet}="3";
 (6,6)*+{\bullet}="3'";
  {\ar "1"; "2"};
  {\ar "3'"; "2"};
  {\ar "1"; "3"};
  {\ar "3"; "3'"};
\endxy$ & Opetopes  \\
  \hline
\end{tabular}
}\end{center} }

 \newcommand{\BUBimpliesTOi}{
  \def\objectstyle{\scriptstyle}
\xy
 (-10,-8)*{\bullet}="L";
 (10,-8)*{\bullet}="R";
 (0,8)*{\bullet}="T";
    "L";"T" **\dir{-};
    "R";"T" **\dir{-};
    "L";"R" **\dir{-};
 \endxy
\quad \xy {\ar^{\txt\bf{2-2}} (-5,0); (5,0)}; \endxy \quad
 \xy
 (-10,-8)*{\bullet}="L";
 (10,-8)*{\bullet}="R";
 (0,8)*{\bullet}="T";
 (0,-8)*{\bullet}="M";
    "L";"T" **\dir{-};
    "R";"T" **\dir{-};
    {\ar@/_1.5pc/@{-}"L";"R"};
    {\ar@{-}@/^1.5pc/ "L";"R"};
    "R";"M" **\dir{-};
    "L";"M" **\dir{-};
 \endxy
\quad \xy {\ar^{\txt\bf{Bubble}} (-8,0); (8,0)}; \endxy \;\;
  \xy
 (-10,-8)*{\bullet}="L";
 (10,-8)*{\bullet}="R";
 (0,8)*{\bullet}="T";
 (0,-8)*{\bullet}="M";
    "L";"T" **\dir{-};
    "R";"T" **\dir{-};
    {\ar@/_1.5pc/@{-}"L";"R"};
    "T";"M" **\dir{-};
    "R";"M" **\dir{-};
    "L";"M" **\dir{-};
 \endxy
}

 \newcommand{\TOimpliesBUBi}{
 \def\objectstyle{\scriptstyle}
  \xy
 (-10,-8)*{\bullet}="L";
 (10,-8)*{\bullet}="R";
 (0,8)*{\bullet}="T";
    "L";"T" **\dir{-};
    "R";"T" **\dir{-};
    "L";"R" **\dir{-};
 \endxy
\quad \xy {\ar^{\txt\bf{3-1}} (-8,0); (8,0)}; \endxy \;\;
  \xy
 (-10,-8)*{\bullet}="L";
 (10,-8)*{\bullet}="R";
 (0,8)*{\bullet}="T";
 (0,-8)*{\bullet}="M";
    "L";"T" **\dir{-};
    "R";"T" **\dir{-};
    {\ar@/_1.5pc/@{-}"L";"R"};
    "T";"M" **\dir{-};
    "R";"M" **\dir{-};
    "L";"M" **\dir{-};
 \endxy
 \quad \xy {\ar^{\txt\bf{2-2}} (-5,0); (5,0)}; \endxy \quad
 \xy
 (-10,-8)*{\bullet}="L";
 (10,-8)*{\bullet}="R";
 (0,8)*{\bullet}="T";
 (0,-8)*{\bullet}="M";
    "L";"T" **\dir{-};
    "R";"T" **\dir{-};
    {\ar@/_1.5pc/@{-}"L";"R"};
    {\ar@{-}@/^1.5pc/ "L";"R"};
    "R";"M" **\dir{-};
    "L";"M" **\dir{-};
 \endxy
}

\newcommand{\TWOTHREEi}{
 \xy
 (0,30)*+{
    \xy 0;/r.10pc/:
    (6.18,19)*{}="t1"; 
    (-16.18,11.74)*{}="t2";
    (-16.18,-11.74)*{}="t3";
    (6.18,-19)*{}="t4";
    (20,0)*{}="t5";
        {\ar@{-}"t1";"t2"};
        {\ar@{-} "t2";"t3"};
        {\ar@{-} "t3";"t4"};
        {\ar@{-} "t4";"t5"};
        {\ar@{-} "t1";"t5"};
        {\ar@{-} "t1";"t3"}; {\ar@{-} "t3";"t5"};
   \endxy}="t";
 (-30,15)*+{
    \xy 0;/r.10pc/:
    (6.18,19)*{}="t1"; 
    (-16.18,11.74)*{}="t2";
    (-16.18,-11.74)*{}="t3";
    (6.18,-19)*{}="t4";
    (20,0)*{}="t5";
        {\ar@{-}"t1";"t2"};
        {\ar@{-} "t2";"t3"};
        {\ar@{-} "t3";"t4"};
        {\ar@{-} "t4";"t5"};
        {\ar@{-} "t1";"t5"};
        {\ar@{-} "t1";"t3"}; {\ar@{-} "t3";"t5"};
        {\ar@{-}|<<<<<<<{ \hole } "t2";"t5"};
   \endxy}="l1";
 (-30,-15)*+{
    \xy 0;/r.10pc/:
    (6.18,19)*{}="t1"; 
    (-16.18,11.74)*{}="t2";
    (-16.18,-11.74)*{}="t3";
    (6.18,-19)*{}="t4";
    (20,0)*{}="t5";
        {\ar@{-}"t1";"t2"};
        {\ar@{-} "t2";"t3"};
        {\ar@{-} "t3";"t4"};
        {\ar@{-} "t4";"t5"};
        {\ar@{-} "t1";"t5"};
        {\ar@{-} "t1";"t3"}; {\ar@{-} "t3";"t5"};
        {\ar@{-}|<<<<<<<{ \hole } "t2";"t5"};
        {\ar@{-}|<<<<<<<{  \hole}|>>>>>>{  \hole} "t2";"t4"};
   \endxy}="l2";
 (0,-30)*+{
    \xy 0;/r.10pc/:
    (6.18,19)*{}="t1"; 
    (-16.18,11.74)*{}="t2";
    (-16.18,-11.74)*{}="t3";
    (6.18,-19)*{}="t4";
    (20,0)*{}="t5";
        {\ar@{-}"t1";"t2"};
        {\ar@{-} "t2";"t3"};
        {\ar@{-} "t3";"t4"};
        {\ar@{-} "t4";"t5"};
        {\ar@{-} "t1";"t5"};
        {\ar@{-} "t1";"t3"}; {\ar@{-} "t3";"t5"};
        {\ar@{-}|<<<<<<<{ \hole \; \hole}|>>>>>>>{ \hole \; \hole} "t1";"t4"};
        {\ar@{-}|<<<<<<<{ \hole } "t2";"t5"};
        {\ar@{-}|<<<<<<<{  \hole}|>>>>>>{  \hole} "t2";"t4"};
   \endxy}="b";
 (30,15)*+{
    \xy 0;/r.10pc/:
    (6.18,19)*{}="t1"; 
    (-16.18,11.74)*{}="t2";
    (-16.18,-11.74)*{}="t3";
    (6.18,-19)*{}="t4";
    (20,0)*{}="t5";
        {\ar@{-}"t1";"t2"};
        {\ar@{-} "t2";"t3"};
        {\ar@{-} "t3";"t4"};
        {\ar@{-} "t4";"t5"};
        {\ar@{-} "t1";"t5"};
        {\ar@{-} "t1";"t3"}; {\ar@{-} "t3";"t5"};
        {\ar@{-}|<<<<<<<{ \hole } "t2";"t5"};
   \endxy}="r1";
 (30,-15)*+{
    \xy 0;/r.10pc/:
    (6.18,19)*{}="t1"; 
    (-16.18,11.74)*{}="t2";
    (-16.18,-11.74)*{}="t3";
    (6.18,-19)*{}="t4";
    (20,0)*{}="t5";
        {\ar@{-}"t1";"t2"};
        {\ar@{-} "t2";"t3"};
        {\ar@{-} "t3";"t4"};
        {\ar@{-} "t4";"t5"};
        {\ar@{-} "t1";"t5"};
        {\ar@{-} "t1";"t3"}; {\ar@{-} "t3";"t5"};
        {\ar@{-}|<<<<<<<{ \hole } "t2";"t5"};
        {\ar@{-}|<<<<<<<{  \hole}|>>>>>>{  \hole} "t2";"t4"};
   \endxy}="r2";
    {\ar@3{->} "t";"r1" };
    {\ar@3{->} "t";"l1" };
    {\ar@3{->} "l1";"l2" };
    {\ar@3{->} "r1";"r2" };
    {\ar@3{->} "r2";"b" };
    {\ar@3{->} "l2";"b" };
 \endxy
 }

 \newcommand{\CATchartI}{
 \begin{center}  \makebox[0pt]{
\begin{tabular}{|p{2.2in}|p{2.4in}|}
  \hline
  \textbf{Set-based mathematics} & \textbf{Category-based mathematics} \\
  \hline \hline
  \textbf{Sets}
 $ \xy
 (0,0)*{\includegraphics{blob30.eps}};
 (0,0)*{\scs \bullet};
 (2,4)*{\scs \bullet};
 (-5,-4)*{\scs \bullet};
 (6,-1)*{\scs \bullet};
 \endxy $ &  \textbf{Categories}
 $ \xy
 (0,0)*{\includegraphics{blob30.eps}};
 (2,0)*{\scs \bullet}="1";
 (-3,6)*{\scs \bullet}="2";
 (-5,-4)*{\scs \bullet}="3";
 (2,-7)*{\scs \bullet}="4";
  "1";"2" **\crv{(-2,6)}?(.25)*\dir{<};
  "1";"3" **\crv{(-2,-4)}?(.25)*\dir{<};
  "3";"2" **\crv{(-6,0)}?(.5)*\dir{<};
  "1";"4" **\crv{}?(.25)*\dir{<};
 \endxy $\\
 \textbf{Functions} \quad$ \xy
 (0,0)*{\includegraphics{blob30.eps}};
 (30,0)*{\includegraphics{blobII30.eps}};
 (0,-5)*{\scs \bullet}="1";
 (2,4)*{\scs \bullet}="2";
 (6,-1)*{\scs \bullet}="3";
 (33,-2)*{\scs \bullet}="1'";
 (27,-4)*{\scs \bullet}="2'";
 (28,1)*{\scs \bullet}="3'";
 (31,5)*{\scs \bullet}="4'";
 {\ar@{>} "1";"1'"};
 {\ar@{>} "2";"4'"};
 {\ar@{>} "3";"3'"};
 \endxy$  & \textbf{Functors} \quad $ \xy
 (0,0)*{\includegraphics{functor30.eps}};
 (-14,2)*{\scs \bullet}="2";
 (-8,-3)*{\scs \bullet}="3";
 (21,-2)*{\scs \bullet}="1'";
 (17,-4)*{\scs \bullet}="3'";
 (18,5)*{\scs \bullet}="4'";
  "2";"3" **\crv{}?(.5)*\dir{>};
  "3'";"4'" **\crv{}?(.5)*\dir{>};
  "1'";"4'" **\crv{(22,2)}?(.55)*\dir{>};
 {\ar@{>} "2";"4'"};
 {\ar@{>} "3";"3'"};
 \endxy$ \\
   &  \textbf{Natural Transformations}
 $ \xy
 (0,0)*{\includegraphics{natural40.eps}};
 (-17,1)*{\scs \bullet}="1";
 (-17,-7)*{\scs \bullet}="2";
 (12,4)*{\scs \bullet}="tl";
 (20,4)*{\scs \bullet}="tr";
 (12,-5)*{\scs \bullet}="bl";
 (20,-5)*{\scs \bullet}="br";
  "1";"2" **\crv{}?(.58)*\dir{>};
  "tl";"tr" **\crv{}?(.5)*\dir{>};
  "bl";"br" **\crv{}?(.55)*\dir{>};
  "tl";"bl" **\crv{}?(.55)*\dir{>}+(-3,0)*{\scs \alpha_y};
  "tr";"br" **\crv{}?(.55)*\dir{>}+(3,0)*{\scs \alpha_x};
 \endxy $\\
  \hline
\end{tabular} }
\end{center}
}

\title{Lectures on $n$-Categories and Cohomology}
\author{Talks by John Baez, Notes by Michael Shulman}
\begin{document}
\maketitle

\tableofcontents 

\eject

\section*{Preface}
\label{sec:preface}
The goal of these talks was to explain how cohomology and other
tools of algebraic topology are seen through the lens of 
$n$-category theory.  
The talks were extremely informal, glossing over the difficulties
involved in making certain things precise, just trying to sketch the
big picture in an elementary way.  A lot of the material is
hard to find spelled out anywhere, but nothing new is due to me:
anything not already known by experts was invented by James Dolan,
Toby Bartels or Mike Shulman (who took notes, fixed lots of mistakes,
and wrote the Appendix).

The first talk was one of the 2006 Namboodiri Lectures
in Topology at the University of Chicago.
It's a quick introduction to the relation 
between Galois theory, covering spaces, cohomology, and higher 
categories.  The remaining talks, given in the category theory 
seminar at Chicago, were more advanced.
Topics include nonabelian cohomology, Postnikov towers, the 
theory of `$n$-stuff', and $n$-categories for $n = -1$ and $-2$.
Some questions from the audience have been included.  

Mike Shulman's extensive Appendix (\S\ref{sec:enrichm-posets-fiber})
clarifies many puzzles raised in
the talks.  It also ventures into deeper waters, such as the role of posets
and fibrations in higher category theory, alternate versions of the 
periodic table of $n$-categories, and the theory of higher topoi.
For readers who want more details, we have added an annotated
bibliography. --- JB

\vfill
\eject

\section{The Basic Principle of Galois Theory}
\label{sec:basic-principle-galois}

\subsection{Galois theory}
\label{sec:galois-erlangen}

Around 1832, Galois discovered a basic principle:

\begin{quote}
\textbf{We can study the ways a little thing $k$ can sit 
in a bigger thing $K$:
\[               k \hookrightarrow K   \]
by keeping track of the symmetries of $K$ that fix $k$.
These form a subgroup of the symmetries of $K$:
\[             \Gal(K|k) \subseteq \Aut(K).   \]
}
\end{quote}
For example, a point $k$ of a set $K$ is completely determined by the
subgroup of permutations of $K$ that fix $k$.  More generally, we
can recover any subset $k$ of a set $K$ from the subgroup of permutations 
of $K$ that fix $k$.

However, Galois applied his principle in a trickier example, namely
commutative algebra.  He took $K$ to be a field and $k$ to be a
subfield.  He studied this situation by looking at the subgroup
$\Gal(K|k)$ of automorphisms of $K$ that fix $k$.  Here this subgroup
does not determine $k$ unless we make a further technical assumption,
namely that $K$ is a `Galois extension' of $k$.  In general, we just
have a map sending each subfield of $K$ to the subgroup of $\Aut(K)$ that
fixes it, and a map sending each subgroup of $\Aut(K)$ to the subfield
it fixes.  These maps are not inverses; instead, they satisfy some
properties making them into what what is called a `Galois connection'.

When we seem forced to choose between extra technical assumptions or
less than optimal results, it often means we haven't fully understood the
general principle we're trying to formalize.  But, it can be very hard to 
take a big idea like the basic principle of Galois theory and express 
it precisely without losing some of its power.  That is not my goal
here.  Instead, I'll start by considering a weak but precise version 
of Galois' principle as applied to a specific subject: not commutative 
algebra, but {\it topology}.  

Topology isn't really separate from commutative algebra.
Indeed, in the mid-1800s, Dedekind, Kummer and Riemann realized that
commutative algebra is a lot like topology, only backwards.  Any space
$X$ has a commutative algebra $\scrO(X)$ consisting of functions on
it.  Any map
\[               f \maps X \to Y  \]
gives a map 
\[               f^\ast \maps \scrO(Y) \to \scrO(X)  .\]
If we're clever we can think of any commutative 
ring as functions on some space --- or on some `affine scheme':
\[   [\textrm{Affine Schemes}] = [\textrm{Commutative Rings}]^{\rm op}  .\]
Note how it's backwards: the {\it inclusion} of commutative rings
\[      p^\ast \maps \C[z] \hookrightarrow \C[\sqrt{z}]  \]
corresponds to the {\it branched cover} of the complex plane
by the Riemann surface for $\sqrt{z}$:
\[ 
\begin{array}{ccl}
         p \maps   \C &\to & \C    \\
              z &\mapsto& z^2 
\end{array}
\] 
So: classifying how a little commutative algebra can {\it sit inside} a 
big one amounts to classifying how a big space can {\it cover} a little one. 
Now the Galois group gets renamed the group of {\bf deck transformations}:
in the above example it's $\bbZ/2$:
\[              \sqrt{z} \mapsto -\sqrt{z}  .\]

The theme of `branched covers' became very important in later work on
number theory, where number fields are studied by analogy to function
fields (fields of functions on Riemann surfaces).  However, it's the
simpler case of unbranched covers where the basic principle of Galois
theory takes a specially simple and pretty form, thanks in part to
Poincar\'e.  This is the version we'll talk about now.  Later we'll 
generalize it from unbranched covers to `fibrations' of various sorts ---
fibrations of spaces, but also fibrations of $n$-categories.  Classifying
fibrations using the basic principle of Galois theory will eventually 
lead us to cohomology.

\subsection{The fundamental group}
\label{sec:fundamental group}

Around 1883, Poincar\'e discovered that any nice connected
space $B$ has a connected covering space 
that covers all others: its {\bf universal cover}.  This has
the biggest deck transformation group of all: the {\bf fundamental
group} $\pi_1(B)$.  

The idea behind Galois theory --- turned backwards! ---
then says that:

\begin{quote}
\textbf{Connected covering spaces of $B$
are classified by subgroups 
\[ H \subseteq \pi_1(B).\] 
}
\end{quote}

\noindent
This is the version we all learn in grad school.
To remove the `connectedness' assumption, we can start by
rephrasing it like this:

\begin{quote}
\textbf{Connected covering spaces of $B$
with fiber $F$ are classified by transitive actions of $\pi_1(B)$ 
on $F$.}
\end{quote}

\noindent
This amounts to the same thing, since 
transitive group actions are basically the same as subgroups:
given a subgroup $H \subseteq \pi_1(B)$ we can define $F$
to be $\pi_1(B)/H$, and given a transitive action of $\pi_1(B)$
on $F$ we can define $H$ to be the stabilizer group of a point.
The advantage of this formulation is that we can generalize it
to handle covering spaces where the total space isn't connected:

\begin{quote}
\textbf{Covering spaces of $B$ with fiber $F$ are classified
by actions of $\pi_1(B)$ on $F$:
\[          \pi_1(B) \to \Aut(F).   \]}
\end{quote}

\noindent
Here $F$ is any set and $\Aut(F)$ is the group of permutations
of this set:

\begin{center}
\begin{picture}(100,150)
  \includegraphics[scale=.5]{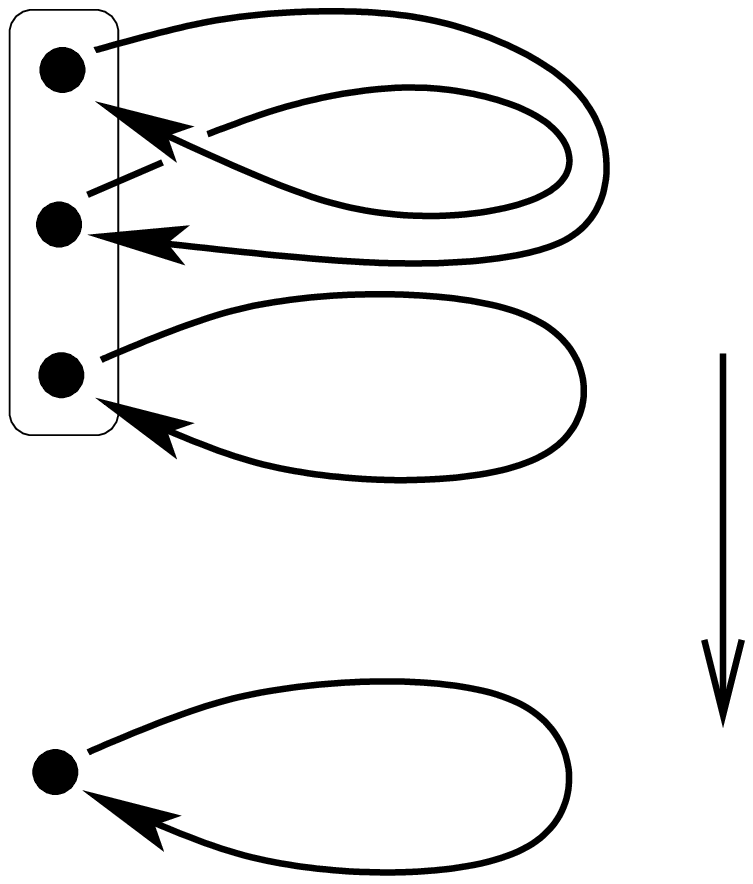}
  \put(-20,7){$B$}
  \put(-20,82){$E$}
  \put(-10,55){$p$}
  \put(-110,130){$F$}
  \end{picture}
\end{center}

\noindent
You can see how a loop in the base space gives a permutation of
the fiber.  The basic principle of Galois theory has become
`visible'!

\subsection{The fundamental groupoid}
\label{sec:fundamental groupoid}

So far the base space $B$ has been connected.  
What if $B$ is not connected?  For this, we should replace
$\pi_1(B)$ by $\Pi_1(B)$: the {\bf fundamental groupoid} of $B$.
This is the category where:

\begin{itemize}
\item
objects are points of $B$: 
$  \quad \bullet {\textstyle{{\textit{x}}}} $
\item
morphisms are homotopy classes of paths in $B$:  
\[
\xymatrix{
 {\textstyle{{\textit{x}}}}  \,
   \bullet \ar@/^1pc/[rr]^{f}
&& \bullet \, {\textstyle{{\textit{y}}}}   
}
\]
\end{itemize}

\vskip 1em
\noindent
The basic principle of Galois theory then says this:

\begin{quote}
\textbf{Covering spaces
$F \hookrightarrow E \to B$ 
are classified by actions of $\Pi_1(B)$ on $F$: that
is, functors
\[        \Pi_1(B) \to \Aut(F)  .\]
}
\end{quote}
 
\noindent
Even better, we can let the fiber $F$ be different
over different components of the base $B$:

\begin{quote}
\textbf{Covering spaces $E \to B$ 
are classified by functors
\[        \Pi_1(B) \to \Set  .\]
}
\end{quote}
 
What does this mean?  It says a lot in a very terse way.  
Given a covering space $ p \maps E \to B$, we can uniquely lift
any path in the base space to a path in $E$, given a lift
of the path's starting point.  Moreover, this lift depends
only on the homotopy class of the path.  So, our covering
space assigns a set $p^{-1}(b)$ to each point $b \in B$, and 
a map between these sets for any homotopy class of paths in $B$.
Since composition of paths gets sent to composition of maps,
this gives a {\it functor} from $\Pi_1(B)$ to $\Set$.  

Conversely, given any functor $F \maps \Pi_1(B) \to \Set$, we can use
it to cook up a covering space of $B$, by letting the fiber over
$b$ be $F(b)$, and so on.  So, with some work, we get a
one-to-one correspondence between isomorphism classes of
covering spaces $E \to B$ and natural isomorphism classes of
functors $\Pi_1(B) \to \Set$.  

But we actually get more: we get an {\it equivalence of categories.} 
The category of covering spaces of $B$ is equivalent to the
category where the objects are functors $\Pi_1(B) \to \Set$ 
and the morphisms are natural transformations between these guys.
This is what I've really meant all along by saying ``$X$'s are classified 
by $Y$'s.''  I mean there's a category of $X$'s, a category of $Y$'s,
and these categories are equivalent.  

\subsection{Eilenberg--Mac Lane Spaces}
\label{sec:Eilenberg-MacLane}

In 1945, Eilenberg and Mac Lane published their famous paper about
categories.  They {\it also} published a paper showing that any group
$G$ has a `best' space with $G$ as its fundamental group:
the {\bf Eilenberg-Mac Lane space} $K(G,1)$.  

In fact their idea is easiest to understand if we describe it a
bit more generally, not just for groups but for groupoids.
For any groupoid $G$ we can build a space 
$K(G,1)$ by taking a vertex for each object of $G$:
\[         \bullet\; x   \]
an edge for each morphism of $G$:
\[   \xy 0;/r.25pc/:
  (-8,0)*+{\bullet}="1";
  (0,5)*+{}; 
  (0,-5)*+{}; 
  (8,0)*+{\bullet}="2";
  {\ar "1";"2"^{f}};
 \endxy
\]
a triangle for each composable pair of morphisms:
\[   \vcenter{\xy 0;/r.25pc/:
   (-10,-5)*+{ \scriptstyle \bullet}="1";
   (10,-5)*+{ \scriptstyle \bullet}="2";
   (0,12)*+{ \scriptstyle \bullet}="3";
    {\ar "1"; "2"_{fg}};
    {\ar "3"; "2"^{g}};
    {\ar "1"; "3"^{f}};
    \endxy}
\]
a tetrahedron for each composable triple:
\[ \xy 0;/r.30pc/:
    (-10,-5 )*+{ \scriptstyle \bullet}="1";
    (8,-10)*+{ \scriptstyle \bullet}="2";
    (15,0)*+{ \scriptstyle \bullet}="3";
   (1,12)*+{ \scriptstyle \bullet}="4";
       {\ar "1";"2"_{fg} };
       {\ar "2";"3"_{h} };
       {\ar "4";"3"^{gh} };
    {\ar "1";"4"^{f} };
    {\ar "4";"2"_<<<<<<<<{g} };
       {\ar @{.>} "1";"3"^<<<<<<<<<{fgh}};
       \endxy
\]
and so on.  This space has $G$ as its fundamental
groupoid, and it's a {\bf homotopy
1-type}: all its homotopy groups above the 1st vanish.
These facts characterize it. 

Using this idea, one can show a portion of topology is just
groupoid theory: homotopy 1-types are the same as groupoids!
To make this precise requires a bit of work.  It's not
true that the category of homotopy 1-types and maps between
them is equivalent to the category of groupoids and functors
between them.  But, they form Quillen equivalent model categories.
Or, if you prefer, they form 2-equivalent 2-categories.

\subsection{Grothendieck's dream}
\label{sec:Grothendieck}

Since the classification of covering spaces
\[       E \to B  \]
only involves the fundamental groupoid of $B$, we might
as well assume $B$ is a homotopy 1-type.  Then $E$ will
be one too.

So, we might as well say $E$ and $B$ are {\it groupoids!}  
The analogue of a covering space for groupoids is a 
{\bf discrete fibration}: a functor 
$p \maps E \to B$ such that for any morphism 
$f \maps x \to y$ in $B$ and object $\tilde x \in E$ 
lifting $x$, there's a unique morphism 
$\tilde{f} \maps \tilde{x} \to \tilde{y}$ lifting $f$:

\begin{center}
\begin{picture}(100,150)
  \includegraphics[scale=.5]{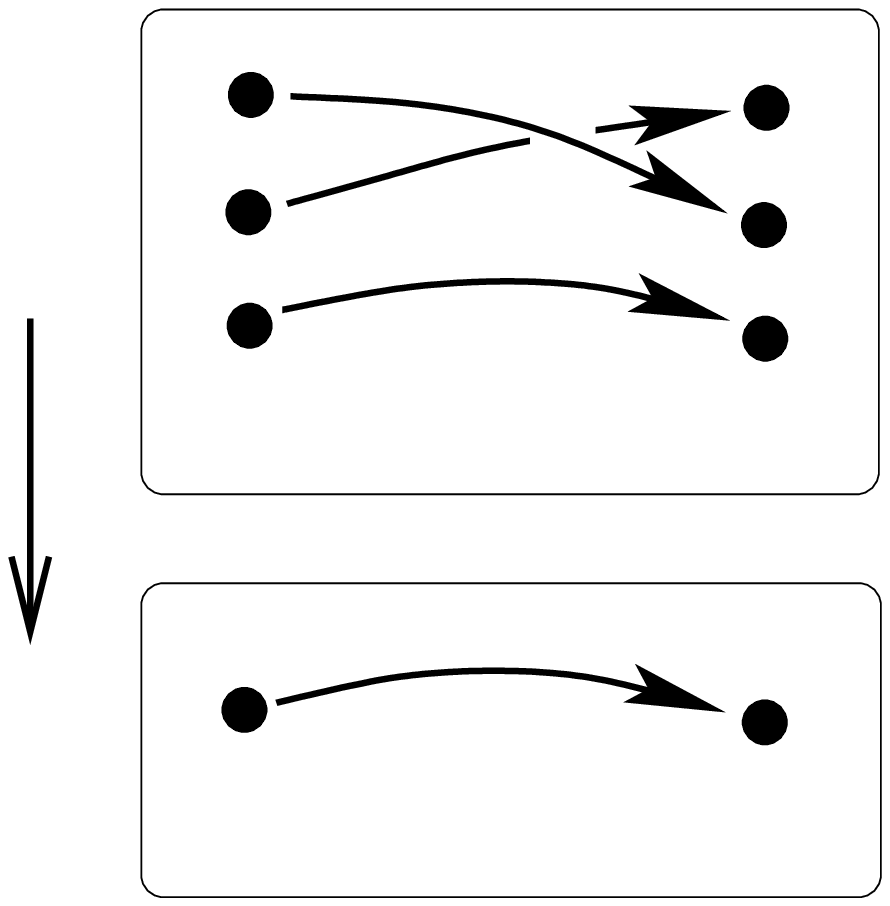}
  \put(-60,15){$f$}
  \put(-60,70){$\tilde f$}
  \put(-98,10){$x$}
  \put(-98,65){$\tilde x$}
  \put(-20,10){$y$}
  \put(-20,65){$\tilde y$}
  \put(-130,23){$B$}
  \put(-130,92){$E$}
  \put(-135,65){$p$}
  \end{picture}
  \vskip 1em
\end{center}

The basic principle of Galois theory then becomes:

\begin{quote}
\textbf{Discrete fibrations $E \to B$ are classified by functors
\[  B \to \textbf{\Set}. \]
}
\end{quote}
This is true even when $E$ and $B$ are categories, 
though then people use the term `opfibrations'.
This --- and much more --- goes back to Grothendieck's
1971 book {\it \'Etale Coverings and the Fundamental Group},
usually known as SGA1.

Grothendieck dreamt of a much bigger generalization of Galois 
theory in his 593-page letter to Quillen, {\it Pursuing Stacks}.  
Say a space is a \textbf{homotopy \textit{n}-type} 
if its homotopy groups above the $n$th all vanish.  Since
homotopy 1-types are `the same' as groupoids,
maybe homotopy $n$-types are `the same' as $n$-groupoids!
It's certainly true if we use Kan's simplicial approach to
$n$-groupoids --- but we want it to emerge from a general theory
of $n$-categories.  

For $n$-groupoids, the basic principle of Galois
theory should say something like this:

\begin{quote}
\textbf{Fibrations $E \to B$ where $E$ and $B$ are 
$n$-groupoids are classified by weak $(n+1)$-functors 
\[     B \to n\Gpd. \]}
\end{quote}

\noindent
Now when we say `classified by' we mean there's an
equivalence of $(n+1)$-categories.  `Weak' $n$-functors
are those where everything is preserved {\it up to equivalence}. 
I include the adjective `weak' only for emphasis: we
need all $n$-categories and $n$-functors to be weak 
for Grothendieck's dream to have any chance of coming true,
so for us, everything is weak by default.

Grothendieck made the above statement precise and proved it
for $n = 1$; later Hermida did it for $n = 2$.
Let's see what it amounts to when $n = 1$.  To keep things
really simple, suppose $E,B$ are 
just {\it groups}, and fix the fiber $F$, also a group.
With a fixed fiber our classifying 2-functor will land not in 
all of $\Gpd$, but in $\AUT(F)$, which is the `automorphism 2-group' 
of $F$ --- I'll say exactly what that is in a minute.
In this simple case fibrations are just {\it extensions},
so we get a statement like this:

\begin{quote}
\textbf{Extensions of the group $B$ by the group 
$F$, that is, short exact sequences 
\[          1 \to F \to E \to B \to 1,  \]
are classified by weak 2-functors 
\[               B \to \AUT(F) .\]
}
\end{quote}

\noindent
This is called {\bf Schreier theory}, since a version of this result
was first shown by Schreier around 1926.  The more familiar classifications 
of abelian or central group extensions using ${\rm Ext}$ or $H^2$ are just
watered-down versions of this.

$\AUT(F)$ is the {\bf automorphism 2-group} of $F$, a 2-category 
with: 
\begin{itemize}
\item $F$ as its only object:
$  \quad \bullet {\textstyle{{\textit{F}}}} $
\item automorphisms of $F$ as its morphisms:
\[
\xymatrix{
 {\textstyle{{\textit{F}}}}  \,
   \bullet \ar@/^1pc/[rr]^{\alpha}
&& \bullet \, {\textstyle{{\textit{F}}}}   
}
\]
\item elements $g \in F$ 
with $g\alpha(f)g^{-1} = \beta(f)$ as its 2-morphisms:
\[
\xymatrix{
 {\textstyle{{\textit{F}}}} \, 
   \bullet \ar@/^1pc/[rr]^{\alpha}_{}="0"
           \ar@/_1pc/[rr]_{\beta}^{}="1"
&& \bullet \, {\textstyle{{\textit{F}}}}   
\ar@{=>}"0";"1"^{g}
}
\]
\end{itemize}
In other words, we get $\AUT(F)$ by taking $\Gpd$ and
forming the sub-2-category with $F$ as its only object,
all morphisms from this to itself, and all 2-morphisms
between these.

\noindent
Given a short exact sequence of groups, we classify it by 
choosing a {\bf set-theoretic section}:
\[
\xymatrix{  
1 \ar[r] & F \ar[r]_{i} & E \ar[r]_{p} & B 
\ar@/_1pc/[l]_{s}
\ar[r] & 1 }, \]
meaning a function $s \maps B \to E$ with $p(s(b)) = b$
for all $b \in B$.  This gives for any $b \in B$ an automorphism 
$\alpha(b)$ of $F$:
\[   \alpha(b)(f) = s(b) f s(b)^{-1} .\]
Since $s$ need not be a homomorphism, we may not have
\[    \alpha(b) \, \alpha(b') = \alpha(bb')   \]
but this holds {\it up to conjugation} by an element
$\alpha(b,b') \in F$.  That is, 
\[    \alpha(b,b') \, [\alpha(b) \, \alpha(b') f] \, \alpha(b,b')^{-1} = 
\alpha(bb')f   \]
where 
\[           \alpha(b,b') = s(bb') \, (s(b) \, s(b'))^{-1} .\]
This turns out to yield a {\it weak} 2-functor
\[       \alpha \maps B \to \AUT(F)  .\]
If we consider two weak 2-functors equivalent when there's
a `weak natural isomorphism' between them, different choices of 
$s$ will give equivalent 2-functors.  Isomorphic extensions of $B$ by $F$ 
also give equivalent 2-functors.  

The set of equivalence classes of weak 2-functors $B \to \AUT(F)$ is
often called the {\bf nonabelian cohomology} $H(B,\AUT(F))$.  
So, we've described a map sending isomorphism classes of short
short exact sequences 
\[  1 \to F \to E \to B \to 1 \]
to elements of $H(B,\AUT(F))$.   And, 
this map is one-to-one and onto.

This is {\it part} of what we mean by saying extensions of 
$B$ by $F$ are classified by weak 2-functors $B \to \AUT(F)$.  
But as usual, we really have something much better: an 
equivalence of 2-categories.  

There's a well-known category of extensions of $B$ by $F$,
where the morphisms are commutative diagrams like this:
\[
\xymatrix{
  && E \ar[dr] \ar[dd]\\
  1 \ar[r] & F\ar[ur]\ar[dr]  && B \ar[r] & 1\\
  && E' \ar[ur]
}\]
But actually, we get a 2-category of extensions using
the fact that groups are special groupoids and $\Gpd$
is a 2-category.  Similarly, $\hom(B,\AUT(F))$ is a
2-category, since $\AUT(F)$ is a 2-category and $B$ is 
a group, hence a special sort of category, hence a special
sort of 2-category.  And, the main result of Schreier 
theory says the 2-category of extensions of $B$ by $F$ is
equivalent to $\hom(B,\AUT(F))$.  This implies the earlier 
result we stated, but it's much stronger.  

In short, generalizing the fundamental principle of Galois
theory to fibrations where everything is a {\it group}
gives a beautiful classification of group extensions
in terms of nonabelian cohomology.  In the rest of these
lectures, we'll explore how this generalizes as we go from
groups to $n$-groupoids and $n$-categories.

\vfill 
\eject

\section{The Power of Negative Thinking}
\label{sec:power-neg-thinking}

\subsection{Extending the periodic table}
\label{sec:extend-perod-table}

Now I want to dig a lot deeper into the relation between
fibrations, cohomology and $n$-categories.  At this point
I'll suddenly assume that you have some idea of what 
$n$-categories are, or at least can fake it.
The periodic table of $n$-categories shows what various degenerate
versions of $n$-category look like.   We can think of
an $(n+k)$-category with just one $j$-morphism for $j < k$ as
a special sort of $n$-category.  They look like this:

\vskip 1em

\vbox{
\begin{center}
{\bf THE PERIODIC TABLE }
\vskip 1em
\end{center}
\noindent \periodictable
\vskip 1em
}

\noindent
For example, in the $n = 0$, $k= 1$ spot we have
$1$-categories with just one $0$-morphism, or in normal
language, categories with just one object --- i.e., monoids.
Indeed a `monoid' is the perfect name for a one-object category,
because `monos' means `one' --- but that's not where the 
name comes from, of course.   It's a good thing, too, or else
Eilenberg and Mac Lane might have called categories `polyoids'.

We'll be thinking about all these things in the most weak manner
possible, so `2-category' means `weak 2-category', aka
`bicategory', and so on.  Everything I'm going to tell you
\emph{should} be true, once we really understand what is going on.
Right now it's more in the nature of dreams and speculations, but I
don't think we'll be able to prove the theorems until we dream enough.

Eckmann and Hilton algebraicized a topological argument going back
to Hurewicz, which proves that strict monoidal categories with one object are
commutative monoids.  Eugenia Cheng and Nick Gurski have studied
this carefully, and they've shown that things are a little more
complicated when we consider \emph{weak} monoidal categories, 
but I'm going to proceed in a robust spirit and leave such issues
to smart young people like them.

Things get interesting in the second column, when we get \emph{braided}
monoidal categories.  These are not the most obvious sort of
`commutative' monoidal categories that Mac Lane first wrote down,
namely the symmetric ones.  James Dolan and I were quite
confused about why braided and symmetric both exist, until we started
getting the hang of the periodic table.  

Noticing that in the first column we stabilize after 2 steps at
commutative monoids, and in the second after 3 steps at symmetric
monoidal categories, we enunciated the {\bf stabilization hypothesis},
which says that the $n$th column should stabilize at the $(n+2)$nd
row.  We believed this because of the Freudenthal suspension theorem
in homotopy theory, which says that if you keep suspending and looping
a space, it gets nicer and nicer, and if it only has $n$ nonvanishing
homotopy groups, eventually it's as nice as it can get, and it
stabilizes after $n+2$ steps.  This is related, because any space
gives you an $n$-groupoid of points, paths, paths-of-paths, and so on.

Such a beautiful pattern takes the nebulous, scary subject of
$n$-categories and imposes some structure on it.  There are all sorts
of operations that take you hopping between different squares of this
chart.

We call this chart the `Periodic Table' --- not because it's periodic, 
but because we can use it to predict new phenomena, like Mendeleev used
the periodic table to predict new elements.

After we came up with the periodic table
I showed it to Chris Isham, who does quantum gravity at Imperial College.
I was incredibly happy with it, but he said: ``That's 
obviously not right --- you didn't start the chart at the right place.
First there should be a column with just one interesting row,
then a column with two, and {\it then} one with three.''  

I thought he was crazy, but it kept nagging me.  It's sort of 
weird to start counting at three, after all.
But there are no $(-1)$-categories or $(-2)$-categories!  Are there?

It turns out there are!
Eventually Toby Bartels and James Dolan figured out what they are.
And they realized that Isham was right.  The periodic table really
looks like this:

\vskip 2em
\vbox{
\begin{center}
{\bf THE EXTENDED PERIODIC TABLE }
\vskip 1em
\end{center}
\noindent \extendedperiodictable
\vskip 1em
}

\noindent
You should be dying to know what fills in those question marks. 
Just for fun, I'll tell you what two of them are now.  You 
probably won't think these answers are obvious --- but you will 
soon:

\begin{itemize}
\item
$(-1)$-categories are just truth values: there are only two of them,
True and False.  \vskip 1em
\item
$(-2)$-categores are just `necessarily true' truth values: there is
only one of them, which is True.  
\end{itemize}

\noindent
I know this sounds crazy, but it sheds lots of light on many
things.  Let's see why $(-1)$- and $(-2)$-categories really work this 
way.

\subsection{The categorical approach}
\label{sec:categorical-approach}

Before describing $(-1)$-categories and $(-2)$-categories, we need to
understand a couple of facts about the $n$-categorical world.

The first is that in the $n$-categorical universe, every $n$-category
is secretly an $(n+1)$-category with only identity $(n+1)$-morphisms.
It's common for people to talk about sets as discrete categories, for
example.  A way to think about it is that these identity
$(n+1)$-morphisms are really \emph{equations}.  

When you play the $n$-category game, there's a rule that you should
never say things are equal, only isomorphic.  This makes sense up
until the top level, the level of $n$-morphisms, when you break down 
and allow yourself to say that $n$-morphisms are equal.  But actually
you aren't breaking the rule here, if you think of your $n$-category
as an $(n+1)$-category with only identity morphisms.  
Those equations are really isomorphisms: it's just that the only
isomorphisms existing at the $(n+1)$st level are identities.  Thus when we
assert equations, we're refusing to think about things still more
categorically, and saying ``all I can take today is an $n$-category''
rather than an $(n+1)$-category.

We can iterate this and go on forever, so every $n$-category is really
an $\infty$-category with only identity $j$-morphisms for $j>n$.

The second thing is that big $n$-categories have lots of little
$n$-categories inside them.  For example, between two objects 
$x, y$ in a 3-category, there's a little 2-category
$\hom(x,y)$.  James and I jokingly call this sort of thing 
a `microcosm', since it's like a little world within a world:

\vskip 1em
\vbox{
\begin{center}
   $\xy 0;/r.22pc/:
  (-20,0)*+{x \bullet}="1";
  (0,0)*+{\bullet y}="2";
  {\ar@/^2pc/ "1";"2"};
  {\ar@/_2pc/ "1";"2"};
   (-10,6)*+{}="A";
   (-10,-6)*+{}="B";
  {\ar@{=>}@/_.7pc/ "A"+(-2,0) ; "B"+(-1,-.8)};
  {\ar@{=}@/_.7pc/ "A"+(-2,0)  ; "B"+(-2,0)};
  {\ar@{=>}@/^.7pc/ "A"+(2,0)  ; "B"+(1,-.8)};
  {\ar@{=}@/^.7pc/ "A"+(2,0)  ; "B"+(2,0)};
  {\ar@3{->} (-12,0)*{}; (-8,0)*{}};
 \endxy$ 

\medskip
A microcosm
\end{center}
}

\vskip 1em
\noindent
In general, given objects $x$, $y$ in an $n$-category,
there is an $(n-1)$-category called $\hom(x,y)$ (because it's the
`thing' of morphisms from $x$ to $y$) with the morphisms $f\maps x\to y$
as its objects, etc..  

We can iterate this and look at microcosms of microcosms.  A couple of
objects in $\hom(x,y)$ give an $(n-2)$-category, and so on.  It's 
handy to say that two $j$-morphisms $x$ and $y$ are {\bf parallel} 
if they look like this:
\[\xymatrix{
  z\bullet \rtwocell^{x}_{y}{} & \bullet z'
}\]
This makes sense if $j > 0$; if $j = 0$ we decree that all 
$j$-morphisms are parallel.  The point is that it only makes sense
to talk about something like $f \maps x \to y$ when $x$ and $y$ are parallel.
Given parallel $j$-morphisms $x$ and $y$, we get an $(n-j-1)$-category 
$\hom(x,y)$ with $(j+1)$-morphisms $f\maps x\to y$ as objects, and so on.
This is a little microcosm.

In short: {\it given parallel $j$-morphisms $x$ and $y$ in an $n$-category,
$\hom(x,y)$ is an $(n-j-1)$-category.}
Now take $j=n$.  Then we get a $(-1)$-category!  If $x$ and $y$ are
parallel $n$-morphisms in an $n$-category, then $\hom(x,y)$ is a
$(-1)$-category.  What is it?

You might say ``that's cheating: you're not allowed to go that
high.''  But it isn't really cheating, 
since, as we said, every $n$-category is
secretly an $\infty$-category.  We just need to work out the answer,
and that will tell us what a $(-1)$-category is.

The objects of $\hom(x,y)$ are $(n+1)$-morphisms, which here are just
identities.  So, if $x=y$ there is one object in $\hom(x,y)$,
otherwise there's none.  So there are really just two possible
$(-1)$-categories.  There aren't any $(-1)$-categories that don't
arise in this way, since in general any $n$-category can be stuck 
in between two objects to make an $(n+1)$-category.

Thus, there are just two $(-1)$-categories.  You could think of them as
the 1-element set and the empty set, although they're not exactly
sets.  We can also call them `$=$' and `$\neq$', or `True' and
`False.'  The main thing is, there are just two.

I hope I've convinced you this is right.  You may think it's silly,
but you should think it's right.

Now we should take $j= n+1$ to get the $(-2)$-categories.  If we have
two parallel $(n+1)$-morphisms in an $n$-category, they are both
identities, so being parallel they must both be 
$1_z \maps z\to z$ for some $z$, so they're equal.
At this level, they \emph{have} to be equal, so there \emph{is} an
identity from one from to the other, necessarily.  It's
like the previous case, except we only have one choice.

So there's just one $(-2)$-category.  When there's just one of
something, you can call it anything you want, since you don't have to
distinguish it from anything.  But it might be good to call 
this guy the 1-element set, or `$=$', or `True'.  Or maybe we should
call it `Necessarily True', since there's no other choice.

So, we've worked out $(-1)$-categories and $(-2)$-categories.  
We could keep on going, but it stabilizes past this point:
for $n>2$, $(-n)$-categories are all just `True.'
I'll leave it as a puzzle for you to figure out what a monoidal
$(-1)$-category is.

\subsection{Homotopy $n$-types}
\label{sec:homotopy-n-types}

I really want to talk about what this has to do with topology.  We're
going to study very-low-dimensional algebra and apply it to
very-low-dimensional topology---in fact, so low-dimensional that they
never told you about it in school.  Remember Grothendieck's dream:

\begin{hyp}[Grothendieck's Dream, aka the Homotopy Hypothesis]
$n$-Groupoids are the same as homotopy $n$-types.
\end{hyp}

Here `$n$-groupoids' means \emph{weak} $n$-groupoids, in which
everything is invertible up to higher-level morphisms, in the
weakest possible way.  Similarly `the same' is meant in the weakest
possible way, which we might make precise using something
called `Quillen equivalence'.  A {\bf homotopy $n$-type} is
a nice space (e.g.\ a CW-complex) with vanishing homotopy groups
$\pi_j$ for $j>n$.  (If it's not connected, we have to take $\pi_j$ at
every basepoint.)  We could have said `all spaces' instead of `nice
spaces,' but then we'd need to talk about \emph{weak} homotopy
equivalence instead of homotopy equivalence.

People have made this precise and shown that it's true for various low
values of $n$, and we're currently struggling with it for higher
values.  It's known (really well) for $n=1$, (pretty darn well) for
$n=2$, (partly) for $n=3$, and somewhat fuzzier for higher $n$.  Today
we're going to do it for \emph{lower} values of $n$, like $n = -1$
and $n= -2$.  

Now, they never told you about the negative second homotopy group of
a space---or `homotopy thingy', since after all we know $\pi_0$ is 
only a set, not a group.   In fact we won't define these negative 
homotopy thingies as thingies yet, but only what it means for them to vanish:

\begin{defn}
  We say {\bf $\pi_j(X)$ vanishes for all basepoints} if given any
  $f\maps S^j\to X$ there exists $g\maps D^{j+1}\to X$ extending $f$.
  We say $X$ is a {\bf homotopy $n$-type} if $\pi_j(X)$ vanishes
  for all basepoints whenever $j > n$.
\end{defn}

We'll use this to figure out what a homotopy 0-type is, then use it
to figure out what a homotopy $(-1)$-type and a homotopy $(-2)$-type
are.  

$X$ is a \emph{homotopy 0-type} when all circles and higher-dimensional
spheres mapped into $X$ can be contracted.  So, $X$ is just a disjoint
union of connected components, all of which are contractible.  (The
fact that being able to contract all spheres implies the space is
contractible, for nice spaces, is Whitehead's theorem.)  From the
point of view of a homotopy theorist, such a space might as well just
be a set of points, i.e.\ a discrete space, but the points could be
`fat'.  This is what Grothendieck said should happen; all
0-categories are 0-groupoids, which are just sets.

Now let's figure out what a homotopy $(-1)$-type is.
If you pay careful attention, you'll see the following argument 
is sort of the same as what we did
before to figure out what $(-1)$-categories are.

By definition,
a \emph{homotopy $(-1)$-type} is a disjoint union of contractible spaces
(i.e.\ a homotopy 0-type) with the extra property
that maps from $S^0$ can be
contracted.  What can $X$ be now?  It can have just one contractible
component (the easy case)---or it can have none (the sneaky case).  So
$X$ is a disjoint union of 0 or 1 contractible components.  From the
point of view of homotopy theory, such a space might as well be an
empty set or a 1-point set.  This is the same answer that we got
before: the absence or presence of an equation, `False' or ``True'.

Finally, a \emph{homotopy $(-2)$-type} is thus a space like this such
that any map $S^{-1}\to X$ extends to $D^0$.  Now we have to remember
what $S^{-1}$ is.  The $n$-sphere is the unit sphere in $\bbR^{n+1}$,
so $S^{-1}$ is the unit sphere in $\bbR^0$.  It consists of all unit
vectors in this zero-dimensional vector space, i.e.\ it is the empty
set, $S^{-1}=\emptyset$.  Well, a map from the empty set into $X$ is a
really easy thing to be given: there's always just one.  And $D^0$ is
the unit disc in $\bbR^0$, so it's just the origin, $D^0=\{0\}$.  So
this extension condition says that $X$ has to have at least one point
in it.  Thus a homotopy $(-2)$-type is a disjoint union of precisely
one contractible component.  Up to homotopy, it is thus a one-point
set, or `True'.  So Grothendieck's idea works here too.

Now I think you all agree; I gave you \emph{two} proofs!

What about $(-3)$-types?  Even I get a little scared about $S^{-2}$
and $D^{-1}$.  

\vskip 1em
[Peter May: {\sl They're both empty, so it stabilizes here.}]

\subsection{Stuff, structure, and properties}
\label{sec:stuff-struct-prop}

What's all this nonsense about?  In math we're often interested in
equipping things with extra structure, stuff, or properties, and
people are often a little vague about what these mean.  For example, a
group is a set (\emph{stuff}) with operations (\emph{structure}) such
that a bunch of equations hold (\emph{properties}).

You can make these concepts very precise by thinking about forgetful
functors.  It always bugged me when reading books that no one ever
defined `forgetful functor'.  Some functors are more forgetful than
others.  Consider a functor $p\maps E\to B$ (the notation reflects that
later on, we're going to turn it into a fibration when we use
Grothendieck's idea).  There are various amounts of forgetfulness
that $p$ can have:

\begin{itemize}
\item $p$ \textbf{forgets nothing} if it is an equivalence of categories,
  i.e.\ faithful, full, and essentially surjective.  For example
   the identity functor $\mathbf{AbGp}\to\mathbf{AbGp}$ forgets nothing.
\item $p$ \textbf{forgets at most properties} if it is faithful and
  full.  E.g.\ $\mathbf{AbGp}\to\mathbf{Gp}$, which forgets the
  property of being abelian, but a homomorphism of abelian groups is
  just a homomorphism between groups that happen to be abelian.
\item $p$ \textbf{forgets at most structure} if it is faithful.  E.g.\
  the forgetful functor from groups 
  to sets, $\mathbf{AbGp}\to\mathbf{Sets}$, forgets the structure of
  being an abelian group, but it's still faithful.
\item $p$ \textbf{forgets at most stuff} if it is arbitrary.  E.g.
  $\textbf{Sets}^2\to\mathbf{Sets}$, where we just throw out the
  second set, is not even faithful.
\end{itemize}

There are different ways of slicing this pie.
For now, we are thinking of each level of forgetfulness as
subsuming the previous ones, 
so `forgetting at most structure' means forgetting structure 
and/or properties and/or nothing, but we can also try to make
them completely disjoint concepts.  Later I'll define
a concept of `forgetting {\it purely} structure' and so on.

What's going on here is that in every case, what you can do is take
objects downstairs and look at their \emph{fiber} or really
\emph{homotopy fiber} upstairs.  An object in the homotopy fiber
upstairs is an object together with a morphism from its image to the
object downstairs, as follows.

Given $p\maps E\to B$ and $x\in B$, its \textbf{homotopy fiber} or
\textbf{essential preimage}, which we write $p\inv(x)$, has:
\begin{itemize}
\item objects $e\in E$ equipped with isomorphism $p(e)\iso x$
\item morphisms $f\maps  e \to e'$ in $E$ compatible with the 
  given isomorphisms:
  \[\xymatrix{p(e) \ar[rr]^{p(f)} \ar[dr]_\iso && p(e')\ar[dl]^\iso \\ & x}\]
\end{itemize}

It turns out that the more forgetful the functor is, the bigger and
badder the homotopy fibers can be.  In other words, switching to the 
language of topology, they can have bigger \textbf{homotopical dimension}: 
they can have nonvanishing homotopy groups up to dimension $d$ for bigger 
$d$.  

\begin{fact}\label{forgetfulness-by-fibers}
  If $E$ and $B$ are groupoids (we'll consider the case of categories
  in \S\ref{sec:fibers}), then
  \begin{itemize}
  \item $p$ forgets stuff if all $p\inv(x)$ are arbitrary
    (1-)groupoids.  For example, there's a whole \emph{groupoid} 
    of ways to add an extra set to some set.
  \item $p$ forgets structure if all $p\inv(x)$ are 0-groupoids,
    i.e.\ groupoids which are (equivalent to) sets.  For example, 
    there's just a \emph{set} of ways of making a set into a group.
  \item $p$ forgets properties if all $p\inv(x)$ are $(-1)$-groupoids.
    For example, there's just a \emph{truth value} of ways of making 
    a group into an abelian group: either you can or you can't 
    (i.e.\ it is or it isn't).
  \item $p$ forgets nothing if all $p\inv(x)$ are $(-2)$-groupoids.
For example, there's just a \emph{`necessarily true' truth value} of ways
of making an abelian group into an abelian group: you always can, 
in one way.
  \end{itemize}
(In the examples above, we are considering the \emph{groupoids} of
sets, pairs of sets, groups and abelian groups.
We'll consider the case of categories in \S\ref{sec:fibers}.)
\end{fact}

We can thus study how forgetful a functor is by looking at what
homotopy dimension its fibers have.

Note that to make this chart work, we really needed the negative
dimensions.  You should want to also go in the other direction, say if
we had 2-groupoids or 3-groupoids; then we'll have something `even
more substantial than stuff.'  James Dolan dubbed that \emph{eka-stuff},
by analogy with how Mendeleev called elements which were missing in
the periodic table `eka-?', e.g.\ `eka-silicon' for the missing
element below silicon, which now we call germanium.  He
guessed that eka-silicon would be a lot like silicon, but heavier, and
so on.  Like Mendeleev, we can use the periodic table to guess things,
and then go out and check them.

\subsection{Questions and comments}
\label{sec:questions-comments}

\subsubsection{What should forgetful mean?}
\label{sec:what-should-forg}
{\quad}

\vskip 1em
Peter May: \emph{Only functors with left adjoints should really be called
`forgetful.'  Should the free group functor be forgetful?}

\vskip 1em

MS: \emph{A set is the same as a group with the property of 
being free and the structure of specified generators.}

\vskip 1em
JB: Right, you can look at it this way.  Then the free functor 
$F \maps \Set \to \Grp$ `forgets at most structure': it's faithful,
but neither full nor essentially surjective.  Whether you want 
to call it `forgetful' is up to you, but this is how I'm using
the terminology now.

\subsubsection{Monoidal $(-1)$-categories}
\label{sec:monoidal--1}
\quad

\vskip 1em
Puzzle: What's a monoidal $(-1)$-category?

\vskip 1em
Answer: A $(-1)$-category is a truth value, and the only
monoidal $(-1)$-category is True.

\vskip 1em
To figure this out, note that a monoidal $(-1)$-category is 
what we get when we take a $0$-category with just one object,
say $x$, and look at $\hom(x,x)$.  A $0$-category with one object
is just a one-element set $\{x\}$, and $\hom(x,x)$ is just
the equality $x=x$, which is `True.'

Note that monoidal $(-1)$-categories are stable, so we are just adding
a property to the previous one, just like when we pass from monoids to
commutative monoids, or braiding to symmetry.

In general, we have forgetful functors marching \emph{up} the periodic
table, which forget different amounts of things.  We forget nothing
until we get up to the end of the stable range, then we forget a
property (symmetry, or commutativity), then a structure (monoid
structure, or braiding), then stuff (a monoidal structure, of which
there are a whole category of ways to add to a given category), then
eka-stuff (ways to make a 2-category into a monoidal 2-category), and
so on.

\subsubsection{Maps of truth values}
\label{sec:maps-truth-values}

Here's a funny thing.
Note that for categories, we have a composition \emph{function}
$$\hom(x,y)\times\hom(y,z)\to \hom(x,z).$$
For a set, the substitute is transitivity:
$$(x=y) \; \& \; (y=z) \Rightarrow (x=z).$$
In other words, we can `compose equations.'
But here $\Rightarrow$ is acting as a map between truth values.
What sort of morphisms of truth values do we have?  We just have a
0-category of $(-1)$-categories, so there should be only identity morphisms.
The implication $F \Rightarrow T$ doesn't show up in this story.

That's a bit sad.  Ideally, a bunch of
\emph{propositional logic} would show up at the level of $-1$-categories.
Toby Bartels has a strategy to fix this.  In his approach, posets
play a much bigger role in the periodic table, to include the
notion of implication between truth values.   

In fact, it seems that the periodic table is just a slice of
a larger 3-dimensional table relating higher categories and logic...
see \S\ref{sec:enrichment-posets} for more on this. 

\section{Cohomology: The Layer-Cake Philosophy}
\label{sec:fact-postn-towers}

We're going to continue heading in the direction of cohomology, but
we'll get there by a perhaps unfamiliar route.  Last time we led up to 
the concept of \emph{$n$-stuff}, although we stopped right after
discussing ordinary stuff and only mentioned eka-stuff briefly.

\subsection{Factorizations}
\label{sec:factorizations}

What does it mean to forget \emph{just} stuff, or \emph{just}
properties?

\begin{defn}
  An $\infty$-functor $p\maps E\to B$ is \textbf{$n$-surjective} (perhaps
  \textbf{essentially $n$-surjective}) if given any parallel
  $(n-1)$-morphisms $e$ and $e'$ in $E$, and any $n$-morphism $f\maps p(e)\to
  p(e')$, there is an $n$-morphism $\ftil\maps e\to e'$ such that $p(\ftil)$
  is equivalent to $f$: $p(\ftil)\eqv f$.   
\end{defn}

For example, suppose $p\maps E\to B$ is a function between \emph{sets}.  It
is
\begin{itemize}
\item 0-surjective if it is \emph{surjective} in the usual sense, since in
  this case equivalences are just equalities.  (The presence of $e$
  and $e'$ here is a bit confusing, unless you believe that all
  $n$-categories go arbitrarily far down as well.)
\item 1-surjective if it is \emph{injective}; since in this case all
  1-morphisms are identities, 1-surjective means that if $p(e)=p(e')$,
  then $e=e'$.
\end{itemize}

That's a nice surprise: \emph{injective means `surjective on
equations'!}  

Now, another thing that you can do in the case of sets is to take any
old function $p\maps E\to B$ and factor it as first a surjection and then
an injection:
\[\xymatrix{
  E \ar[dr]_{\mathrm{0-surj}} \ar[rr]^p & & B\\
  & E' = \im(p) \ar[ur]_{\mathrm{1-surj}}
}\]
The interesting thing is that this keeps generalizing as we go up to
higher categories.

To see how this works, first suppose $E$ and $B$ are categories.  
A functor $p\maps E\to B$ is:
\begin{itemize}
\item 0-surjective if it is essentially surjective;
\item 1-surjective if it is full;
\item 2-surjective if it is faithful;
\item 3-surjective always, and so on.
\end{itemize}

Do you see why?
\begin{itemize}
\item
Crudely speaking, 0-surjective means `surjective on objects'.  But
you have to be a bit careful: it's sufficient that every object in 
$B$ is \emph{isomorphic} to $p(e)$ for some $e \in E$.  So,
0-surjective really means \emph{essentially} surjective.
\item
Crudely speaking,
1-surjective means `surjective on arrows,' but you just have to be a
little careful: it's not fair to ask that an arrow downstairs be the
image of something unless we already know that its source and target
are the images of something.  So, it really means our functor
is \emph{full}
\item
Similarly, 2-surjective means `surjective on equations
between morphisms,' i.e.\ injective on hom-sets.  So, it means
our functor is \emph{faithful}.  
\end{itemize}
Note that the conjunction of all three of these means our functor
is an equivalence, just as a surjective and injective function is
an isomorphism of sets.

The notions of forgetting at most stuff, structure,
or properties can also be defined using conjunctions of these
conditions, namely:
\begin{itemize}
\item Forgets nothing: 0, 1, and 2-surjective
\item Forgets at most properties: 1 and 2-surjective
\item Forgets at most structure: 2-surjective
\item Forgets stuff: arbitrary
\end{itemize}

As with functions between sets, we get a factorization result for
functors between categories.  Any functor factors like this:
\[\xymatrix{
E\ar[dr]_{\mathrm{0, 1-surj}}\ar[rrrr]^p &&&& B\\
& E' \ar[rr]_{\mathrm{0, 2-surj}} && E'' \ar[ur]_{\mathrm{1, 2-surj}}
}\]
I think this is a well-known result.
You build these other categories as `hybrids' of $E$ and $B$:
$E$ gradually turns into $B$ from the top down.
We start with $E$; then we throw in new 2-morphisms 
(equations between morphisms) that we get from $B$; then 
we throw in new 1-morphisms (morphisms), and finally 
new 0-morphisms (objects).  It's like a horse transforming into 
a person from the head down.   First it's a
horse, then it's a centaur, then it's a faun-like thing that's horse 
from the legs down, and finally it's a person.

In more detail:
\begin{itemize}
\item The objects of $E'$ are the same objects as $E$, but a morphism
  from $e$ to $e'$ is a morphism $p(e)\to p(e')$ in $B$ which is in
  the image of $p$; and
\item The objects of $E''$ are the objects of $B$ in the (essential)
  image of $p$, with all morphisms between them.
\end{itemize}
This is the same thing as is happening for sets, but there it's
happening so fast that you can't see it happening.

The next example will justify this terminology:

\begin{itemize}
\item A functor which is 0- and 1-surjective \textbf{forgets purely stuff};
\item A functor which is 0- and 2-surjective \textbf{forgets purely structure};
\item A functor which is 1- and 2-surjective \textbf{forgets purely properties}.
\end{itemize}

\begin{eg}
  Let's take the category of \emph{pairs of vector spaces} and forget
  down to just the underlying set of the \emph{first} vector space (so
  that we have an interesting process at every stage).
  \[\xymatrix{
    \mathrm{Vect}^2\ar[dr]\ar[rrr]^p &&& \mathrm{Set}\\
    & E' \ar[r] & E'' \ar[ur]
  }\]

  The objects of $E'$ are again pairs of vector spaces, but its
  morphisms are just linear maps between the \emph{first} ones.  We
  write this as 
  \[
  \big[\text{pairs of vector spaces, linear maps between first ones}\big].
  \]
  In fact, this category is equivalent to Vect; the extra vector space
  doesn't participate in the morphisms, so it might as well not be
  there.  Our factorization cleverly managed to forget the stuff (the
  second vector space) but still keep the structure on what remains.

  The objects of $E''$ are sets with the \emph{property} that they can
  be made into vector spaces, and its morphisms are arbitrary
  functions between them.  Here we forgot just the structure of being
  a vector space, but we cleverly didn't forget the \emph{property} of
  being vector-space-izable.  That's sort of cool.
\end{eg}

Having seen how these factorizations work for sets and categories,
we can guess how they go for $\infty$-categories.  Let's say an 
$\infty$-functor $p \maps E \to B$ {\bf forgets purely $j$-stuff}
if it's $i$-surjective for all $i \ne j+1$.
Note the funny `$+1$' in there: we need this to make things
work smoothly.  For ordinary categories, we have:
\begin{itemize}
\item A functor that `forgets purely $1$-stuff' forgets purely stuff;
\item A functor that `forgets purely $0$-stuff' forgets purely structure;
\item A functor that `forgets purely $(-1)$-stuff' forgets purely properties;
\item A functor that `forgets purely $(-2)$-stuff' forgets nothing.
\end{itemize}
If we go up to 2-categories we get a new concept, `$2$-stuff'.   This 
is what we called `eka-stuff' last time.

Given the pattern we're seeing here, and using what
they knew about Postnikov towers, James Dolan and Toby Bartels
guessed a factorization result like this:

\begin{hyp}
  Given a functor $p\maps E\to B$ between $n$-categories, it admits a
  factorization:
\[\xymatrix{
  E_n=E\ar[dr]_{p_{n}} \ar[rrrrr]^p &&&&& B = E_{-2}\\
  & E_{n-1} \ar[r]_{p_{n-1}} &\dots&\dots \ar[r]_{p_0} & E_{-1} 
    \ar[ur]_{p_{-1}} }\]
where $p_j$ forgets purely $j$-stuff.
\end{hyp}

Someone must have already made this precise and proved this
for 2-categories.  If not, someone should
go home and do it tonight; it shouldn't be hard.

In fact, this result is already known for all $n$-groupoids, but only
if you believe Grothendieck that they are the same as homotopy
$n$-types.  In this case, there's a topological result which says that
every map between homotopy $n$-types factors like this.  Such a
factorization is called a \textbf{Moore--Postnikov tower}.  When we
factor a map from a space to a {\it point} this way, it's called a
\textbf{Postnikov tower}.  As we'll see, this lets us view a space as
being made up out of layers, one for each homotopy group.  And this
lets us \emph{classify homotopy types using cohomology} --- at least
in principle.

\subsection{Cohomology and Postnikov towers}
\label{sec:cohomology}

In a minute we'll see that from the viewpoint of homotopy theory,
a space is a kind of `layer cake' with one layer for each dimension.
I claim that cohomology is fundamentally the study of
classifying `layer cakes' like this.  There are many other kinds of
layer cakes, like chain complexes (which are watered-down versions of
spaces), $L_\infty$-algebras and $A_\infty$-algebras (which are chain
complexes with extra bells and whistles), and so on.  But let's start
with spaces.

How does it work?  If $E$ is a homotopy $n$-type, we study it as
follows.  We map it to something incredibly boring, namely a point, 
and then work out the Postnikov tower of this map:

\[\xymatrix{
  E_n=E\ar[dr]_{p_{n}} \ar[rrrrr]^p &&&&& \ast = E_{-2}\\
  & E_{n-1} \ar[r]_{p_{n-1}} &\dots&\dots \ar[r]_{p_0} & E_{-1} 
    \ar[ur]_{p_{-1}} }\]

\noindent
Here we think of $E$ as an $n$-groupoid and the point as the
terminal $n$-groupoid, which has just one object, one morphism
and so on.  The Postnikov tower keeps crushing $E$ down, so
that $E_j$ is really just a $j$-groupoid.  This process is called
{\bf decategorification}.  At the end of the day,
$E_{-2}$ is the only $(-2)$-groupoid there is: the point.

Of course, in the world of topology, they don't use our category-theoretic 
terminology to describe these maps.  Morally, $p_j$
forgets purely $j$-stuff --- but topologists call this `killing the
$j$th homotopy group.'  More precisely, the map
\[p_j \maps E_j \to E_{j-1}\]
induces isomorphisms
\[\pi_i(E_j)\to \pi_i(E_{j-1})\]
for all $i$ \emph{except} $i=j$, in which case it induces the zero map.

We won't get into how topologists actually
construct Postnikov towers.  Once Grothendieck's dream
comes true, it will be a consequence of the result for 
$n$-categories.

It's fun to see how this Postnikov tower works in the shockingly 
low-dimensional cases $j = 0$ and $j = -1$, where the $j$th 
`homotopy thingy' isn't a group --- just a set, or truth value.  
When we get down to $E_0$, our space is just a set, at least
up to homotopy equivalence.  Killing its $0$th `homotopy set' then 
collapses all its points to the same point (if it has any to begin with).  
We're left with either the one-point set or the empty set.
Killing the $(-1)$st `homotopy truth value' then gives us the
one-point set.

But enough of this negative-dimensional madness.  Let's see how
people use Postnikov towers to classify spaces up to homotopy
equivalence.  Consider any simplification step $p_j\maps E_j\to E_{j-1}$ 
in our Postnikov tower.  By
the wonders of homotopy theory, we can describe this as a fibration.
One of the great things about homotopy theory is that even a map
that doesn't look anything like a `bundle' is always 
equivalent to a fibration, so we can think of it as some kind of
bundle-type thing.  Thus we can consider the homotopy fiber
$F_j$ of this map, which can either be constructed directly, the way
we constructed the `essential preimage' for a functor (here using
paths in the base space) or by first converting the map into an actual
fibration and then taking its literal fiber.

(We say `the' fiber as if they were all the same, but if the space
isn't connected they won't necessarily all be the same.  Let's assume
for now that $E$ was connected, for simplicity.)

So, we get a fibration
\[            F_j \to E_j \stackto{p_j} E_{j-1}  \]
Since $p_j$ doesn't mess with any homotopy groups except the $j$th,
the long exact sequence of homotopy groups for a fibration
\[ \dots\too \pi_i(F_j) \too \pi_i(E_j) \too \pi_i(E_{j-1}) 
\too \pi_{i-1}(F_j)\too\dots\]
tells us that the homotopy fiber must have only one non-vanishing
homotopy group: $\pi_i(F_j)=0$ unless $i=j$.  We killed the $j$th
homotopy group, so where did that group go?  It went up into the fiber.

Such an $F_j$, with only one non-vanishing homotopy group, is called
an \textbf{Eilenberg--Mac Lane space}.  The great thing is that a
space with only its $j$th homotopy group nonzero is completely
determined by that group --- up to homotopy equivalence, that is.  For
fancier spaces, the homotopy groups aren't enough to determine the
space: we also have to say how the homotopy groups talk to each other,
which is what this Postnikov business is secretly doing.  The
Eilenberg--Mac Lane space with $G$ as its $j$th homotopy group is
called $K(G,j)$.  (Of course we need $G$ abelian if $j > 1$.)  

So, we've got a way of building any homotopy $n$-type $E$ as a `layer cake' 
where the layers are Eilenberg--Mac Lane spaces, one for each dimension.
At the $j$th stage of this procress, we get a space $E_j$ as the 
total space of this fibration:
\[           K(\pi_j(E), j) = F_j \to E_j \stackto{p_j} E_{j-1} . \]
These spaces $E_j$ become better and better approximations to our space
as $j$ increases, and $E_n = E$. 

If we know the homotopy groups of the space $E$, the main task is
to understand the fibrations $p_j$.  The basic principle of Galois 
theory says how to classify fibrations:

\begin{quote}
\textbf{Fibrations of $n$-groupoids
\[   F\to E\stackto{p} B\]
with a given base $B$ and fiber $F$ are classified by maps
\[   k\colon B\to \AUT(F).\]
}
\end{quote}

Here $\AUT(F)$ is the \textbf{automorphism $(n+1)$-group} of
$F$, i.e.\ an $(n+1)$-groupoid with one object. For example, a set
has an automorphism group, a category has an automorphism 2-group, and
so on.  So, the map $k$ is an $(n+1)$-functor, of the weakest possible sort.

How do we get the map $k$ from the fibration?
Topologists have a trick that involves turning $\AUT(F)$ into a space
called the `classifying space' for $F$-bundles, 
at least when $B$ and $F$ are spaces.  

But now I want you to think about it $n$-categorically.
How does it work?  Think of $B$ as an $\infty$-category.  Then a
path in $B$ (a 1-morphism) lifts to a path in $E$, which when we move
along it, induces some automorphism of the fiber (a 1-morphism in
the one-object $(n+1)$-groupoid $\AUT(F)$).  Similarly, a path of
paths induces a morphism between automorphisms.  This continues all the
way up, which is how we get a map $k\colon B\to\AUT(F)$.

This is a highbrow way of thinking about cohomology
theory.  We may call it `nonabelian cohomology.'  You've probably
seen cohomology with coefficients in some \emph{abelian} group, which
is a special case that's easy to compute; 
here we are talking about a more general version that's supposed to
explain what's really going on.

Specifically, we call the set of $(n+1)$-functors $k \colon B \to \AUT(F)$
modulo equivalence the \textbf{nonabelian cohomology} of $B$ with 
coefficients in $\AUT(F)$.  We denote it like this:
\[H(B,\AUT(F)).\]
We purposely leave off the little superscript $i$'s that people
usually put on cohomology; our point of view is more global.
The element $[k]$ in the cohomology $H(B,\AUT(F))$ corresponding to a
given fibration is called its \textbf{Postnikov invariant}.  

So, to classify a space $E$, we think of it as an $n$-groupoid
and break it down with its Postnikov tower, getting a
whole \emph{list} of guys
\[  k_j \maps E_{j-1} \to \AUT(F_j) \]
and thus a list of Postnikov invariants
\[[k_j]\in H(E_{j-1},\AUT(F_j))  .\]
Together with the homotopy groups of $E$ (which determine 
the fibers $F_j$), these Postnikov
invariants classify the space $E$ up to homotopy equivalence.
Doing this in practice, of course, is terribly hard.  But the 
idea is simple.

Next time we'll examine certain low-dimensional cases of this and see
what it amounts to.  In various watered-down cases
we'll get various famous kinds of cohomology.  The 
full-fledged $n$-categorical version is beyond what anyone
knows how to handle except in low dimensions --- even for $n$-groupoids,
except by appealing to Grothendieck's dream.  Street has a nice
paper on cohomology with coefficients in an $\infty$-category;
he probably knew this stuff I'm talking about way back when I was just a kid.
I'm just trying to bring it to the masses.

\subsection{Questions and comments}
\label{sec:questions-comments-2}

\subsubsection{Internalizing $n$-surjectivity}
\label{sec:internalizing}
\quad

\vskip 1em

Tom Fiore: \emph{You can define epi and mono categorically and apply
them in any category, not just sets.  Can you do a similar thing and
define analogues of 0-, 1-, and 2-surjectivity using diagrams in 
any 2-category, etc.?} 

\vskip 1em

JB: I don't know.  That's a great question. 

\vskip 1em

Eugenia Cheng: \emph{You can define a concept of `essentially epic' in any
2-category, by weakening the usual definition of epimorphism.  But in 
$\Cat$, `essentially epic' turns out to mean 
essentially surjective and full.  I expect that in $n$-categories,
it will give `essentially surjective on $j$-morphisms below level $n$'.
You can't isolate the action on objects from the action on
morphisms, so you can't characterize a property that just refers to
objects.}

\vskip 1em

(For a more thorough discussion see \S\ref{sec:monom-epim}: there is a way to
characterize essentially surjective functors 2-categorically, as the
functors that are left orthogonal to full and faithful functors.)

\subsubsection{How normal people think about this stuff}
\label{sec:normal}
\quad

\vskip 1em

Aaron Lauda: \emph{What do mortals call this $(n+1)$-group
$\AUT(F)$ that you get from an $n$-groupoid $F$?}

\vskip 1em
JB:  Well, suppose we have any $(n+1)$-group, say $G$.  The first
thing to get straight is that there are {\it two ways} to think of
this in terms of topology.

First, by Grothendieck's dream, we can think of $G$ as a 
topological group, say $|G|$, that just happens to be a homotopy 
$n$-type.  Why do the numbers go down one like this?  It's just 
like when people see a category with one object: they call it a 
monoid.   There's a level shift here: the \emph{morphisms} of 
the 1-object category get called \emph{elements} of the monoid.  

Second, we can just admit that our $(n+1)$-group $G$ is
a special sort of $(n+1)$-groupoid.  Following Grothendieck's
dream, we can think of this as a homotopy $(n+1)$-type, called $B|G|$.  
But, we should always think of this as a `connected pointed' 
homotopy $(n+1)$-type.

Why?  Well, an $(n+1)$-group is just an $(n+1)$-groupoid
with one object.
More generally, an $(n+1)$-groupoid is \emph{equivalent} to an $(n+1)$-group
if it's {\bf connected} --- if all its objects are equivalent.  But
to actually turn a connected $(n+1)$-groupoid into an $(n+1)$-group, we 
need to pick a distinguished object, or `basepoint'.  So, an
$(n+1)$-group is essentially the same as a connected pointed 
$(n+1)$-groupoid.  If we translate this into the language of topology, 
we see that an $(n+1)$-group amounts to a connected pointed homotopy 
$(n+1)$-type.  This is usually called $B|G|$, the {\bf classifying space} 
of the topological group $|G|$.  

These two viewpoints are closely related.
The homotopy groups of $B|G|$ are the same as those of $|G|$,
just shifted:
\[           \pi_{j+1}(B|G|) = \pi_j(|G|) . \]

You may think this is unduly complicated.  Why bother thinking about
an $(n+1)$-group in {\it two different ways} using topology?
In fact, both are important.  Given your $n$-groupoid
$F$, it's good to use both tricks just described to study 
the $(n+1)$-group $\AUT(F)$. 

The first trick gives a topological group $|\AUT(F)|$ 
that happens to be a homotopy $n$-type.  This group is often 
called the group of {\bf homotopy self-equivalences} of 
$|F|$, the homotopy $n$-type corresponding to $F$.
The reason is that its elements are homotopy equivalences 
$f\colon |F| \stackrel{\sim}{\longrightarrow} |F|$.

The second trick gives a connected pointed homotopy
$(n+1)$-type $B|\AUT(F)|$.  We can use this to classify fibrations
with $F$ as fiber.

\vskip 1em
Aaron Lauda: \emph{How does that work?}

\vskip 1em
JB: Well, we've seen that fibrations whose base $B$ and fiber $F$
are $n$-groupoids should be classified by $(n+1)$-functors
\[     k \colon B \to \AUT(F)    \]
where $\AUT(F)$ is a $(n+1)$-groupoid with one object.
But normal people think about this using topology.  So, they
turn $B$ into a homotopy $n$-type, say $|B|$. 
They turn $\AUT(F)$ into a connected pointed homotopy 
$(n+1)$-type: the classifying space $B|\AUT(F)|$.  And, they turn
$k$ into a map, say
\[    |k| \colon |B| \to B|\AUT(F)|  .\]
So, instead of thinking of the Postnikov invariant as an
equivalence class of $(n+1)$-functors
\[     [k] \in H(B, \AUT(F))    \]
they think of it as a homotopy class of maps:
\[     [|k|] \in [\,|B|,\, B|\AUT(F)|\,]  \]
where now the square brackets mean `homotopy classes of maps' ---
that's what equivalence classes of $j$-functors become in
the world of topology.
And, they show fibrations with base $|B|$ and fiber $|F|$
are classified by this sort of Postnikov invariant.

Since I'm encouraging you to freely hop back and forth between the 
language of $n$-groupoids and the language of topology, from
now on I won't write $|\cdot |$ to describe the passage from
$n$-groupoids to spaces, or $n$-functors to maps.  I just
wanted to sketch how it worked, here.

\vskip 1em
Aaron Lauda: \emph{So, all this is part of some highbrow 
approach to cohomology... but how does this relate to plain old 
cohomology, like the kind you first learn about in school?}

\vskip 1em
JB:   Right.  Suppose we're playing the Postnikov tower game.  
We have a homotopy $n$-type $E$, and somehow we know 
its homotopy groups $\pi_j$.  So, we get this tower of fibrations
\[           K(\pi_j, j) = F_j \to E_j \stackto{p_j} E_{j-1}  \]
where $E_n$ is the space with started with and $E_{-1}$ is just a
point.  To classify our space $E$ just need to classify
all these fibrations.  That's what the Postnikov invariants do:
\[             [k_j] \in H(E_{j-1},\AUT(K(\pi_j,j)))  .\]
Now I'm using the language of topology, where $\AUT$ stands
for the group of homotopy self-equivalences.
But in the language of topology, the Postnikov invariants are
homotopy classes of maps
\[             k_j \maps E_{j-1} \to B\AUT(K(\pi_j,j)) .  \]

So, in general, our cohomology involves the space 
$B\AUT(K(\pi_j,j))$, which sounds pretty complicated.  
But we happen to know some very nice automorphisms of $K(\pi_j,j)$.
It's an abelian topological group, at least for
$j > 1$, so it can act on itself by left translations.  Thus, sitting inside
$B\AUT(K(\pi_j,j))$ we actually have $BK(\pi_j,j)$, which is
actually the same as $K(\pi_j,j+1)$, since applying $B$ shifts things 
up one level.

If the map $k_j$ happens to land in this smaller space, at least 
up to homotopy, we call our space 
\textbf{simple}.   Then we can write the Postnikov invariant as 
\[          [k_j] \in [E_{j-1} ,K(\pi_j,j+1)]  \]
and the thing on the right is just what people
call the \textbf{ordinary cohomology} of our space $E_{j-1}$ with
coefficients in the group $\pi_j$, at least if $j > 1$.  
So, they write
\[            [k_j] \in H^{j+1}(E_{j-1},\pi_j)  .\]
Note that by now the indices are running all the way from 
$j-1$ to $j+1$, since we've played so many sneaky level-shifting
tricks. 

\vskip 1em
MS: \emph{Actually, Postnikov towers have a nice interpretation 
  in terms of cohomology even for spaces that aren't `simple'.
  The trick is to use `cohomology with local
  coefficients'.  Given a space $X$ and an abelian group $A$
  together with an action $\rho$ of $\pi_1(X)$ on $A$, you can define
  cohomology groups $H^n_\rho(X, A)$ where the coefficients
  are `twisted' by $\rho$.  It then turns out that 
  \begin{align*}
    H(X,\AUT(K(\pi,j))) &= [X, B\AUT(K(\pi,j))]\\
    &= \coprod_{\rho} H^{j+1}_\rho(X,\pi).
  \end{align*}
  So for a space that isn't necessarily simple, a topologist would
  consider its Postnikov invariants to live in some cohomology with
  local coefficients.  In the simple case, the action $\rho$ is trivial,
  so we don't need local coefficients.
}

\section{A Low-Dimensional Example}
\label{sec:a-low-dimensional}

\subsection{Review of Postnikov towers}
\label{sec:review-cohom}

Last time we discussed a big idea; this time let's look at an example.  
Let's start with a single fibration:
\[ F\to E\to B. \]
This means that we have some point $\ast\in B$ and $F = p\inv(\ast)$
is the homotopy fiber, or `essential preimage' over $\ast$.  This
won't depend on the choice of $\ast$ if $B$ is connected.
Let's restrict ourselves to that case: this is no great loss,
since any base is a disjoint union of connected components.

We can then classify these fibrations via their `classifying maps'
\[ k\maps B\to \AUT(F) \]
where $\AUT(F)$ is an $(n+1)$-group if $F$
is an $n$-groupoid.  A lowbrow way to state this classification is
that there's a notion of equivalence for both these guys, and the
equivalence classes of each are in one-to-one correspondence.  We
could also try to state a highbrow version, which asserts that 
$\hom(B,\AUT(F))$ is equivalent to some $(n+1)$-category of fibrations
with $B$ as base and $F$ as fiber.   But let's be lowbrow today.

In both cases, hopefully the notion of equivalence is sort of obvious.
`Equivalence of fibrations' looks a lot like equivalence
for extensions of groups --- which are, in fact, a special case.  In other
words, fibrations are equivalent when there exists a vertical map 
making this diagram commute (weakly):
\[
\xymatrix{
  & E \ar[dr] \ar[dd]^\eqv\\
  F\ar[ur]\ar[dr]  && B\\
  & E' \ar[ur]
}\]
On the other side, the notion of equivalence for classifying maps 
\[
k,k'\maps B\to \AUT(F)
\]
is equivalence of $(n+1)$-functors, if we think of $E$ and $F$
as $n$-groupoids, or homotopy of maps if we think of $E$
and $F$ spaces.

This sounds great, but of course we're using all sorts of concepts
from $n$-category theory that haven't been made precise yet.  So, 
today we'll do an example, where we cut things down to a low enough 
level that we can handle it.

But first --- why are we so interested in this?  I hope you remember
why it's so important.  We have a grand goal: we want to classify
$n$-groupoids.   This is a Sisyphean task.  We'll never
actually complete it --- but nonetheless, we can learn a lot by trying.

For example, consider the case $n = 1$: the classification of groupoids.  
Every groupoid is a disjoint union of groups, so we just need to 
classify groups.  Let's say we start by trying to classify finite
groups.  Well, it's not so easy --- after 10,000 pages of work people
have only managed to classify the finite \emph{simple} groups.
Every finite group can be built up out of those by repeated extensions:
\[ 1\to F\to E\to B \to 1.  \]
These extensions are just a special case of the fibrations we've been talking
about.   So we can classify them using cohomology, at least in principle.  
But it won't be easy, because each time we do an extension we get a 
new group whose cohomology we need to understand.  So, we'll
probably never succeed in giving a useful classification of all
finite groups.  Luckily, even what we've learned so far can help us 
solve a lot of interesting problems. 

Now suppose we want to classify $n$-groupoids for $n > 1$.  
We do it via their Postnikov towers, which are iterated 
fibrations.  Given an $n$-groupoid $E$, we successively 
squash it down, dimension by dimension, until we get a single point: 
\[\xymatrix{
  E_n=E\ar[dr]_{p_{n}} \ar[rrrrr]^p &&&&& \ast = E_{-2}\\
  & E_{n-1} \ar[r]_{p_{n-1}} &\dots&\dots \ar[r]_{p_0} & E_{-1} 
  \ar[ur]_{p_{-1}}
}\]

\noindent
At each step, we're \textbf{decategorifying}: to get the $(j-1)$-groupoid 
$E_{j-1}$ from the $j$-groupoid $E_j$, we promote all the 
$j$-isomorphisms to \emph{equations}.  That's really what's going on, 
although I didn't emphasize it last time.  Last time I emphasized 
that the map 
\[     p_j \maps E_j \to E_{j-1}   \]
`forgets purely $j$-stuff'.  What that means here is that decategorification
throws out the top level, the $j$-morphisms, while doing as little 
damage as possible to the lower levels.  

So, we get fibrations
\[         F_j \to E_j \to E_{j-1}  \]
where each homotopy fiber $F_j$, which records the stuff that's been 
thrown out, is a $j$-groupoid with only nontrivial $j$-morphisms: it 
has at most one $i$-morphism for any $i < j$.  
If you look at the periodic table, you'll see this means $F_j$ is
secretly a \emph{group} for $j = 1$, and an \emph{abelian group} for 
$j \ge 2$.  Another name for this group is $\pi_j(E)$, which is
more intuitive if you think of $E$ as a space.  If you think of 
$F_j$ as a space, then it's the Eilenberg--Mac Lane space 
$K(\pi_j(E), j)$.

What do we learn from this business?  That's where the basic
principle of Galois theory comes in handy.  We take the fibrations
\[F_j \to E_j\too[p_j] E_{j-1}\]
and describe them via their classifying maps
\[ k_j\maps E_{j-1} \to \AUT(F_j). \]
These give cohomology classes
\[   [k_j] \in H(E_{j-1}, \AUT(F_j))  \]
called {\bf Postnikov invariants}.

So, we ultimately classify $n$-groupoids by a list of groups, namely
$\pi_1$, \dots , $\pi_n$, and all these cohomology classes $[k_j]$.
What I want to do is show you how this works in detail, in a very
low-dimensional case. 

\subsection{Example: the classification of 2-groupoids}
\label{sec:n=2-case}

Let's illustrate this for $n=2$ and classify connected 2-groupoids.
Since we're assuming things are connected, we might as well, for the
purposes of classification, consider our connected 2-groupoids to be
2-groups.  These have one object, a bunch of 1-morphisms from it to
itself which are weakly invertible, and a bunch of 2-morphisms from
these to themselves which are strictly invertible.

We can classify these using cohomology.  
Here's how.  Given a 2-group, take a skeletal version of it, say
$E$, and form these four things:

{\bf 1.} The group $G = \pi_1(E) = $ the group of `1-loops', i.e.\
  1-morphisms that start and end at the unique object.  Composition of
  these would, a priori, only be associative up to isomorphism, but we
  said we picked a \emph{skeletal} version, so these isomorphic
  objects have to be, in fact, equal.

{\bf 2.}  The group $A = \pi_2(E) =$ the group of `2-loops', i.e.\
  2-morphisms which start and end at the identity 1-morphism
  $1_\ast$.  They form a group more obviously, and the Eckmann--Hilton
  argument shows this group is abelian.

{\bf 3.}  An action $\rho$ of $G$ on $A$, where $\rho(g)(a)$ is defined
  by `conjugation' or `whiskering':
  \[\xymatrix{ \bullet \ar@/_5mm/[rrr] \ar@/^5mm/[rrr] 
   \rrtwocell\omit{\makebox[0pt][l]{$\scriptstyle\rho(g)(a)$}} &  && \bullet} 
   = 
   \xymatrix{
   \bullet \ar[r]^g & 
   \bullet \rtwocell^{1_\ast}_{1_\ast}{a} & 
   \bullet \ar[r]^{g\inv} & \bullet
  }\]
  You can think of the loops as starting and ending at anything, if
  you like, by doing more whiskering.  Then you have to spend a year
  figuring out whether you want to use left whiskering or right
  whiskering.  This is supposed to be familiar from topology: there
  $\pi_1$ always acts on $\pi_2$.

{\bf 4.} The associator
  \[ \alpha_{g_1g_2g_3}\maps (g_1g_2)g_3) \to g_1(g_2g_3) \]
  gives a map
  \[ \alpha\maps G^3\to A \]
  as follows.  Take three group elements and get an interesting
  automorphism of $g_1g_2g_3$.  Automorphisms of anything can be
  identified with automorphisms of the identity, by whiskering.
  Explicitly, we cook up an element of $A$ as follows:
  \[\xymatrix{
   \bullet \ar[rr]^{(g_1g_2g_3)\inv} && \bullet 
   \rrtwocell^{g_1g_2g_3}_{g_1g_2g_3}{\alpha}&& \bullet
  }\]
  (Don't ask why I put the whisker on the left instead of the
  right; you can do it either way and it doesn't really matter,
  though various formulas work out slightly differently.)

That's the stuff and structure, but there's also a property: the
  associator satisfies the pentagon identity, which means that
  $\alpha$ satisfies some equation.  You all know the pentagon
  identity.  It turns out this equation on $\alpha$ is something
  that people had been talking about for ages, before Mac Lane
  invented the pentagon identity.  In fact, one of the people who'd been
  talking about it for ages was Mac Lane himself, because he'd also
  helped invent cohomology of groups.  It's called the 
  \textbf{3-cocycle equation}
  in group cohomology:  
   \[
  \rho(g_0)\alpha(g_1,g_2,g_3) - \alpha(g_0g_1, g_2, g_3) +  
  \alpha(g_0, g_1g_2, g_3) \\ - \alpha(g_0, g_1, g_2g_3) +
  \alpha(g_0, g_1, g_2) = 0
  \]
  Here we write the group operation in $A$ additively, since it's abelian.

  This equation is secretly just the pentagon identity satisfied by the
  associator; that's why it has 5 terms.  But, people in group cohomology
  often write it simply as $d\alpha=0$, because they know a standard
  trick for getting function of $(n+1)$ elements of $G$ from a function 
  of $n$ elements, and this trick is called the `differential' $d$ in 
  group cohomology.  If you have trouble remembering this trick, just 
  think of a bunch of kids riding a school bus, but today there's one
  more kid than seats on the bus.  What can we do?  Either the first kid
  can jump out and sit on the hood, or the first two kids can squash 
  into the first seat, and so on... or you can throw the last kid out 
  the back window.  That's a good way to remember the formula I just wrote.

Note that our \emph{skeletal} 2-group is not necessarily
\emph{strict}!  Making isomorphic objects equal doesn't mean making
isomorphisms into identities.  The associator isomorphism is still
nontrivial, but it just happens in this case to be an
\emph{automorphism} from one object to the same object.

\begin{thm}[Sinh]
  Equivalence classes of 2-groups are in one-to-one correspondence
  with equivalence classes of 4-tuples
  \[(G, A, \rho, \alpha)\]
  consisting of a group, an abelian group, an action, and a 3-cocycle.
\end{thm}

The equivalence relation on the 4-tuples is via isomorphisms of $G$ and
$A$ which get along with $\rho$, and get along with $\alpha$ up to a
coboundary.  Since cohomology is precisely cocycles modulo
coboundaries, we really get the traditional notion of group cohomology
showing up.  

Just as the 3-cocycle equation comes from the pentagon 
identity for monoidal categories, the 
coboundary business comes from the notion of a monoidal \emph{natural 
transformation}: monoidally equivalent monoidal categories won't have 
the same $\alpha$, but their $\alpha$'s will differ by a coboundary.

This particular kind of group cohomology is called the `third'
cohomology since $\alpha$ is a function of three variables.
We say that
\[
[\alpha] \in H^3_\rho(G,A),
\]
the third group cohomology of $G$ with coefficients in a $G$-module
$A$ (where the action is defined by $\rho$).  Sometimes this is called
`twisted' cohomology.  

So, in short, once we fix $G$, $A$, and $\rho$, the equivalence classes 
of 2-groups we can build are in one-to-one correspondence with 
$H^3_\rho(G,A)$. 

\subsection{Relation to the general case}
\label{sec:relat-gener-case}

Now I'm going to show why this stuff is a special case of the general
notion of cohomology we introduced last time.  Why is $H^3_\rho(G,A)$
a special case of what we were calling $H(B,\AUT(F))$?

Consider a little Postnikov tower where we start with a 2-group
$E$ and decategorify it getting a group $B$.  We get a fibration 
$F\to E\to B$.  To relate this to what we were just talking about,
think of $B$ as the group $G$.  The 2-group $F$ is the 2-category with
one object, one morphism and some abelian group $A$ as 2-morphisms.
So, seeing how the 2-group $E$ is built out of the base $B$ and the
fiber $F$ should be the same as seeing how its built out of $G$ and
$A$.  We want to see how the classifying map
\[k\maps B \to \AUT(F)\]
is the same as an element of $H^3_\rho(G,A)$.

It's good to think about this using a little topology.  
As a \emph{space}, $B$ is called
$K(G,1)$.  It is made by taking one point:
\[         \bullet    \]
one edge for element $g_1 \in G$:
\[   \xy 0;/r.25pc/:
  (-8,0)*+{\bullet}="1";
  (0,5)*+{}; 
  (0,-5)*+{}; 
  (8,0)*+{\bullet}="2";
  {\ar "1";"2"^{g}};
 \endxy
\]
a triangle for each pair of elements:
\[   \vcenter{\xy 0;/r.25pc/:
   (-10,-5)*+{ \scriptstyle \bullet}="1";
   (10,-5)*+{ \scriptstyle \bullet}="2";
   (0,12)*+{ \scriptstyle \bullet}="3";
    {\ar "1"; "2"_{g_1 g_2}};
    {\ar "3"; "2"^{g_2}};
    {\ar "1"; "3"^{g_1}};
    \endxy}
\]
a tetrahedron for each triple:
\[ \xy 0;/r.30pc/:
    (-10,-5 )*+{ \scriptstyle \bullet}="1";
    (8,-10)*+{ \scriptstyle \bullet}="2";
    (15,0)*+{ \scriptstyle \bullet}="3";
   (1,12)*+{ \scriptstyle \bullet}="4";
       {\ar "1";"2"_{g_1g_2} };
       {\ar "2";"3"_{g_3} };
       {\ar "4";"3"^{g_2 g_3} };
    {\ar "1";"4"^{g_1} };
    {\ar "4";"2"_<<<<<<<<{g_2} };
       {\ar @{.>} "1";"3"^<<<<<<<<<{g_1 g_2 g_3}};
       \endxy
\]
and so on.  The fiber $F$ is $A$ regarded as a 2-group with only an identity
1-cell and 0-cell.  So, as a space, $F$ is called $K(A,2)$, built with one
point, one edge, one triangle for every element of $A$, one
tetrahedron whenever $a_1+a_2=a_3+a_4$, and so on.

Now, $\AUT(F)$ is a \emph{3-group} which looks roughly like this.  
It has one object (which we can think of as `being' $F$), one
1-morphism for every automorphism $f\maps F\to F$, one 2-morphism
\[\xymatrix{
  F\rtwocell^{f}_{f'}{\gamma} & F
}\]
for each pseudonatural isomorphism, and one 3-morphism for each
modification.

That seems sort of scary; what you have to do is figure out what that
actually amounts to in the case when $F$ is as above.  Let me just
tell you.  It turns out that in fact, in our case $\AUT(F)$ has
\begin{itemize}
\item one object;
\item its morphisms are just the group $\Aut(A)$;
\item only identity 2-morphisms (as you can check);
\item $A$ as the endo-3-morphisms of any 2-morphism.
\end{itemize}
As a space, this is called `$B(\AUT(F))$', and it is made from one
point, an edge for each automorphism $f$, a triangle for each equation
$f_1f_2=f_3$, and tetrahedrons whose boundaries commute which are
labeled by arbitrary elements of $A$.

Now, we can think about our classifying map as a weak 3-functor
\[k\maps B \to \AUT(F)\]
but we can also think of it as a map of spaces
\[k\maps K(G,1) \to B(\AUT(F))\]
Let's just do it using spaces --- or actually, simplicial sets.  We have
to map each type of simplex to a corresponding type.  Here's how it
goes:
\begin{itemize}
\item This map is boring on 0-cells, since there's only one choice.
\item We get a map $\rho\maps G\to \Aut(A)$ for the 1-cells, which is good
  because that's what we want.
\item The map on 2-cells says that this a group homomorphism, since it
  sends equations $g_1g_2=g_3$ to equations
  $\rho(g_1)\rho(g_2)=\rho(g_3)$.
\item The map on 3-cells sends tetrahedra in $K(G,1)$, which are
  determined by triples of elements of $G$, to elements of $A$.  This
  gives the may $\alpha\maps G^3\to A$.
\item The map on 4-cells is what forces $\alpha$ to be a 3-cocycle.
\end{itemize}

Our map $B\to \AUT(F)$ is \emph{weak}, which is all-important.  Even
though our 2-morphisms are trivial, which makes the action on
1-morphisms actually a strict homomorphism, and our domain has no
interesting 3-morphisms, we also get the higher data which gives
the 3-cocycle $G^3\to A$:
\[\xymatrix{
  \ast \ar[r] & \ast\\
  \ast \ar[r] & A\\
  \ast \ar[r] & \ast\\
  G \ar[r]\ar[ruu] & \Aut(A)\\
  \ast \ar[r] & \ast  
}\]

In fact, all sorts of categorified algebraic
gadgets should be classified as `layer cakes' built using
Postnikov invariant taking values in the cohomology theory for that sort
of gadget.  We get group cohomology when we 
classify $n$-groupoids.   Similarly, to classify categorified Lie
algebras, which are called $L_\infty$ algebras, we need Lie algebra 
cohomology --- Alissa Crans has checked this in the simple case of an
$L_\infty$ algebra with only two nonzero chain groups.  The next case is 
$A_\infty$ algebras, which are categorified associative algebras.
I bet that classifying these involves Hochschild cohomology --- 
but I haven't ever sat down and checked it.  And, it should keep on
going.  There should be a general theorem about this.  That's
what I mean by the `layer cake philosophy' of cohomology.

\subsection{Questions and comments}
\label{sec:questions-comments3}

\subsubsection{Other values of $n$}
\label{sec:other-values-n}  \quad

\vskip 1em
Aaron Lauda: \emph{what about other $H^n$?}

\vskip 1em
JB: Imagine an alternate history of the world in which people knew about
$n$-categories and had to learn about group cohomology from them.

We can figure out $H^3$ using 2-groups.
\[
\begin{array}{c}
  \vdots \\ \ast\\ \pi_2 \\ \pi_1 \\ \ast
\end{array}
\]
To get the classical notion of $H^4$, we would have to think about
classifying 3-groups that only have interesting 3-morphisms and
1-morphisms.
\[
\begin{array}{c}
  \vdots \\\ast \\ \pi_3 \\ \ast \\ \pi_1 \\ \ast
\end{array}
\]
We still get the whiskering action of 1-morphisms on 3-morphisms, and
the pentagonator gives a 4-cocycle.  Group cohomology, as
customarily taught, is about classifying these `fairly wimpy'
Postnikov towers in which there are just two nontrivial groups:
$H^{n+2}(\pi_1,\pi_n)$.  This is clearly just a special case of
something, and that something is a lot more complicated.

\vskip 1em
MS: \emph{What about $H^2$?}

\vskip 1em
JB:  Ah, that's interesting!
In general, $H^n(G,A)$ classifies ways of building an $(n-1)$-group
with $G$ as its bottom layer (what I was just calling $\pi_1$)
and $A$ as its top layer (namely $\pi_{n-1}$).  So, it's all about
layer-cakes with two nontrivial layers: the first and $(n-1)$st layers. 
For simplicity I'm assuming the action of $G$ on $A$ is trivial here.

But the case $n = 2$ is sort of degenerate: now our layer cake 
has only {\it one} nontrivial layer, the first layer, built by 
squashing $A$ right into $G$.  More precisely, $H^2(G,A)$ classifies ways 
of building a 1-group --- an ordinary group --- by taking a {\it 
central extension} of $G$ by $A$.  We've seen that 3-cocycles come from 
the associator, so it shouldn't be surprising that 2-cocycles 
come from something more basic: multiplication, where \emph{two} 
elements of $G$ give you an element of $A$.  

Ironically, this weird degenerate low-dimensional case is the
highest-dimensional case of group cohomology that ordinary textbooks
bother to give any clean conceptual interpretation to.  They say that $H^2$ 
classifies certain ways build a group out of two groups, but they
don't say that $H^n$ classifies certain ways to build an $(n-1)$-group
out of two groups.  They don't say that extensions are just degenerate
layer-cakes.  And, it gets even more confusing when people start
using $H^2$ to classify `deformations' of algebraic structures, 
because they don't always admit that a deformation is just a special 
kind of extension.

\subsubsection{The unit isomorphisms}
\label{sec:unit-isomorphisms} \quad

\vskip 1em
MS: \emph{What happened to the unit isomorphisms?}

\vskip 1em
JB: When you make a 2-group skeletal, you can also make its unit
isomorphisms equal to the identity.  However, you can't make the
associator be the identity --- 
that's why it gives some interesting data in our classification
of 2-groups, namely the 3-cocycle.  This 3-cocycle is the only
obstruction to a 2-group being both skeletal and strict.

\section{Appendix: Posets, Fibers, and $n$-Topoi}
\label{sec:enrichm-posets-fiber}

This appendix is a hodgepodge of (proposed) answers to questions that
arose during the lectures, some musings about higher topos theory, and
some philosophy about the distinction between pointedness and
connectedness.

\subsection{Enrichment and posets}
\label{sec:enrichment-posets}

We observed in \S\ref{sec:maps-truth-values} that while
$(-1)$-categories are truth values, having only a $0$-category (that
is, a set) of them is a bit unsatisfactory, since it doesn't allow us
to talk about \emph{implication} between truth values.  The notetaker
believes the best resolution to this problem is to extend the notion
of `$0$-category' to include not just sets but \emph{posets}.  Then we
can say that truth values form a poset with two elements, true and
false, and a single nonidentity implication
$\mathrm{false}\Impl\mathrm{true}$.  A set can then be regarded as a
discrete poset, or equivalently a poset in which every morphism is
invertible; that is, a `$0$-groupoid'.

One way to approach a general theory including posets is to start from
very low-dimensional categories and build up higher-dimensional ones
using \emph{enrichment}.  We said in \S\ref{sec:power-neg-thinking}
that one of the general principles of $n$-category theory is that big
$n$-categories have lots of little $n$-categories inside them.
Another way of expressing this intuition is to say that \emph{An
  $n$-category is a category enriched over $(n-1)$-categories}.

What does `enriched' mean?  Roughly speaking, a \textbf{category
  enriched over $V$} consists of
\begin{itemize}
\item A collection of objects $x,y,z,\dots$;
\item For each pair of objects $x,y$, an object in $V$ called $\hom(x,y)$;
\item For each triple of objects $x,y,z$, a morphism in $V$ called
  \[\circ \maps \hom(x,y)\times \hom(y,z)\too\hom(x,z)\]
\item Units, associativity, etc.\ etc.\
\end{itemize}
Remember than in the world of $n$-categories, it doesn't make sense to
talk about anything being strictly equal except at the top level.  So
when $V$ is, for instance, the $(n+1)$-category of $n$-categories, the
composition in a $V$-enriched category should only be associative and
unital up to coherent equivalence.  For example:
\begin{itemize}
\item 1-categories are categories enriched over sets;
\item weak 2-categories are categories weakly enriched over
  categories;
\item weak 3-categories are categories weakly enriched over
  weak 2-categories.
\end{itemize}
Making this precise for $n>2$ is tricky, but it's a good intuition.

We may also say that a \textbf{$V$-enriched functor} $p\maps E\to B$
between such categories consists of:
\begin{itemize}
\item A function sending objects of $E$ to objects of $B$;
\item Morphisms of $V$-objects $\hom(x,y)\to \hom(px,py)$;
\item various other data (again, as weak as appropriate).
\end{itemize}
And that's as far up as we need to go.  We'll say that a
\textbf{$V$-enriched groupoid} is a $V$-enriched category such that
`every morphism is invertible' in a suitably weak sense.

Let's investigate this notion in our very-low-dimensional world,
starting with $(-2)$-categories, which we take to all be trivial by
definition.  What is a category enriched over $(-2)$-categories?
Well, it has a collection of objects, together with, for every two
objects, the unique $(-2)$-category as $\hom(x,y)$, and composition
maps which are likewise unique.  Thus a $(-2)$-category-enriched
category is either:
\begin{itemize}
\item empty (has no objects), or
\item has some number of objects, each of which is uniquely isomorphic
  to every other.
\end{itemize}
Thus it is either empty (false) or contractible (true), agreeing with
the notion of $(-1)$-category that we got from topology.  In
particular, every $(-1)$-category is a groupoid.

Our general notion of functor now says that there should be a
$(-1)$-functor from false to true, which we can call `implication'.
This is in line with topology: there is also a continuous map from the
empty space to a contractible one.

Continuing on, a category enriched over $(-1)$-categories has a
collection of objects together with, for every pair of objects $x,y$,
a truth value $\hom(x,y)$, and for every triple $x,y,z$, a morphism
\[
\hom(x,y) \times \hom(y,z) \too \hom(x,z)
\]
Now, the product of two $(-1)$-categories is empty (that is,
false) if and only if one of the factors is.  Thus, when we interpret
$(-1)$-categories as truth values, the product $\times$ becomes the
logical operation `and', so if we interpret the truth of $\hom(x,y)$
as meaning $x\le y$, we see that a category enriched over
$(-1)$-categories is precisely a \textbf{poset}.  (A non-category
theorist would call this a \emph{preordered set} since we don't have
antisymmetry, but from a category theorist's perspective that's asking
for equality of objects instead of isomorphism, which is perverse.)

Thus, from the enrichment point of view, perhaps `0-category' should
mean a poset, rather than a set.  As remarked above, we can view a set
as a discrete poset (that is, one in which $x\le y$ only when $x=y$).
Categorically, a poset is \emph{equivalent} to a discrete one whenever
$x\le y$ implies $y\le x$, which is essentially the condition for it
to be a \emph{0-groupoid}: since composites are uniquely determined in
a poset, a morphism $x\le y$ is an isomorphism precisely when $y\le x$
as well.

In a way, it's not surprising that our intuition may have been a
little off in this regard, since a lot of it was coming from topology
where \emph{everything} is a groupoid.

What happens at the next level?  Well, if we enrich over sets
(0-groupoids), we get what are usually called categories, or
1-categories.  If we enrich instead over posets, we get what could
variously be called \textbf{poset-enriched categories}, \textbf{locally
  posetal 2-categories}, or perhaps \textbf{2-posets}.

At the 0-level, we had one extra notion arising: instead of just sets,
we got posets as well.  At the 1-level, in addition to categories and
poset-enriched categories, we also have a third notion: groupoids.
These different levels correspond to the different levels of
invertibility one can impose.  If we start with a 2-poset and make its
2-morphisms all invertible, we get just a category.  Then if we go
ahead and make its 1-morphisms also invertible, we end up with a
groupoid.

We can then go ahead and consider categories enriched over each of
these three things, obtaining respectively 3-posets, 2-categories, and
locally groupoidal 2-categories.  And again there is an extra level
that comes in: if we also make the 1-morphisms invertible, we get
2-groupoids.

All these various levels of invertibility can be fit together into the
`enrichment table' below.  A dotted arrow $\xymatrix{X\ar@{.>}[r]&Y}$
means that a $Y$ is a category enriched over $X$s.  The horizontal
arrows denote inclusions; as we move to the left along any given line,
we make more and more levels of morphisms invertible, coming from the
top down.  In general, the $n$th level will have $n+2$ different
levels of invertibility stretching off to the left.

\vskip 1em
\vbox{
\begin{center}
{\bf THE ENRICHMENT TABLE }
  \[\small\xymatrix@=1pc{
    \ddots & \vdots & \vdots & \vdots & \vdots \\
    &\text{2-groupoids} \ar@{^(->}[r] \ar@{.>}[u] & 
    \text{\parbox{1in}{locally groupoidal\\2-categories}} \ar@{^(->}[r]  \ar@{.>}[u] &
    \text{2-categories} \ar@{^(->}[r]  \ar@{.>}[u] &
    \text{3-posets}  \ar@{.>}[u]\\
    && \text{groupoids} \ar@{^(->}[r] \ar@{.>}[u] & \text{categories}
    \ar@{^(->}[r] \ar@{.>}[u] &
    \text{2-posets} \ar@{.>}[u]\\
    &&& \text{sets} \ar@{^(->}[r] \ar@{.>}[u] &
    \text{posets} \ar@{.>}[u]\\
    &&&&
    \text{truth values} \ar@{.>}[u]\\
    &&&& \text{triviality} \ar@{.>}[u] }\]\centering
\end{center}
}
\vskip 1em

What does the enrichment table have to do with the periodic table?
Recall that the $n$ in `$n$-categories' labels the \emph{columns} of
the periodic table, while the rows are labeled with the amount of
monoidal structure.  Thus we could, if we wanted to, combine the two
into a three-dimensional table, replacing the line across the top of
the periodic table with the whole table of enrichment.

The enrichment table only includes $n$-categories for finite $n$, but
we can obtain various different `$\infty$' notions by thinking about
passing to some sort of `limit' in various directions.  Of course,
these aren't actually limits in any formal sense.  For example, it
makes intuitive sense to say that the `vertical limit' along the
column
\[\xymatrix{\text{sets} \ar@{.>}[r] &
  \text{1-categories} \ar@{.>}[r] &
  \text{2-categories} \ar@{.>}[r] &
  \dots}
\]
should be the $\infty$-categories.  Moreover, this should also be the
limit along any other column.  This is because in (say) the $m$th
column, all the cells of the top $m$ dimensions are invertible, but in
the limit all these invertible cells get pushed off to infinity and we
end up with noninvertible cells of all dimensions.

We can also consider `diagonal limits'.  It makes intuitive sense to
say that the limit along the far-left diagonal, consisting of
$n$-groupoids for increasing $n$, is the $\infty$-groupoids, aka
homotopy types (\`a la Grothendieck).  The limit along the next diagonal
will be the $\infty$-categories with all morphisms above level 1
invertible.  These are often called $(\infty,1)$-categories (but
sometimes also $(1,\infty)$-categories); see Bergner's survey article
for an introduction to them.

By the way, the term `$(\infty,1)$-categories' may sound strange, but
it is just the most frequently used case of a general terminology.  An
\textbf{(\textit{n,m})-category} is an $n$-category all of whose
$j$-morphisms for $j>m$ are invertible.  Thus a $n$-category may also
be called an $(n,n)$-category, an $n$-groupoid may be called an
$(n,0)$-category, and a locally groupoidal 2-category may be called a
$(2,1)$-category.

To stretch this terminology to its logical limit, we can call a
poset-enriched category a $(1,2)$-category, a poset a
$(0,1)$-category, and so on for the right-hand column of the
enrichment table.  If, instead of regarding an $n$-category as
enriched over $(n-1)$-categories, we return to regarding it as an
$\infty$-category in which all cells of dimension $>n$ are identities,
we can give the following characterization of $(n,m)$-categories which
includes the case of posets as well.

\begin{defn}
  An \textbf{(\textit{n},\textit{m})-category} is an $\infty$-category 
  such that
  \begin{itemize}
  \item All $j$-morphisms for $j>n+1$ exist and are unique wherever
    possible.  In particular, this implies that all parallel
    $(n+1)$-morphisms are `equal'.
  \item All $j$-morphisms for $j>m$ are invertible.
  \end{itemize}
\end{defn}

In the next section we'll consider at length one reason that including
the $n$-posets in the periodic table is important.  Here's a
different, simpler reason.  Let $E$ be a category, and consider its
Postnikov tower:
\[\xymatrix{
E_1=E\ar[dr]_{\mathrm{0, 1-surj}}\ar[rrrr]^p &&&& \ast=B=E_{-2}\\
& E_0 \ar[rr]_{\mathrm{0, 2-surj}} && E_{-1} \ar[ur]_{\mathrm{1, 2-surj}}
}\]
As we said in \S\ref{sec:factorizations}, $E_0$ is what we get by
making parallel morphisms in $E$ equal if they become equal in $B$;
but here $B=\ast$, so this just means we identify \emph{all} parallel
morphisms.  This precisely makes $E$ into a \emph{poset}---not
necessarily a set.  Thus in order for $E_j$ to be a $j$-category in
the factorization of an $n$-category which isn't a groupoid, we have
to consider posets as a sort of 0-category, poset-enriched categories
as a sort of 1-category, and generally $j$-posets as a sort of
$j$-category.

\subsection{Fibers and fibrations}
\label{sec:fibers}

Consider the fibers (or, rather, homotopy fibers) of a functor $p\maps
E\to B$; we saw in \S\ref{sec:stuff-struct-prop} that their
`dimension' should reflect how much the functor $p$ forgets.  We'd
like a generalization of Fact~\ref{forgetfulness-by-fibers} there that
applies to categories in addition to groupoids, but it turns out that
for this we'll need to include the $n$-posets again.  Consider first
the following examples.

\begin{eg}
  We know that the functor $p\maps  \mathbf{AbGp}\to \mathbf{Gp}$ forgets
  only properties.  What is the (essential) preimage $p\inv(G)$ for
  some group $G$?  It is the category of all abelian groups equipped
  with isomorphisms to $G$, and morphisms which preserve the given
  isomorphisms.  This category is contractible if $G$ is abelian, and
  empty otherwise; in other words, it is essentially a $(-1)$-category.
\end{eg}

\begin{eg}
  Even more simply, consider an equivalence of categories $p\maps E\to B$,
  which forgets nothing.  The the preimage $p\inv(b)$ is nonempty
  (since $p$ is essentially surjective), and contractible (since $p$
  is full and faithful); thus it is essentially a $(-2)$-category.
\end{eg}

These examples, along with the groupoid case we considered in
Fact~\ref{forgetfulness-by-fibers}, lead us to guess that a functor
will forget `at most $n$-stuff' precisely when its essential preimages
are all $n$-categories.  We consider properties to be $(-1)$-stuff,
structure to be 0-stuff, ordinary stuff to be 1-stuff, eka-stuff to be
2-stuff, and so on.

However, this guess is not quite right, as we can see by considering
some examples that forget structure.

\begin{eg}
  Consider the usual forgetful functor $p\maps
  \mathbf{Gp}\to\mathbf{Set}$, which we know forgets at most
  structure.  Given a set, such as the 4-element set, its essential
  preimage $p\inv(4)$ is the category of 4-element \emph{labeled}
  groups (since their underlying sets are equipped with isomorphisms
  to the given set $4$), and homomorphisms that preserve the labeling.

  What does this look like?  Well, given two labeled 4-element groups,
  there's exactly one function between them that preserves the
  labeling and either it's a group homomorphism or it isn't.  Since
  the function preserving labeling is necessarily a bijection, if it
  is a homomorphism, then it is in fact a group isomorphism; thus this
  category is (equivalent to) a set.
\end{eg}

In this example, we got what we expected, but we had to use a
special property of groups: that a bijective homomorphism is an
isomorphism.  For many other types of structure, this won't be the
case.

\begin{eg}
  Consider the forgetful functor $p\maps \mathbf{Top}\to\Set$ sending
  a topological space to its underlying set of points, which also
  forgets at most structure (in fact, purely structure).  In this
  case, the essential preimage of the 4-element set is the collection
  of labeled 4-point topological spaces and continuous maps that
  preserve the labeling.  Again, between any two there is exactly one
  function preserving the labeling, and either it is continuous or it
  isn't, so this category is a \emph{poset}.  In general, however, it
  won't be a set, since a continuous bijection is not necessarily a
  homeomorphism.
\end{eg}

Thus, in order to get a good characterization of levels of
forgetfulness by using essential preimages, we really need to include
the $n$-posets as $n$-categories.

Let's look at a couple of examples involving higher dimensions.

\begin{eg}
  We have a forgetful 2-functor
  \[[\text{monoidal categories}]\too{} [\text{categories}]\]
  which forgets at most stuff (since it is locally faithful, i.e.\
  3-surjective).  Here the fiber over a category $C$ is the
  category of ways to add a monoidal structure to $C$.  There are
  lots of different ways to do this, and in between them we have
  monoidal functors that are the identity on objects (up to a
  specified equivalence, if we use the essential preimage), and in
  between \emph{those} we have monoidal transformations whose
  components are identities (or specified isomorphisms).  Now, there's
  at most one natural transformation from one functor to another whose
  components are identities, and either it's monoidal or it isn't.
  This shows that this collection is in fact a locally posetal
  2-category, or a `2-poset', but in fact these monoidal natural
  transformations are automatically invertible when they exist, so it
  is in fact it is a 1-category.
\end{eg}

\begin{eg}
  Let \V\ be a nice category to enrich over, and consider the
  `underlying ordinary category' functor
  \[(-)_0\maps \VCat \too \Cat.\]
  The category $C_0$ has the same objects as C, and $C_0(X,Y) =
  \V(I,C(X,Y))$.  What this functor forgets depends a lot on \V:
  \begin{itemize}
  \item In many cases, such as topological spaces, simplicial sets,
    categories, it is 0-surjective (any ordinary category can be
    enriched), but in others, such as abelian groups, it is not.
  \item In general it is not 1-surjective: not every ordinary functor
    can be enriched.
  \item In general, it is not 2-surjective: not every natural
    transformation is \V-natural.  It is 2-surjective, however,
    whenever the functor $\V(I,-)$ is faithful, as for topological
    spaces and abelian groups.  But when \V\ is, say, simplicial sets,
    the functor $\V(I,-)$ is not faithful, since a simplicial map is
    not determined by its action on vertices.
  \item It is always 3-surjective: a \V-natural transformation is
    determined uniquely by its underlying ordinary natural
    transformation.
  \end{itemize}
  Thus in general, $(-)_0$ forgets at most stuff, but when $\V(I,-)$
  is faithful, it forgets at most structure.

  Now, what is the fiber over an ordinary category $C$?  Its objects
  are enrichments of $C$, its morphisms are \V-functors whose
  underlying ordinary functors are the identity, and its 2-cells are
  \V-natural transformations whose components are identities.  Such a
  2-cell is merely the assertion that two \V-functors are equal, so in
  general this is a 1-category.  This is what we expect, since $(-)_0$
  forgets at most stuff.  However, when $\V(I,-)$ is faithful, a
  \V-functor is determined by its underlying functor, so the fiber is
  in fact a poset, as we expect it to be since in this case the
  functor forgets at most structure.
\end{eg}

It would be nice to have a good example of a 2-functor which forgets
at most stuff and whose fibers are 2-posets that are not 1-categories,
but I haven't thought of one.

\begin{eg}
  In order to do an example that forgets 2-stuff, consider the
  forgetful 2-functor
  \[[\text{pairs of categories}]\too{} [\text{categories}].\]
  This functor is not $j$-surjective for any $j\le 3$, so it forgets
  at most 2-stuff.  And here the (essential) fiber over a category
  is a genuine 2-category: we can have arbitrary functors and natural
  transformations living on that extra category we forgot about.
\end{eg}

Can we make this formal and use it as an alternate characterization of
how much a functor forgets?  The answer is: `sometimes.'  Here's
what's true always in dimension one:

\begin{itemize}
\item If a functor is an equivalence, then all its essential fibers
  are contractible ($(-2)$-categories);
\item If it is full and faithful, then all its essential fibers are empty
  or contractible ($(-1)$-categories); 
\item If it is faithful, then all its essential fibers are posets;
\item and of course, if it is arbitrary, then its essential fibers can
  be arbitrary categories.
\end{itemize}

However, in general none of the implications above can be reversed.
This is because a statement about the essential fibers really tells us
only about the arrows which live over isomorphisms, while full and
faithful tell us something about \emph{all} the arrows.

There are, however, two cases in which the above implications
\emph{are} reversible:
\begin{enumerate}
\item When all categories involved are \emph{groupoids}.  This is
  because in this case, all arrows live over isomorphisms, since they
  all \emph{are} isomorphisms.
\item If the functor is a \emph{fibration} in the categorical sense.
\end{enumerate}

Being a fibration in the categorical sense is like being a fibration
in the topological sense, except that (1) we allow ourselves to lift
arrows that have direction, since our categories have such arrows, and
(2) we don't allow ourselves to take just any old lift, but require
that the lift satisfy a nice universal property.  I won't give the
formal definition here, since you can find it in many places; instead
I want to try to explain what it means.

The notion essentially means that the extra properties, structure, or
stuff that lives upstairs in $E$ can be `transported' along arrows
downstairs in $B$ in a \emph{universal} way.  When we're transporting
along arrows downstairs that are invertible, like paths in topology or
arrows in an $n$-groupoid, this condition is unnecessary since the
invertibility guarantees that we aren't making any irreversible
changes.  My favorite example is the following.

\begin{eg}
  Let $B$ be the category of rings and ring homomorphisms.  Let $E$ be
  the category whose objects are pairs $(R,M)$ where $R$ is a ring and
  $M$ is an $R$-module, and whose morphisms are pairs $(f,\ph)\maps
  (R,M)\to (S,N)$ where $f\maps R\to S$ is a ring homomorphism and
  $\ph\maps M\to N$ is an `$f$-equivariant map', i.e.\
  $\ph(rm)=f(r)\ph(m)$.  Then if $f\maps R\to S$ is a ring
  homomorphism and $N$ is an $S$-module, there is a canonical
  associated $R$-module $f^*N$---namely, $M$ with $R$ acting through
  $f$---and a canonical $f$-equivariant map $f^*N\to N$---namely the
  identity map.  This map is `universal' in a suitable sense, and is
  clearly what we should mean by `transporting' $N$ backwards along
  $f$.
\end{eg}

The formal definition of fibration simply makes this notion precise.

The introduction of directionality here also means
that we get different things by transporting objects along arrows
backwards and forwards.  In the above example, the dual construction
would be to take an $R$-module $M$ and construct an $S$-module $f_!M =
S\ten_R M$ by `extending scalars' to $S$.  Again this comes with a
canonical $f$-equivariant map $M\to f_! M$.  Thus there are actually
two notions of categorical fibration; for historical reasons, the
`backwards' one is usually called a \textbf{fibration} and the
`forwards' one an \textbf{opfibration} (or a `cofibration', but we
eschew that term because it carries the wrong topological intuition).
Either one works equally well for the characterization of
forgetfulness by fibers.

Another nice thing about the notion of categorical fibration is that
while the principle of Galois theory does not apply, in general, to
functors between arbitrary categories, it does apply to fibrations.
Recall that in the groupoid case, fibrations over a base space
($n$-groupoid) $B$ with fiber $F$ are equivalent to functors $B\to
\AUT(F)$.  One can show that for a base category $B$, fibrations over
$B$ are equivalent to (weak) functors $B\op\to\Cat$.  The way to think
of this is that since our arrows are no longer necessarily invertible,
the induced morphisms of fibers are no longer necessarily
automorphisms, nor are all the fibers necessarily the same.  Thus
instead of the automorphism $n$-group of `the' fiber, we have to use
the whole \emph{category} of possible fibers: in this case, \Cat,
since the fibers are categories.

Fibrations also have the nice property that the essential preimage is
equivalent to the literal or `strict' preimage.  Since many forgetful
functors, like those above, are fibrations, in such cases we can use
the strict preimage instead of the essential one.  In fact, a much
weaker property than being a fibration is enough for this; it suffices
that objects upstairs can be transported along `equivalences'
downstairs (which coincides with the notion of fibrations in the
$n$-groupoid case, when all morphisms are equivalences).  This is true
in many examples which are not full-fledged fibrations.

This advantage also implies, however, that there is a sense in which
the notion of categorical fibration is `not fully weak'.  Ross Street
has defined a weaker notion of fibration which does not have this
property, and which makes sense in any (weak) 2-category.  It is easy
to check that this weaker notion of fibration also suffices for the
characterization of forgetfulness via fibers.  Of course, unlike for
traditional fibrations, in this case it is essential that we use
essential preimages, rather than strict ones, since the two are no
longer equivalent.

There ought to be a notion of categorical fibration for
$n$-categories.  Some people have studied particular cases of this.
Claudio Hermida has studied $2$-fibrations between $2$-categories.
Andr\'e Joyal, Jacob Lurie, and others have studied various notions of
fibration between quasi-categories, which are one model for
$(\infty,1)$-categories; see Joyal's introduction to quasi-categories
in this volume for more details.

Let's end this section by formulating a hypothesis about the behavior
of fibers for $n$-categories.  Generalizing an idea from the first
lecture, let's say that a functor \textbf{forgets at most $k$-stuff}
if it is $j$-surjective for $j>k+1$.

\begin{hyp}
  If a functor between $n$-categories forgets at most $k$-stuff, then
  its fibers are $k$-categories (which we take to include
  poset-enriched $k$-categories).  The converse is true for
  $n$-groupoids and for $n$-categorical fibrations.
\end{hyp}

We've checked this hypothesis above for $n=1$ and for $n$-groupoids
(modulo Grothendieck).  Actually, we only checked it for
$(1,1)$-categories, while to be really consistent, we should check it
for $(1,2)$-categories too, but I'll leave that to you.  For $n=0$ it
says that an isomorphism of posets has contractible fibers (obvious)
and that an inclusion of a sub-poset has fibers which are empty or
contractible (also obvious).  Surely someone can learn about
$2$-fibrations and check this hypothesis for $n=2$ as well.

As one last note, recall that in the topological case, when we studied
Postnikov towers in \S\ref{sec:cohomology}, we were able, by the magic
of homotopy theory, to convert all the maps in our factorization into
fibrations.  It would be nice if a similiar result were true for
categorical fibrations.  It isn't true as long as we stick to plain
old categories, but there's a sense in which it becomes true once we
generalize to things called `sites' and their corresponding `topoi'.
I won't say any more about this, but it leads us into the next topic.

\subsection{$n$-Topoi}
\label{sec:n-topoi}

Knowing about the existence of $n$-posets and how they fit into the
enrichment table also clarifies the notion of topos, and in particular
of $n$-topos.

Topos theory (by which is usually meant what we would call 1-topos
theory and 0-topos theory; I'll explain later) is a vastly beautiful
and interconnected edifice of mathematics, which can be quite
intimidating for the newcomer, not least due to the lack of a unique
entry point.  In fact, the title of Peter Johnstone's epic compendium
of topos theory, {\sl Sketches of an Elephant}, compares the many
different approaches to topos theory to the old story of six blind men
and an elephant.  (The six blind men had never met an elephant before,
so when one was brought to them, they each felt part of it to
determine what it was like.  One felt the legs and said ``an elephant
is like a tree,'' one felt the ears and said ``an elephant is like a
banana leaf,'' one felt the trunk and said ``an elephant is like a
snake,'' and so on.  But of course, an elephant is all of these things
and none of them.)

So what is a topos anyway?  For now, I want you to think of a
(1-)topos as \emph{a 1-category that can be viewed as a generalized
  universe of sets}.  What this turns out to mean is the following:
\begin{itemize}
\item A topos has limits and colimits;
\item A topos is cartesian closed; and
\item A topos has a `subobject classifier'.
\end{itemize}

It turns out that this much does, in fact, suffice to allow us to more
or less replace the category of sets with any topos, and build all of
mathematics using objects from that topos instead of our usual notion
of sets.  (There's one main caveat I will bring up below.)

Now, how can we generalize this to $n$-topoi for other values of $n$?
I'm going to instead ask the more general question of how we can
generalize it to $(n,m)$-topoi for other values of $n$ and $m$.  I
claim that a sensible generalization should allow us to assert that:

\begin{quote}
  \emph{An $(n,m)$-topos is an $(n,m)$-category that can be viewed as
    a generalized universe of $(n-1,m-1)$-categories.}
\end{quote}

One thing this tells us is that we shouldn't expect to have much of a
notion of topos for $n$-groupoids: we don't want to let ourselves drop
off the enrichment table.  Inspecting the definition of 1-topos
confirms this: groupoids generally don't have limits or colimits, let
alone anything fancier like exponentials or subobject classifiers.
The only groupoid that is a topos is the trivial one.  This is a bit
unfortunate, since it means we can't test our hypotheses using
homotopy theory in any obvious way, but we'll press on anyway.

Let's consider our lower-dimensional world.  What should we mean by a
`$(0,1)$-topos'?  (I'm going to abuse terminology and call this a
`0-topos', since as we saw above, we expect the only $(0,0)$-topos to
be trivial.)  Well, our general philosophy tells us that it should be
a poset that can be viewed as a generalized universe of truth values.
At this point you may think you know what it's going to turn out to
be---and you may be right, or you may not be.

What do the characterizing properties of a 1-topos say when
interpreted for posets?  Limits in a poset are meets (greatest lower
bounds), and colimits are joins (least upper bounds), so our 0-topoi
will be complete lattices.  Being cartesian closed for a poset means
that for any elements $b,c$ there exists an object  $b\Impl c$ such
that
\[
a \le (b\Impl c) \quad\text{if and only if}\quad (a\meet b) \le c
\]
As we expected, this structure makes our poset look like a generalized
collection of truth values: we have a conjunction operation $\meet$, a
disjunction operation $\join$, and an implication operation $\Impl$.
We can define a negation operator by $\neg a = (a \Impl \bot)$, which
turns out to behave just as we expect, except that in general
$\neg\neg a \neq a$.  Thus the logic we get is not \emph{classical}
logic, but \emph{constructive} logic, in which the principle of double
negation is denied (as are equivalent statements such as the `law of
excluded middle', $a\join \neg a$).  Boolean algebras, which model
classical logic, are a special case of these cartesian closed posets,
which are called \textbf{Heyting algebras}.  Thus, a 0-topos is
essentially just a complete Heyting algebra.

Now, one of the most exciting things in topos theory is that Heyting
algebras turn up in topology!  Namely, the lattice of open sets
$\scrO(X)$ of any topological space $X$ is a complete Heyting algebra,
and any continuous map $f\maps X\to Y$ gives rise to a map of posets
$f^*\maps \scrO(Y)\to \scrO(X)$ which preserves finite meets and
arbitrary joins (but not $\Impl$).  Thus, we can view complete Heyting
algebras as a sort of `generalized topological space'.  When we do
this, we call them \textbf{locales}.  So a more correct thing to say is
that 0-topoi are the same as locales.

Now, I didn't mention the subobject classifier.  In fact, it doesn't
turn out to mean anything interesting for posets.  This makes us
wonder what its appearance for 1-topoi means.  One answer is that
\emph{it allows us to apply the principle of Galois theory inside a
  1-topos}.

What should that mean?  Well, what does the principle of Galois theory
(suitably generalized to nonidentical fibers) say for sets?  It says,
first of all, that functions $p\maps E\to B$, for sets $E$ and $B$,
are equivalent to functors $B\to\Set$.  This is straightforward: we
take each $b\in B$ to the fiber over it.

But what if we reduce the dimension of the fibers?  A function $p\maps
E\to B$ whose fibers are $(-1)$-categories, i.e.\ truth values, is
just a subset of $B$, and the principle of Galois theory says that
these should be equivalent to functors from $B$ to the category of
truth values, which is the poset $\{\text{false} \le \text{true}\}$,
often written $\mathbf{2}$.  This is just the correspondence between
subsets and their \emph{characteristic functions}.

Now, a subobject classifier is a categorical way of saying that you
have on object $\Omega$ which acts like $\mathbf{2}$: it is a target
for characteristic functions of subobjects (monomorphisms).  Thus,
this condition in the definition of 1-topos essentially tells us that
we can apply the principle of Galois theory inside the topos.  (It
turns out that the unrestricted version for arbitrary functions
$p\maps E\to B$ is also true in a topos, once you figure out how to
interpret it correctly.)

Now, since a 1-topos is a generalized universe of sets and contains an
object $\Omega$ which acts as a generalization of the poset
$\mathbf{2}$ of truth values, we naturally expect $\Omega$ to be a
generalized universe of truth values, i.e.\ a 0-topos.  This is in
fact the case, although there are couple of different ways to make
this precise.

One such way is to consider the \textbf{subterminal objects} of the
topos, which are the objects $U$ such that for any other object $E$
there is \emph{at most} one map $E\to U$.  They are called
`subterminal' because they are the subobjects of the terminal object
$1$, which are by definition the same as the maps $1\to \Omega$, or
the `points' of $\Omega$.  They can also be described as the objects
which are `representably $(-1)$-categories', since each hom-set
$C(E,U)$ has either 0 or 1 element, so it is precisely a truth value.
Thus it makes sense that the collection of subterminal objects turns
out to be a 0-topos, whose elements are the `internal truth values' in
our given 1-topos.  Since in general the logic of a 0-topos is
constructive, not classical, the internal logic of a 1-topos is also
in general constructive; this is the one caveat I mentioned earlier
for our ability to redo all of mathematics in an arbitrary topos.

Thus every 1-topos, or universe of sets, contains inside it a 0-topos,
or universe of truth values.  We can also go in the other direction:
given a locale $X$ (a 0-topos), we can construct its category $\Sh(X)$
of \textbf{sheaves}, by an obvious generalization of the notion of
sheaves on a topological space, and this turns out to be a 1-topos,
which we regard as `the category of sets in the universe parametrized
by $X$'.  As we expect, the subobject classifier in $\Sh(X)$ turns out
to be $\scrO(X)$.  In fact, this embeds the category of locales in the
(2-)category of topoi, which leads us to consider any 1-topos as a
vastly generalized kind of topological space.

As a side note, recall that in \S\ref{sec:enrichment-posets} we
observed that a set is a groupoid enriched over truth values.  Thus
you might expect that the objects of $\Sh(X)$, which intuitively are
`sets in the universe where the truth values are $\scrO(X)$', could be
defined as `groupoids enriched over $\scrO(X)$'.  This is almost
right; the problem is that all the objects of such a groupoid turn out
to have `global extent', while an arbitrary sheaf can have objects
which are only `partially defined'.  We can, however, make it work if
we consider instead groupoids enriched over a suitable
\emph{bi}category constructed from $\scrO(X)$.

Anyway, these relationships between 0-topoi and 1-topoi lead us to
hope that in higher dimensions, each $(n,m)$-topos will contain within
it topoi of lower dimensions, and in turn will embed in topoi of
higher dimensions via a suitable categorification of sheaves (usually
called `$n$-stacks' or simply `stacks').  Notions of $(n,m)$-topos
have already been studied for a few other values of $n$ and $m$.  For
instance, there has also been a good deal of interest lately in
something that people call `$\infty$-topoi', although from our point
of view a better name would be $(\infty,1)$-topoi.  These are special
$(\infty,1)$-categories that can be considered as a generalized
universe of homotopy types (i.e.\ $\infty$-groupoids).  And for a long
time algebraic geometers have been studying `stacks of groupoids',
which are pretty close to what we would call a `$(2,1)$-topos'.

Where does the interest in higher topoi come from?  In topology, the
principle of Galois theory already works very nicely, and people were
working with fibrations, homotopy groups, Postnikov towers, and
cohomology long before Grothendieck came along to tell them they were
really working with $\infty$-groupoids.  A fancy way to say this is
that the category of spaces is already an $(\infty,1)$-topos.

But in algebraic geometry, the Galois theory fails, because the
category under consideration is `too rigid'.  The $n$-groups $\AUT(F)$
just don't exist.  So what the algebraic geometers do is to take their
category and \emph{embed} it in a larger category in which the desired
objects \emph{do} exist; we would say that they embed it in a
$(2,1)$-topos (if they're only interested in one level of
automorphisms) or an $(\infty,1)$-topos (if they're interested in the
full glory of homotopy theory).  The way they do this is with a
suitable generalization of the sheaf construction to arbitrary
categories.

In the case $m>1$, it is not clear whether there exists a single
notion of $(n,m)$-topos that shares most of the good properties of
$1$-topoi.  Several people have studied this question, however, with
some partial encouraging results; the `fibrational cosmoi' of Ross
Street can be viewed as generalized universes of 1-categories, and
more recently Mark Weber has studied certain special cosmoi under the
name `2-topos'.  One of the defining properties of a cosmos is the
existence of `presheaf objects' which allow the application of the
principle of Galois theory to internal fibrations in the 2-category
(suitably defined).  Some people speak of this as ``considering sets
to be generalized truth values''.

\subsection{Geometric morphisms, classifying topoi, and $n$-stuff}
\label{sec:geom-morph}

In this section we'll see that morphisms between topoi admit similar
`Postnikov' factorizations, which in turn tell us interesting things
about the logical theories they `classify'.  This section will
probably be most interesting to readers with some prior acquaintance
with topos theory, but I've tried to make it as accessible as
possible.

Recall that a continuous map $f\maps X\to Y$ of topological spaces
gives rise to a function $f^*\maps \scrO(Y)\to \scrO(X)$ which
preserves finite meets and arbitrary joins.  Let $X$ and $Y$ be
locales and $\scrO(X)$ and $\scrO(Y)$ the corresponding complete
Heyting algebras; we \emph{define} a map of locales $f\maps X\to Y$ to
be a function $f^*\maps \scrO(Y)\to \scrO(X)$ preserving finite meets
and arbitrary joins.  We distinguish notationally between the locale
$X$ and its poset of `open sets' $\scrO(X)$ because the maps go in the
opposite direction, even though the locale $X$ technically consists of
nothing but $\scrO(X)$.

Similarly, let $X$ and $Y$ be topoi, and $\scrS(X)$ and $\scrS(Y)$
their corresponding 1-categories.  We define a map of 1-topoi $f\maps
X\to Y$ to be a functor $f^*\maps \scrS(Y)\to \scrS(X)$ which
preserves finite limits and arbitrary colimits; these maps are called
\textbf{geometric morphisms} for historical reasons.

Now, it turns out that for any (small) category $C$, the category
$\Set^C$ of functors from $C$ to \Set\ is a topos, and functors
$f\maps C\to D$ give rise to geometric morphisms $\fhat\maps \Set^C\to
\Set^D$.  We can thus ask how properties of the functor $f$ are
reflected in properties of the geometric morphism $\fhat$.  It turns
out that we have the following dictionary (at least, `modulo splitting
idempotents', which is something I don't want to get into---just
remember that this is all morally true, but there are some details.)

\begin{center}
  \begin{tabular}{ccc}
    $f$ is full and faithful & $\sim$ & $\fhat$ is an `inclusion'\\
    $f$ is essentially surjective & $\sim$ & $\fhat$ is a `surjection'\\
    $f$ is faithful & $\sim$ & $\fhat$ is `localic'\\
    $f$ is full and essentially surjective & $\sim$ & $\fhat$ is `hyperconnected'
  \end{tabular}
\end{center}

What do all those strange terms on the right mean?  I'm certainly not
going to define them!  But I'll try to give you some idea of how to
think about them.  The notions of `inclusion' and `surjection' are
suitable generalizations of the correspondingly named notions for
topological spaces.  Moreover, just as is true for spaces, any
geometric morphism factors uniquely as a surjection followed by an
inclusion; this also parallels one of our familiar factorizations for
functors.  This part of the correspondence should make some intuitive
sense.

To explain the term `localic', consider a geometric morphism $p\maps
E\to S$.  It turns out that we can think of this either as a map
between two topoi in the universe of sets, or we can use it to think
of $E$ as an \emph{internal} topos in the generalized universe
supplied by the topos $S$.  We say that the morphism $p$ is `localic'
if this internal topos is equivalent to the sheaves on some internal
locale in $S$.  It turns out that there is another sort of morphism
called `hyperconnected' such that every geometric morphism factors
uniquely as a hyperconnected one followed by a localic one, and this
too corresponds to a factorization we know and love for functors.

Moreover, every inclusion is localic, and every hyperconnected
morphism is a surjection, and it follows that every geometric morphism
factors as a hyperconnected morphism, followed by a surjective localic
map, followed by an inclusion.  This should also look familiar in the
world of functors.

Now I want to explain why these classes of geometric morphisms in fact
have an \emph{intrinsic} connection to the notions of properties,
structure, and stuff, but to do that I have to talk about
`classifying topoi'.

The basic idea of classifying topoi is that we can apply the principle
of Galois theory once again, only this time we apply it in the
2-category \emph{of topoi}, and we apply it to classify \emph{models
of logical theories}.  Let $\bbT$ be a \textbf{typed logical
theory}; thus it has some collection of `types', some `function and
relation symbols' connecting these types, and some `axioms' imposed on
the behavior of these symbols.  An example is the theory of
categories, which has two types $O$ (`objects') and $A$ (`arrows'),
three function symbols $s,t\maps A\to O$, $i:O\to A$, a relation
symbol $c$ of type $A\times A \times A$ (here $c(f,g,h)$ is intended
to express the assertion that $h=g\circ f$), and various axioms, such
as
\[(t(f) = s(g)) \Rightarrow \exists! h \, c(f,g,h)\] 
(which says that any
two composable arrows have a unique composite).  A model of such a
theory assigns a set to each type and a function or relation to each
symbol, such that the axioms are satisfied; thus a model of the theory
of categories is just a small category.

The fact that a topos is a generalized universe of sets implies that
we can consider models of such a theory in \emph{any} topos, not just
the usual topos of sets.  It turns out that for suitably nice theories
$\bbT$ (called `geometric' theories), there exists a topos $[\bbT]$
such that for any other topos $B$, the category of models of $\bbT$ in
$B$ is equivalent to $\hom(B,[\bbT])$, the category of geometric
morphisms from $B$ to $[\bbT]$ (remember that 1-topoi form a 2-category).
Thus, once again, some structure `in' or `over' $B$ can be classified
by functors from $B$ to a `classifying object'.

Now suppose that we have two theories $\bbT$ and $\bbT'$ such that
$\bbT'$ is $\bbT$ with some extra types, symbols, and/or axioms added.
Since this means that any model of $\bbT'$ gives, by neglect of
structure, a model of $\bbT$, by the Yoneda lemma we have a geometric
morphism $p\maps [\bbT']\to{} [\bbT]$.  It turns out that

\begin{center}
  \begin{tabular}{ccp{3in}}
    $p$ is an inclusion & when & $\bbT'$ adds only extra axioms to $\bbT$\\
    $p$ is localic & when & $\bbT'$ adds extra functions, relations, and
    axioms to $\bbT$, but no new types\\
    $p$ is a surjection & when & $\bbT'$ adds extra types, symbols, and
    axioms to \bbT, but no new properties of the existing types and
    symbols in \bbT\ are implied by this new structure.\\
    $p$ is hyperconnected & when & $\bbT'$ adds extra types to \bbT,
    along with symbols and axioms relating to these new types, but no
    new functions, relations, or axioms on the existing types in \bbT\
    are implied by this new structure.
  \end{tabular}
\end{center}

(There are various ways to make these notions precise, which I'm not
going to get into.)  Thus, these classes of geometric morphisms
actually directly encode the notions of forgetting properties
(axioms), structure (function and relation symbols), and/or stuff
(types).

Notice that localic morphisms are those that add no new types; this is
consistent with the fact that locales are 0-topoi, and 0-categories
know only about properties ($(-1)$-stuff) and structure (0-stuff), not
stuff (1-stuff).  In particular, a classifying topos $[\bbT]$ is
equivalent to a topos of sheaves on a locale precisely when the theory
\bbT\ has no types.  Such a theory, which consists only of
propositions and axioms, is called a \textbf{propositional theory}; from
our point of view, we might also call it a `0-theory', with the more
general typed theories considered above being `1-theories'.  As far as
I know, there has been very little work on notions of $n$-theories for
higher values of $n$.

Now, given the correspondence between theories and classifying topoi,
any factorization for geometric morphisms leads to a factorization for
geometric theories.  These factorizations are mostly what we would
expect, but can be slightly different due to the requirement that all
theories in sight be geometric.

\begin{eg}
  Consider the forgetful map from monoids to semigroups.  (A semigroup
  is a set with an associative binary operation.)  Considered as a
  functor $\mathbf{Mon}\to\mathbf{SGp}$, it is faithful, but not
  essentially surjective (since not every semigroup has an identity)
  or full (since not every semigroup homomorphism between monoids
  preserves the identity).  If we factor it into a
  full-and-essentially-surjective functor followed by a
  full-and-faithful one, the intermediate category we obtain is the
  category of `semigroups with identity', i.e.\ the category whose
  objects are monoids but whose morphisms do not necessarily preserve
  the identity.

  Now, the theories \bbM\ of monoids and \bbS\ of semigroups are both
  geometric, so they have classifying topoi $[\bbM]$ and $[\bbS]$, and
  as we expect there is a geometric morphism $[\bbM]\to [\bbS]$ which
  is localic.  If we factor it into a surjection followed by an
  inclusion, however, the intermediate topos we obtain is not the
  classifying topos for semigroups with identity, because that theory
  is not geometric.  Instead, the intermediate topos we get is the
  classifying topos for semigroups such that for any finite set of
  elements, there is an element which behaves as an identity for them.
  In general, however, the `identities' for different finite sets
  could be different.

  This theory is, in a sense, the `closest geometric approximation' to
  the theory of semigroups with identity.  This notion is in accord
  with the general principle (which we have not mentioned) that
  geometric logic is the `logic of finite observation'.  In this case,
  it is evident that if we can only `observe' finitely many elements
  of the semigroup, we can't tell the difference between such a model
  of our weird intermediate geometric theory and a semigroup that has
  an actual identity.
\end{eg}

These considerations may lead us to speculate that morphisms of higher
topoi, once defined, will have similar `Postnikov factorizations'.
However, in the absence of confidence that good notions of $(n,m)$-topos
exist for $m>1$, this must remain a speculation.

\subsection{Monomorphisms and epimorphisms}
\label{sec:monom-epim}

A question was asked at one point (in \S\ref{sec:internalizing}) about
whether notions like essential surjectivity can be defined purely
2-categorically, and thereby interpreted in any 2-category, the way
that epimorphisms and monomorphisms make sense in any 1-category.
This section is an attempt to partially answer that question.

The definitions of monomorphism and epimorphisms in 1-categories are
`representable' in the following sense:
\begin{itemize}
\item $m\maps A\to B$ is a monomorphism if for all $X$, the
  function
  \[C(X,m)\maps C(X,A)\to C(X,B)\]
  is injective.
\item $e\maps E\to B$ is an epimorphism if for all $X$,
  the function 
  \[C(e,X)\maps C(B,X)\to C(A,X)\]
  is injective.
\end{itemize}

Note that both notions invoke \emph{injectivity} of functions of sets.
Thus, the natural notions to consider first are functors which are
`representably' faithful or full-and-faithful.  It is easy to check
that this works in the covariant direction:
\begin{itemize}
\item A functor $p\maps A\to B$ is faithful if and only if it is
  representably faithful, i.e.\ all functors
  \[\Cat(X,p)\maps \Cat(X,A)\to \Cat(X,B)\]
  are faithful; and
\item A functor is full and faithful if and only if it representably
  full and faithful.
\item A functor is an equivalence if and only if it is representably
  an equivalence.
\end{itemize}
Thus, it makes sense to define a 1-morphism in a 2-category to be
\textbf{faithful} or \textbf{full and faithful} when it is representably
so.

We may generalize this (hypothetically) by saying that a functor
between $n$-categories is \textbf{$j$-monic} if it is $k$-surjective for
all $k > j$ (note that this is equivalent to saying that it `forgets
at most $(j-1)$-stuff'), and that a 1-morphism $p\maps A\to B$ in an
$(n+1)$-category $C$ is \textbf{$j$-monic} if all functors $C(X,p)$ are
$j$-monic.  By analogy with the above observation, we expect that
these definitions will be equivalent for the $(n+1)$-category of
$n$-categories.

Thus, every functor is 2-monic, the 1-monic functors are the faithful
ones, the 0-monic functors are the full and faithful ones, and the
$(-1)$-monic functors are the equivalences.  More degenerately, in a
1-category, every map is 1-monic, the 0-monic morphisms are the usual
monomorphisms, and the $(-1)$-monic morphisms are the isomorphisms.

We may define, dually, a 1-morphism $p\maps E\to B$ in an
$(n+1)$-category $C$ to be \textbf{$j$-epic} if all the functors
\[C(p,X)\maps C(B,X) \to C(E,X)
\]
are $(n-1-j)$-monic.  For example, in a 1-category, every morphism is
$(-1)$-epic, the 0-epic morphisms are the usual epimorphisms, and the
1-epic morphisms are the isomorphisms.

The `inversion' of numbering here may look a little strange if we
remember that every $n$-category is secretly an $\infty$-category;
when did it suddenly start to matter which $n$ we are using?  But it
turns out that the transformation above is actually exactly what is
required to make the notion independent of $n$.  For example, if
$p\maps E\to B$ is a surjective function in \Set, then $\Set(p,X)$ is
injective, hence 0-monic, for any set $X$; thus $p$ is 0-epic in the
1-category \Set.  But now consider $p$ as a functor between discrete
categories.  When $X$ is a nondiscrete category, $\Cat(p,X)$ is
faithful, but not full; hence it is only 1-monic, but by our
definition this is just what is required so that $p$ is again 0-epic.

It is easy to check that in a 2-category, every morphism is
$(-1)$-epic, and the 2-epic morphisms are the equivalences.  However,
even in the 2-category \Cat, the 0-epic and 1-epic morphisms are not
that well-behaved.  Here is what is true (proofs are left to the
reader):
\begin{itemize}
\item If a functor $p\maps E\to B$ is essentially surjective, then it
  is 0-epic.
\item Similarly, if it is full and essentially surjective, then it is
  1-epic.
\end{itemize}

However, neither implication is reversible.  For example, the
inclusion of the category $2$, which has two objects and one
nonidentity morphism between them, into the category $\mathscr{I}$,
which has two uniquely isomorphic objects, is 1-epic, but not full.
And if $p\maps E\to B$ has the property that every object of $B$ is a
retract of an object in the image of $p$, then $p$ is 0-epic, but it
need not be essentially surjective.

We thus seek for other characterizations of surjective functions in
\Set\ which will generalize better to \Cat.  It turns out that the
best-behaved notion is the following:

\begin{defn}
  An epimorphism $p\maps E\to B$ in a 1-category is a \textbf{strong
    epimorphism} if it is `left orthogonal' to monomorphisms, i.e.\
  for any monomorphism $m\maps X\to Y$, every commutative square
  \[\xymatrix{E\ar[d]_e \ar[r] & X \ar@{ >->}[d]^m\\ B\ar[r] \ar@{-->}[ur] & Y}\]
  has a unique diagonal filler.  
\end{defn}

In a category with equalizers, the orthogonality property implies that
$p$ is already an epimorphism.  In \Set, every epimorphism is strong,
but in general this is not true.

Notice that saying $p\maps E\to B$ is left orthogonal to $m\maps X\to
Y$ in the category $C$ is equivalent to saying that the following
square is a pullback:
\[\xymatrix{C(B,Y) \ar[r]\ar[d] & C(E,Y) \ar[d]\\
  C(B,X) \ar[r] & C(E,X)}\]
Therefore, we generalize this to 2-categories as follows.

\begin{defn}
  A 1-morphism $p\maps E\to B$ in a 2-category $C$ is \textbf{left
    orthogonal} to another $m\maps X\to Y$ if the square
  \[\xymatrix{C(B,Y) \ar[r]\ar[d] & C(E,Y) \ar[d]\\
    C(B,X) \ar[r] & C(E,X)}\]
  is a pullback (in a suitable 2-categorical sense).
\end{defn}

We can now check that
\begin{itemize}
\item A functor $p\maps E\to B$ is essentially surjective if and only
  if it is left orthogonal to all full and faithful functors, and
\item It is full and essentially surjective if and only if it is left
  orthogonal to all faithful functors.
\end{itemize}

The forward directions are exercises in category theory.  The idea is
that we must progressively `lift' objects, morphisms, and equations
(to show functoriality and naturality) from `downstairs' to
`upstairs'.  In both cases, for each $j$, one of the two functors is
$j$-surjective, so we can use that functor to lift the $j$-morphisms.
The reverse directions are easy using the Postnikov factorization.

We are thus motivated to define, hypothetically, a $j$-epimorphism in
an $(n+1)$-category to be a \textbf{strong $j$-epimorphism} if it is
left orthogonal (in a suitably weak sense) to all $j$-monic morphisms.
We have just shown that in \Cat, the strong 1-epics are precisely the
full and essentially surjective functors, while the strong 0-epics are
the essentially surjective functors.  Clearly all functors are strong
$(-1)$-epic, while only equivalences are strong 2-epic.

We can also prove that in 2-categories with finite limits, any
morphism which is left orthogonal to $j$-monomorphims is automatically
a $j$-epimorphism; we use 2-categorical limits such as `inserters' and
`equifiers' to take the place of equalizers in the 1-dimensional
version.  As we have seen, even in \Cat, not every $j$-epimorphism is
strong.

This leads us to formulate the following hypothesis.

\begin{hyp} \label{epimono}
  The strong $j$-epics in $n\Cat$ (that is, functors which are left
  orthogonal to all $j$-monic functors) are precisely the functors
  which are $k$-surjective for $k\le j$.  Not every $j$-epic is
  strong, even in $n\Cat$.
\end{hyp}

Since $j$-monic functors are those that forget `at most
$(j-1)$-stuff', we might say that the strong $j$-epics are the
functors which `forget no less than $j$-stuff'.  For example, the
strong 0-epics in \Cat\ are the essentially surjective functors, which
do not forget properties ($(-1)$-stuff), although they may forget
structure (0-stuff) and 1-stuff.  Similarly, the strong 1-epics, being
essentially surjective and full, do not forget properties or
structure, although they may forget 1-stuff.

What does this look like for $n$-groupoids?  For a functor between
$n$-groupoids, being $k$-surjective is equivalent to inducing a
surjection on $\pi_k$ and an \emph{in}jection on $\pi_{k-1}$.  Why?
Well, remember that $\pi_k$ of an $n$-groupoid consists of the
automorphisms of the identity $(k-1)$-morphism, modulo the
$(k+1)$-morphisms.  Thus being surjective on $k$-morphisms implies
being surjective on $\pi_k$ (although there might be new
$(k+1)$-morphisms appearing preventing it from being an isomorphism),
but also being injective on $\pi_{k-1}$, since everything we quotient
by downstairs has to already be quotiented by upstairs (although here
there might be entirely new $(k-1)$-morphisms appearing downstairs).
This is also equivalent to saying that $\pi_k$ of the homotopy fiber
is trivial.

Thus, a functor between $n$-groupoids is $j$-monic if it induces
isomorphisms on $\pi_k$ for $k > j$ and an injection on $\pi_{j}$.
Our above conjecture then translates to say that the strong $j$-epics
should be the maps $A\to W$ inducing isomorphisms on $\pi_k$ for $k<j$
and a surjection on $\pi_{j}$.  This is precisely what topologists
call a \textbf{$j$-equivalence} or a \textbf{$j$-connected map},
since it corresponds to the vanishing of the `relative homotopy
groups' $\pi_k(W,A)$ for $k\le j$.

Using this identification, we can then prove our conjecture for
$n$-groupoids (modulo Grothendieck).  Suppose we have a square
\[\xymatrix{A \ar[r]\ar[d]_p & E \ar[d]^m\\
  W\ar[r] & B}\]
of maps between $n$-groupoids (i.e.\ topological spaces), in which $p$
is a $j$-equivalence and $m$ is $j$-monic.

By magic of homotopy theory, similar to the way we can transform any
map into a fibration, we can transform the map $p$ into a `relative
cell complex'.  This means that $W$ is obtained from $A$ by attaching
`cells' $D^k$ along their boundaries $S^{k-1}$.  Since $p$ is a
$j$-equivalence, we can assume that we are only attaching cells of
dimension $k>j$.  Thus our problem is reduced to defining a lift on
each individual $D^k$.  But our assumption on $m$ guarantees that
since (inductively) we have a lift of the boundary $S^{k-1}$, the
whole cell $D^k$ must also lift, up to homotopy.

This only shows that single maps lift, but we can also show that the
appropriate square is a homotopy pullback by considering various
modified squares.  (The fanciest way to make this precise is to
construct a `Quillen model structure'.)  Thus $j$-connected maps are
in fact left `homotopy' orthogonal to $j$-monic maps.

Just as before, we can prove the converse using our knowledge of
factorizations.  Suppose we have a map $p\maps A\to W$ which is left
orthogonal to all $j$-monic maps.  By picking out a particular part of
the Postnikov factorization of $p$, we get a factorization $A\too[f]
E\too[g] W$ in which $g$ is $j$-monic and $f$ is $j$-connected.  Then
in the square
\[\xymatrix{A \ar[r]^f\ar[d]_p & E \ar[d]^g\\
  W \ar@{-->}[ur]^h \ar@{=}[r] & W}\]
there exists a diagonal lift $h$.  Consider any $k\le j$; then
$\pi_k(f)$ is either an isomorphism (if $k< j$) or surjective (if
$k=j$), by assumption.  This implies that $\pi_k(h)$ must also be
surjective.  But $\pi_k(g)\pi_k(h)$ is the identity, so in fact
$\pi_k(h)$ must be an isomorphism.  Therefore, since $hp\eqv f$,
$\pi_k(p)$ must be an isomorphism if $k<j$ and surjective if $k=j$,
since $\pi_k(f)$ is so.  Thus $p$ was already $j$-connected.

So, modulo Grothendieck's dream, we have proved
Hypothesis \ref{epimono} in the case of $n$-groupoids.
Joyal's paper on quasi-categories, in this volume, includes a theory
of factorization systems in $(\infty,1)$-categories, generalizing the
above arguments for $\infty$-groupoids to objects of any
$(\infty,1)$-category.

\subsection{Pointedness versus connectedness}
\label{sec:pointed-vs-connected}

This final section will be even more philosophical, and perhaps
controversial, than the others.  The central point I wish to make is
that the operations of looping and delooping should only be
applied to \emph{pointed} $n$-categories, just as they are only
applied in homotopy theory to pointed topological spaces.  When we do
this, various problems with the periodic table resolve themselves.

What sort of problems?  It's well-known that the hypothesis ``a
$k$-monoidal $n$-category is a $k$-degenerate $(n+k)$-category'' is
false, even in low dimensions, if you interpret `is' as referring to a
fully categorical sort of equivalence.  The simplest example is that
while a monoid `is' a one-object category in a certain sense, the
category of monoids is not equivalent to the (2-)category of
categories-with-one-object.  Similarly, the 2-category of monoidal
categories is not equivalent to the (3-)category of one-object
bicategories, and so on.  In general, the objects and morphisms turn
out mostly correct, but the higher-level transformations and so on are
wrong.  Eugenia Cheng and Nick Gurski have investigated in detail what
happens and how you can often carefully chop things off at a
particular level in the middle to get an equivalence, but here I want
to consider a different point of view.

Let's consider the case of groupoids.  The topological version of a
one-object groupoid is a $K(G,1)$, so the periodic table leads us to
expect that the homotopy theory of $K(G,1)$s should be equivalent to
the category of groups.  This is true, but only if we interpret the
$K(G,1)$s as \emph{pointed} spaces and the corresponding homotopy
theory likewise.  Otherwise, we get the theory of groups and group
homomorphisms modulo conjugation.

A related issue is that the homotopy groups $\pi_n$, and in particular
$\pi_1$, are really only defined on \emph{pointed} spaces.  While it's
true that different choices of basepoint give rise to isomorphic
groups (at least for a connected space), the isomorphism is not
canonical.  In particular, this means that $\pi_n$ is not
\emph{functorial} on the category of unpointed spaces.

Thus, by analogy with topology, we are motivated to consider `pointed
categories'.  A \textbf{pointed $n$-category} is an $n$-category $A$
equipped with a functor $1\to A$ from the terminal $n$-category (which
has exactly one $j$-morphism for every $j$).  Note that this is
essentially the same as choosing an object in $A$.  A \textbf{pointed
  functor} between two pointed $n$-categories is a functor $A\to B$
such that
\[\xymatrix{1 \ar[rr] \ar[dr] && B \lltwocell\omit{<-3>\eqv} \\ & A \ar[ur]}\]
commutes up to a specified natural equivalence.  A \textbf{pointed
  transformation} is a transformation $\alpha$ such that
\[\xymatrix{1 \ar[rr] \ar[dr] && B \\ & A \urtwocell}\]
commutes with the specified equivalences up to an invertible
modification.  And so on for higher data.

What does this look like in low dimensions?  A pointed set is just a
set with a chosen element, and a pointed function between such sets is
a function preserving the basepoints.  More interestingly, a pointed
category is a category $A$ with a chosen base object $*\in A$, a
pointed functor is a functor $f\maps A\to B$ equipped with an
isomorphism $f*\iso *$, and a pointed natural transformation is a
natural transformation $\alpha\maps f\to g$ such that
\[\xymatrix{f\ast \ar[rr]^\alpha \ar[dr]_\iso && g\ast\\ & \ast\ar[ur]_\iso}\]
commutes.

We have only required our basepoints to be preserved up to coherent
equivalence, in line with general $n$-categorical philosophy, but we
now observe that we can always `strictify' a pointed functor to
preserve the basepoints on the nose.  Define $f'$ to be $f$ on all
objects except $f'*=*$, with the action on arrows defined by
conjugating with the given isomorphism $f*\iso *$.  Then $f'*=*$, and
$f'$ is isomorphic to $f$ via a pointed natural isomorphism.  Thus the
2-category of pointed categories, pointed functors, and pointed
transformations is biequivalent to the 2-category of pointed
categories, \emph{strictly} pointed functors, and pointed
transformations.  We expect this to be true in higher dimensions as
well.  Observe that if $f$ and $g$ are strictly pointed functors, then
a pointed natural transformation $\alpha\maps f\to g$ is just a
natural transformation $f\to g$ such that the component $\alpha_* =
1_*$.

Now we can define a functor $\Omega$ from the 2-category of pointed
categories to the category of monoids.  We take our pointed functors
to be strict for convenience, since as we just saw there is no loss in
doing so; otherwise we would just have to conjugate by the
isomorphisms $f*\iso *$.  We define $\Omega A = A(*,*)$ on objects,
and for a strictly pointed functor $f\maps A\to B$, we get a monoid
homomorphism $A(*,*)\to B(f*,f*) = B(*,*)$.  Finally, since as we
observed above a pointed transformation between strictly pointed
functors is the identity on $*$, these transformations induce simply
identities, which is good since those are the only 2-cells we've got
in our codomain!

In the other direction, we construct a functor $B$ from the category
of monoids to the 2-category of pointed categories, sending a monoid
$M$ to the category $BM$ with one object $*$ and $BM(*,*)=M$, and a
monoid homomorphism to the obvious (strictly) pointed functor.  We now
observe that $B$ is left adjoint to $\Omega$ (in a suitable sense),
and that moreover the adjunction restricts to an adjoint biequivalence
between the category of monoids (regarded as a locally discrete
2-category) and the 2-category of pointed categories with exactly one
isomorphism class of objects (which we may call `pointed and
connected').

Similarly, one can construct `adjoint' functors $B$ and $\Omega$
between monoidal categories and pointed bicategories, and show that
they restrict to inverse `triequivalences' between monoidal categories
and pointed bicategories with one equivalence class of objects
(`pointed connected bicategories').  We can also do this for
commutative monoids and `pointed monoidal categories', but in fact
here the word `pointed' becomes redundant: every monoidal category has
an essentially unique basepoint, namely the unit object.  (Similarly,
any monoid has a unique basepoint, namely its identity.  This is
because the terminal monoid and the terminal monoidal category are
also `initial' in a suitable sense.)  We then obtain an adjoint
biequivalence between commutative monoids and connected monoidal
categories.

Composing these two adjunctions, we obtain an adjoint pair $B^2\adj
\Omega^2$ between commutative monoids and pointed bicategories.  This
restricts to a biequivalence between commutative monoids and pointed
bicategories with one equivalence class of objects and one isomorphism
class of 1-morphisms (`pointed 1-connected bicategories').

All of these equivalences carry over in an obvious way to the groupoid
cases, so that groups are equivalent to pointed groupoids with one
isomorphism class of objects, groupal groupoids (2-groups) are
equivalent to pointed 2-groupoids with one equivalence class of
objects, and abelian groups are equivalent to 2-groups with one
isomorphism class of objects, and also to 2-groupoids with one
equivalence class of objects and one isomorphism class of morphisms.
These are well-known topological results.

This suggests the following pointed version of the correspondence
described in the periodic table.  Say that an $n$-category is
\textbf{$i$-connected} if it has exactly one equivalence class of
$j$-morphisms for $0\le j \le i$.

\begin{hyp} [Delooping Hypothesis]\label{hyp:delooping}
  There is an adjoint pair $B^i\adj \Omega^i$ between $k$-monoidal
  $n$-categories and (pointed) $(k-i)$-monoidal $(n+i)$-categories,
  which restricts to an equivalence between $k$-monoidal
  $n$-categories and (pointed) $(i-1)$-connected $(k-i)$-monoidal
  $(n+i)$-categories.
\end{hyp}
\noindent
We have placed ``pointed'' in parentheses because it is 
expected to be redundant for $k>i$.

We have called these functors $\Omega$ and $B$ by analogy with the
corresponding topological constructions of loop space and delooping
(or `classifying space').  Note that topologists usually say that the
left adjoint of $\Omega$ is the `suspension' functor $\Sigma$, rather
than the `delooping' functor $B$.  This is because they often consider
the functor $\Omega$ to take its values just in spaces, rather than
monoidal spaces (say, $A_\infty$-spaces).  We would get a
corresponding adjoint pair in our situation by composing the two
adjunctions
\[\xymatrix{\text{$n$-categories} \rrtwocell^F_U{'\bot} &&
  \parbox{2cm}{monoidal\\$n$-categories} \rrtwocell^B_\Omega{'\bot} &&
  \parbox{3cm}{pointed \\ $(n+1)$-categories}}\]
where $F\adj U$ is the free-forgetful adjunction.  Then 
$\Sigma A = BF(A)$, the delooping of the free monoidal $n$-category on
$A$, is what deserves to be called the `suspension' of $A$.

Note that there is also a forgetful functor from pointed
$(n+1)$-categories to unpointed $(n+1)$-categories, which has a
\emph{left} adjoint $(-)_+$ called `adding a disjoint basepoint'.

Let us investigate further the question of `connectedness'.  Recall
from \S\ref{sec:homotopy-n-types} that for a space $X$ we say that
\textbf{$\pi_j(X)$ vanishes for all basepoints} if given any $f\maps
S^j\to X$, there exists $g\maps D^{j+1}\to X$ extending $f$.  When $X$
is nonempty, this is equivalent to requiring that the actual groups
$\pi_j(X)$ vanish for all base points.  Topologists define a nonempty
space $X$ to be \textbf{$k$-connected} if $\pi_i(X)$ is trivial for $j
\le k$ and all basepoints.  (We'll deal with the empty set later.)

We can generalize this to $n$-categories in a straightforward way, but
we use a different terminology because unlike $n$-groupoids,
$n$-categories are not characterized by a list of homotopy groups.  We
say that an $n$-category \textbf{has no $j$-homotopy} when any two
parallel $j$-morphisms are equivalent.  Another way to say this, which
is closer to the topology, is to define $S^j$ to be the $n$-category
consisting of two parallel $j$-morphisms, and $D^{j+1}$ to consist of
a $(j+1)$-equivalence between two parallel $j$-morphisms; then $X$ has
no $j$-homotopy just when all maps $S^j\to X$ extend to $D^{j+1}$.  We
can then define an $n$-category to be \textbf{$k$-connected} if
it has no $j$-homotopy for $j\le k$.

Note that since we don't know what a $(-1)$-morphism is, the category
$S^{-1}$ can only be empty.  And since in general $D^j$ is generated
by a single $j$-equivalence, $D^0$ should just consist of a single
object.  Thus an $n$-category has no $(-1)$-homotopy just when it is
nonempty.  Similarly, $S^{-2}$ and $D^{-1}$ should both be empty, so
every $n$-category has no $(-2)$-homotopy.  Therefore, just as for
groupoids, an $n$-category is always $(-2)$-connected, is
$(-1)$-connected when it is nonempty, and for $k\ge 0$ it is
$k$-connected when it has precisely one isomorphism class of $j$-cells
for $0\le j\le k$.  Thus, this definition of connectedness agrees with
the one we gave just before Hypothesis~\ref{hyp:delooping}.  Moreover,
the new definition allows that hypothesis to make sense even for
$i=0$, in which case it says that all $k$-monoidal $n$-categories are
nonempty and come equipped with an essentially unique basepoint (the
unit object).

Now, what about that pesky empty set?  By classical topological
definitions, the empty set is unquestionably both connected (it is not
the disjoint union of two nonempty open sets) and path-connected (any
two points in it are connected by a path).  But by our definitions,
although it has no $0$-homotopy, it is not $0$-connected, because it
does not have no $(-1)$-homotopy.  (Of course, the empty set is the
only space with no $0$-homotopy which is not $0$-connected.)

This disagreement is perhaps a slight wart on our definitions.
However, it is worth pointing out that $\pi_0$ of the empty set must
also be empty.  In particular, it is not equal to $0=\{0\}$.  So the
only way the empty set can be $0$-connected, if we use the topological
definition that $X$ is $k$-connected $\pi_j(X)=*$ for all $0\le j\le
k$, is if we maintain that $\pi_0(X)$, just like the other $\pi_j$,
requires a base point to be defined (in which case it is a pointed
set).  In this case, since the empty set has no basepoints, it is
still vacuously true that $\pi_0(\emptyset)=0$ for all basepoints.

We would like to emphasize the crucial distinction between
\emph{connected} (having precisely one equivalence class of objects)
and being \emph{pointed} (being \emph{equipped} with a chosen object).
Clearly, every connected $n$-category can be pointed in a way which is
unique up to equivalence, \emph{but not up to unique equivalence}.
Similarly, functors between connected $n$-categories can be made
pointed, but not in a unique way, while transformations and higher
data can \emph{not} in general be made pointed at all.  Thus the
$(n+1)$-categories of connected $n$-categories and of pointed
connected $n$-categories are not equivalent; the latter is equivalent
to the $n$-category of monoidal $(n-1)$-categories, but the former is
not.

This distinction explains an observation due to David Corfield that
the periodic table seems to be missing a row.  If in the periodic
table we replace `$k$-monoidal $n$-categories' by `$(k-1)$-connected
$(n+k)$-categories', then the first row is seen to be the
$(-2)$-connected things (that is, no connectivity imposed) while the
second row is the $0$-connected things.  Thus there appears to be a
row missing, consisting of the $(-1)$-connected, or nonempty, things,
and morover the top row should be shifted over one to keep the
diagonals moving correctly.  So we should be looking at a table like
this:

\vbox{
\vskip 1em
\begin{center}
{\bf THE CONNECTIVITY PERIODIC TABLE }
\vskip 1em
\begin{tabular}{|l|c|c|c|}  \hline
    & $\mathbf{\mathit n = -1}$ & $\mathbf{\mathit n = 0}$ &
    $\mathbf{\mathit n = 1}$\\ \hline
    $\mathbf{\mathit k = -1}$
    & truth values & sets      & categories       \\     \hline
    $\mathbf{\mathit k = 0}$
    & nonempty  & nonempty     & nonempty         \\
    & sets      & categories   & 2-categories     \\     \hline
    $\mathbf{\mathit k = 1}$
    & connected  & connected    & connected       \\
    & categories & 2-categories & 3-categories    \\     \hline
    $\mathbf{\mathit k = 2}$
    & 1-connected & 1-connected  & 1-connected    \\
    & 2-categories& 3-categories & 4-categories   \\     \hline
  \end{tabular}
\end{center}
\vskip 1em
}

In this table, the objects in the spot labeled $k$ and $n$ have
nontrivial $j$-homotopy for only $n+2$ consecutive values of $j$,
starting at $j=k$.  The column $n=-2$, which is not shown, consists
entirely of trivialities, since if you have nontrivial $j$-homotopy
for zero consecutive values of $j$, it doesn't matter at what value of
$j$ you start counting.

So we have two different periodic tables, and it isn't that one is
right and one is wrong, but rather that one is talking about monoidal
structures (or equivalently, by the delooping hypothesis, pointed
\emph{and} connected things) and the other is talking about
connectivity.  Note that unlike the monoidal periodic table, the
connectivity periodic table does not stabilize.

Finally, here's another reason to make the distinction between
`connected' and `pointed'.  We observed above that for ordinary
$n$-categories in the universe of sets, every connected $n$-category
can be made pointed in a way unique up to (non-unique) equivalence.
However, this can become false if we pass to $n$-categories in some
other universe (topos), such as `sheaves' of $n$-categories over some
space.

Consider, for instance, the relationship between groups and connected
groupoids.  A `sheaf of connected groupoids' is something called a
`gerbe' (a ``locally connected locally nonempty stack in groupoids''),
while a `sheaf of groups' is a well-known thing, but very different.
Every sheaf of groups gives rise to a gerbe, by delooping (to get a
prestack of groupoids) and then `stackifying', but it's reasonably
fair to say that the whole interest of gerbes comes from the fact that
most of them \emph{don't} come from a sheaf of groups.  The ones that
do are called `trivial', and a gerbe is trivial precisely when it has
a basepoint (a global section).  So the equivalence of groups with
\emph{pointed} connected groupoids is true even in the world of
sheaves, but in this case not every `connected' groupoid can be given
a basepoint.

If we move down one level, this corresponds to the statement that not
every well-supported sheaf has a global section.  Thus in the world of
sheaves, not every `nonempty' set can be given a basepoint.  So one
cause of the confusion between connectedness and pointedness is what
we might call `\Set-centric-ness': the two notions are quite similar
in the topos of sets, but in other topoi they are much more distinct.

\section{Annotated Bibliography}
\label{sec:annotated-bibliography}

The following bibliography should help the reader find
more detailed information about some topics mentioned
in the talks and Appendix.  It makes no pretense to completeness,
and we apologize in advance to all the authors whose work we
fail to cite.  In the spirit of `something for everybody', we 
include references with wildly different prerequisites: some
are elementary, while others even we don't understand.

\vskip 1em \noindent
\textbf{1.1  Galois Theory.}  For a
gentle introduction to Galois theory, try these:

\vskip 1em \noindent
Ian Stewart, {\sl Galois Theory}, 3rd edition, Chapman and Hall, 
New York, 2004.

\vskip 1em \noindent
Jean-Pierre Escofier, {\sl Galois Theory}, Springer, Berlin, 2000.

\vskip 1em \noindent
For more of the history, try:

\vskip 1em \noindent
Jean-Pierre Tignol, {\sl Galois' Theory of Algebraic Equations},
World Scientific, 2001.

\vskip 1em \noindent
For a treatment that emphasizes the analogy to covering spaces, try:

\vskip 1em \noindent
Adrien Douady and R\'egine Douady, {\sl Alg\`ebre et 
Th\'eories Galoisiennes}, Cassini, Paris, 2005. 

\vskip 1em \noindent
To see where the analogy between commutative algebras and
spaces went after the work of Dedekind and Kummer, try this:

\vskip 1em \noindent
Igor R.\ Shafarevich, {\sl Basic Algebraic Geometry I, II},
trans.\ M.\ Reid, Springer, Berlin, 1995/1994.

\vskip 1em \noindent
and then these more advanced but still very friendly texts:

\vskip 1em \noindent
Dino Lorenzini, {\sl An Invitation to Arithmetic Geometry}, 
American Mathematical Society, Providence, Rhode Island, 1996.

\vskip 1em \noindent
David Eisenbud and Joe Harris, {\sl The Geometry of Schemes},
Springer, Berlin, 2006.

\vskip 1em \noindent
Finally, for a very general treatment of Galois theory, try this:

\vskip 1em \noindent
Francis Borceux and George Janelidze, {\sl Galois Theories}, 
Cambridge Studies in Advanced Mathematics {\bf 72}, Cambridge
U.\ Press, Cambridge, 2001.

\vskip 1em \noindent
\textbf{1.2  The fundamental group.}  
The fundamental group is covered in almost every basic textbook on
algebraic topology.  This one is freely available and starts at a
basic level:

\vskip 1em \noindent
Allen Hatcher, {\sl Algebraic Topology}, Cambridge U.\ Press, Cambridge, 
2002.  Also available at 
\href{http://www.math.cornell.edu/~hatcher/AT/ATpage.html}
{$\langle$\texttt{http://www.math.cornell.edu/$\sim$hatcher/AT/ATpage.html}$\rangle$}.

\vskip 1em \noindent
Chapter 1 is a detailed treatment of the fundamental group
and covering spaces.

For the reader with enough prior background in category theory and
topology, the following book provides a more conceptual and
categorical approach, but it can be hard going for the novice:

\vskip 1em \noindent
Peter May, {\sl A Concise Course in Algebraic Topology}, Chicago U.\
Press, Chicago, 1999.  Also available at \hfill \break
\href{http://www.math.uchicago.edu/~may/CONCISE/ConciseRevised.pdf}
{$\langle$\texttt{http://www.math.uchicago.edu/$\sim$may/CONCISE/ConciseRevised.pdf}$\rangle$}.

\vskip 1em \noindent
A less traditional approach, also more categorical than Hatcher's, and
more closely related to modern abstract homotopy theory, can be found
in the following book:

\vskip 1em \noindent
Marcelo Aguilar, Samuel Gitler, and Carlos Prieto, {\sl Algebraic
Topology from a Homotopical Viewpoint}, Springer, Berlin, 2002.

\vskip 1em \noindent
It also provides a good concrete introduction to classifying spaces,
via covering spaces and then vector bundles.  However, the treatment
of some topics (such as homology) may strike more traditional
algebraic topologists as perverse.

\vskip 1em \noindent
\textbf{1.3  The fundamental groupoid.}
It is possible that a good modern introduction to algebraic
topology should start with the fundamental groupoid rather than
the fundamental group.  Ronnie Brown has written a text
that takes this approach:

\vskip 1em \noindent
Ronald Brown, {\sl Topology and Groupoids}, 
Booksurge Publishing, North Charleston, South Carolina, 2006. 

\vskip 1em \noindent
\textbf{1.4  Eilenberg--Mac Lane spaces.}  There are a lot of
interesting ideas packed in Eilenberg and Mac Lane's original 
series of papers on the cohomology of groups, starting around 1942 and
going on until about 1955:

\vskip 1em \noindent
Samuel Eilenberg and Saunders Mac Lane, {\sl 
Eilenberg--Mac Lane: Collected Works}, Academic Press, Orlando, 
Florida, 1986. 

\vskip 1em \noindent
These papers are a bit tough to read, but they repay the effort even 
today.  The spaces $K(G,n)$ appear implicitly in their 1945 paper 
`Relations between homology and the homotopy groups of spaces',
though much more emphasis is given on the corresponding chain
complexes.  The concept of $k$-invariant, so important for Postnikov
towers, shows up in the 1950 paper `Relations between homology and 
the homotopy groups of spaces, II'. The three papers entitled `On
the groups $H(\Pi,n)$, I, II, III' describe the bar construction
and how to compute, in principle, the cohomology groups of any
space $K(G,n)$ (where of course $G$ is abelian for $n > 1$).

The basic facts on Eilenberg--Mac Lane spaces are nicely explained
in Hatcher's {\sl Algebraic Topology} (see above).  

\vskip 1em \noindent
\textbf{1.5 Grothendieck's dream.}
The classification of general extensions of groups
goes back to Schreier:

\vskip 1em \noindent
O.\ Schreier, \"Uber die Erweiterung von Gruppen I, 
{\sl Monatschefte f\"ur Mathematik and Physik} {\bf 34} (1926), 165--180. 
\"Uber die Erweiterung von Gruppen II, 
{\sl Abh.\ Math.\ Sem.\ Hamburg} {\bf 4} (1926), 321--346.

\vskip 1em \noindent
But, the theory was worked out more thoroughly by Dedecker:

\vskip 1em \noindent
P.\ Dedecker, Les foncteuers ${\rm Ext}_\Pi, H^2_\Pi$ and $H^2_\Pi$
non abeliens, {\sl C.\ R.\ Acad.\ Sci.\ Paris} {\bf 258} (1964), 
4891--4895.

\vskip 1em \noindent
To really understand our discussion of Schreier theory, one needs
to know a bit about 2-categories.  These are good introductions:

\vskip 1em \noindent
G.\ Maxwell Kelly and Ross Street, Review of the elements of 
2-categories, Springer Lecture Notes in Mathematics {\bf 420}, 
Springer, Berlin, 1974, pp.\ 75-103. 

\vskip 1em \noindent
Ross Street, Categorical structures, in {\sl Handbook of Algebra}, 
vol.\ 1, ed.\ M.\ Hazewinkel, Elsevier, Amsterdam, 1996, pp.\ 529--577.

\vskip 1em \noindent
What we are calling `weak 2-functors' and `weak natural transformations',
they call `pseudofunctors' and `pseudonatural transformations'.

Our treatment of Schreier theory used a set-theoretic section
$s \maps B \to E$ in order to get 
an element of $H(B,\AUT(F))$ from an exact sequence
$1 \to F \to E \to B \to 1$.  The arbitrary choice of section is annoying,
and in categories other than $\Set$ it may not exist.  
Luckily, Jardine has given a construction that avoids 
the need for this splitting:

\vskip 1em \noindent
J.\ F.\ Jardine, Cocycle categories, sec.\ 4: Group extensions
and 2-groupoids, available at
\href{http://www.math.uiuc.edu/K-theory/0782/}{$\langle$\texttt{http://www.math.uiuc.edu/K-theory/0782/}$\rangle$}.

\vskip 1em 
The generalization of Schreier theory to higher dimensions has
a long and tangled history.  Larry Breen generalized it 
`upwards' from groups to 2-groups:

\vskip 1em \noindent
Lawrence Breen, Theorie de Schreier superieure, 
{\sl Ann.\ Sci.\ Ecole Norm.\ Sup.\ }{\bf 25} (1992), 465-514. 
Also available at  \hfill \break
\href{http://www.numdam.org/numdam-bin/feuilleter?id=ASENS_1992_4_25_5}
{$\langle$\texttt{http://www.numdam.org/numdam-bin/feuilleter?id=ASENS$\underline{}$1992$\underline{}$4$\underline{}$25$\underline{}$5}$\rangle$}.

\vskip 1em \noindent
It has also been generalized `sideways' from groups to groupoids:

\vskip 1em \noindent
V.\ Blanco, M.\ Bullejos and E.\ Faro, Categorical non abelian 
cohomology, and the Schreier theory of groupoids, available as 
\href{http://www.arxiv.org/abs/math.CT/0410202}{math.CT/0410202}.

\vskip 1em \noindent
However, the latter generalization is already implicit in the
work of Grothendieck: he classified all groupoids fibered over a 
groupoid $B$ in terms of weak 2-functors from $B$ to $\Gpd$, 
the 2-groupoid of groupoids.  The point is that $\Gpd$ 
contains $\AUT(F)$ for any fixed groupoid $F$:

\vskip 1em \noindent
Alexander Grothendieck, 
{\sl Rev\^etements \'Etales et Groupe Fondamental (SGA1)}, chapter VI: 
Cat\'egories fibr\'ees et descente, Lecture Notes in Mathematics 224, 
Springer, Berlin, 1971.  Also available as 
\href{http://www.arxiv.org/abs/math.AG/0206203}{math.AG/0206203}.

\vskip 1em \noindent
A categorified version of Grothendieck's result can be found here:

\vskip 1em \noindent
Claudio Hermida, Descent on 2-fibrations and strongly 2-regular
2-categories, {\sl Applied Categorical Structures}, {\bf 12} (2004),
427--459. Also available at \hfill \break
\href{http://maggie.cs.queensu.ca/chermida/papers/2-descent.pdf}
{$\langle$\texttt{http://maggie.cs.queensu.ca/chermida/papers/2-descent.pdf}$\rangle$}.

\vskip 1em \noindent
While Grothendieck was working on fibrations and `descent', 
Giraud was studying a closely related topic: nonabelian cohomology
with coefficients in a gerbe:

\vskip 1em \noindent
Jean Giraud, {\sl Cohomologie Non Ab\'elienne}, Die Grundlehren der
mathematischen Wissenschaften {\bf 179}, Springer, Berlin, 1971.

\vskip 1em \noindent
Nonabelian cohomology and $n$-categories came together in 
Grothendieck's letter to Quillen.  This is now available online along
with many other works by Grothendieck, thanks to the `Grothendieck
Circle': 

\vskip 1em \noindent
Alexander Grothendieck, {\sl Pursuing Stacks}, 1983.  Available at 
\hfill \break
\href{http://www.grothendieckcircle.org/}
{$\langle$\texttt{http://www.grothendieckcircle.org/}$\rangle$}.

\vskip 1em \noindent
Unfortunately we have not explained how these ideas are 
related to `$n$-stacks' (roughly weak sheaves of $n$-categories) and
`$n$-gerbes' (roughly weak sheaves of $n$-groupoids that are locally
connected and nonempty).  So, let us simply quote above letter:
\begin{quote}
At first sight it had seemed to me that the Bangor group had indeed
come to work out (quite independently) one basic intuition of the
program I had envisioned in those letters to Larry Breen --- namely, that
the study of $n$-truncated homotopy types (of semisimplicial sets, or of
topological spaces) was essentially equivalent to the study of so-called
$n$-groupoids (where $n$ is any natural integer).  This is expected to be
achieved by associating to any space (say) $X$ its `fundamental
$n$-groupoid' $\Pi_n(X)$, generalizing the familiar Poincar\'e fundamental
groupoid for $n = 1$.  The obvious idea is that 0-objects of $\Pi_n(X)$
should be the points of $X$, 1-objects should be `homotopies' or paths
between points, 2-objects should be homotopies between 1-objects, etc.
This $\Pi_n(X)$ should embody the $n$-truncated homotopy type of X, in much
the same way as for $n = 1$ the usual fundamental groupoid embodies the
1-truncated homotopy type.  For two spaces $X, Y,$ the set of
homotopy-classes of maps $X \to Y$ (more correctly, for general $X, Y$, the
maps of $X$ into $Y$ in the homotopy category) should correspond to
$n$-equivalence classes of $n$-functors from $\Pi_n(X)$ to $\Pi_n(Y)$ --- etc.
There are some very strong suggestions for a nice formalism including a
notion of geometric realization of an $n$-groupoid, which should imply
that any $n$-groupoid is $n$-equivalent to a $\Pi_n(X)$.  Moreover when the
notion of an $n$-groupoid (or more generally of an $n$-category) is
relativized over an arbitrary topos to the notion of an $n$-gerbe (or more
generally, an $n$-stack), these become the natural `coefficients' for a
formalism of non commutative cohomological algebra, in the spirit of
Giraud's thesis.
\end{quote}
The ``Bangor group'' led by Ronald Brown were working on
$\infty$-groupoids, but only strict ones.  For young readers, it may
be worth noting that Grothendieck's ``semisimplicial sets'' are now
called simplicial sets.

\vskip 1em \noindent
For more modern work on $n$-stacks, nonabelian cohomology and
their relation to Galois theory, try these and the
many references therein:

\vskip 1em \noindent
Andr\'e Hirschowitz and Carlos Simpson, Descente pour les 
$n$-champs, available as 
\href{http://www.arxiv.org/abs/math.AG/9807049}{math.AG/9807049}. 

\vskip 1em \noindent
Bertrand Toen, Toward a Galoisian interpretation of homotopy theory, 
available as 
\href{http://www.arxiv.org/abs/math.AT/0007157}{math.AT/0007157}.

\vskip 1em \noindent
Bertrand Toen, Homotopical and higher categorical structures in 
algebraic geometry, Habilitation thesis, Universit\'e de Nice,
2003, available as 
\href{http://www.arxiv.org/abs/math.AG/0312262}{math.AG/0312262}.

\vskip 1em \noindent
Also see the material on topoi and higher topoi
in the bibliography for Section \S\ref{sec:enrichm-posets-fiber}.

\vskip 1em \noindent
\textbf{2. The Power of Negative Thinking.}
The theory of weak $n$-categories (and $\infty$-categories) is 
in a state of rapid and unruly development, with many alternate
approaches being proposed.   For a quick sketch of the basic ideas,
try:

\vskip 1em \noindent
John C.\ Baez, An introduction to $n$-categories, in 
{\sl 7th Conference on Category Theory and Computer Science,} 
eds.\ E.\ Moggi and G.\ Rosolini, Lecture Notes in Computer Science 
{\bf 1290}, Springer, Berlin, 1997. 

\vskip 1em \noindent
John C.\ Baez and James Dolan, 
Categorification, in {\sl Higher Category Theory}, eds.\ E.\ 
Getzler and M.\ Kapranov, Contemp. Math. {\bf 230}, American 
Mathematical Society, Providence, Rhode Island, 1998, pp.\ 1--36.

\vskip 1em \noindent
For a tour of ten proposed definitions, try:

\vskip 1em \noindent
Tom Leinster, A survey of definitions of $n$-category,
available as 
\href{http://www.arxiv.org/abs/math.CT/0107188}{math.CT/0107188}.

\vskip 1em  \noindent
For more intuition on these definitions work, see this book:

\vskip 1em \noindent
Eugenia Cheng and Aaron Lauda, {\sl Higher-Dimensional Categories:
an Illustrated Guide Book}, available at 
\href{http://www.dpmms.cam.ac.uk/~elgc2/guidebook/}
{$\langle$\texttt{http://www.dpmms.cam.ac.uk/$\sim$elgc2/guidebook/}$\rangle$}.

\vskip 1em \noindent
Another useful book on this nascent subject is:

\vskip 1em \noindent
Tom Leinster, {\sl Higher Operads, Higher Categories}, 
London Math.\ Soc.\ Lecture Note Series {\bf 298}, Cambridge 
U.\ Press, Cambridge, 2004.  Also available as 
\href{http://www.arxiv.org/abs/math.CT/0305049}{math.CT/0305049}.

\vskip 1em \noindent
In the present lectures we implicitly make use of the `globular'
weak $\infty$-categories developed by Batanin:

\vskip 1em \noindent
Michael A.\ Batanin, 
Monoidal globular categories as natural environment for the 
theory of weak $n$-categories, {\sl Adv.\ Math.\ }{\bf 136} 
(1998), 39--103.

\vskip 1em \noindent
For recent progress on the homotopy hypothesis in this approach, see:

\vskip 1em \noindent
Denis-Charles Cisinski, Batanin higher groupoids and homotopy types,
in {\sl Categories in Algebra, Geometry and Mathematical Physics},
eds.\ M.\ Batanin {\it et al}, Contemp.\ Math.\ {\bf 431},
American Mathematical Society, Providence, Rhode Island, 2007, pp.\
171--186.  Also available as 
\href{http://www.arxiv.org/abs/math.AT/0604442}{math.AT/0604442}.

\vskip 1em \noindent
However, the most interesting questions about weak $n$-categories,
including the stabilization hypothesis, homotopy hypothesis and 
other hypotheses mentioned in these lectures, should ultimately
be successfully addressed by every `good' approach to the subject.
At the risk of circularity, one might even argue that this 
constitutes part of the criterion for which approaches count as `good'.

The stabilization hypothesis is implicit
in Larry Breen's work on higher gerbes:

\vskip 1em \noindent
Lawrence Breen, On the classification of 2-gerbes and 2-stacks, 
{\sl Ast\'erisque} {\bf 225}, Soci\'et\'e Math\'ematique de France, 1994.

\vskip 1em \noindent
but a blunt statement of this hypothesis, together with the 
Periodic Table, appears here:

\vskip 1em \noindent
John C.\ Baez and James Dolan, Higher-dimensional algebra and
topological quantum field theory, 
{\sl Jour.\ Math.\ Phys.\ }{\bf 36} (1995), 6073--6105.
Also available as \href{http://www.arxiv.org/abs/q-alg/9503002}
{q-alg/9503002}.

\vskip 1em \noindent
There has been a lot of progress recently toward
precisely formulating and proving the stabilization hypothesis
and understanding the structure of $k$-tuply monoidal $n$-categories
and their relation to $k$-fold loop spaces:

\vskip 1em \noindent
Carlos Simpson, On the Breen--Baez--Dolan stabilization hypothesis for 
Tamsamani's weak $n$-categories, available as 
\href{http://www.arxiv.org/abs/math.CT/9810058}{math.CT/9810058}.

\vskip 1em \noindent
Michael A.\ Batanin, The Eckmann--Hilton argument and higher operads, 
available as 
\href{http://www.arxiv.org/abs/math.CT/0207281}{math.CT/0207281}.

\vskip 1em \noindent
Michael A.\ Batanin, The combinatorics of iterated loop spaces,
available as \break 
\href{http://www.arxiv.org/abs/math.CT/0301221}{math.CT/0301221}.

\vskip 1em \noindent
Eugenia Cheng and Nick Gurski, The periodic table of $n$-categories for
low dimensions I: degenerate categories and degenerate bicategories, in
{\sl Categories in Algebra, Geometry and Mathematical Physics},
eds.\ M.\ Batanin {\it et al}, Contemp.\ Math.\ {\bf 431},
American Mathematical Society, Providence, Rhode Island, 2007, pp.\
143--164.
Also available as \href{http://www.arxiv.org/abs/0708.1178}{arXiv:0708.1178}.

\vskip 1em \noindent
Eugenia Cheng and Nick Gurski, The periodic table of $n$-categories for
low dimensions II: degenerate tricategories, available as
\href{http://www.arxiv.org/abs/0706.2307}{arXiv:0706.2307}.

\vskip 1em \noindent
The mathematical notion of `stuff' was introduced here: 

\vskip 1em \noindent
John C.\ Baez and James Dolan, From finite sets to Feynman diagrams,
in {\sl Mathematics Unlimited - 2001 and Beyond}, vol.\ 1, eds.\
Bj\o rn Engquist and Wilfried Schmid, Springer, Berlin, 2001, pp.\ 29--50.
Also available as \href{http://www.arxiv.org/abs/math.QA/0004133}
{math.QA/0004133}.

\vskip 1em \noindent
where `stuff types' (groupoids over the groupoid of finite sets
and bijections) were used to explain the combinatorial underpinnings
of the theory of Feynman diagrams.  A more detailed study of this
subject can be found here:

\vskip 1em \noindent
John C.\ Baez and Derek Wise, {\sl Quantization and Categorification},
Quantum Gravity Seminar, U.\ C.\ Riverside, Spring 2004 lecture notes, 
available at
\hfill \break
\href{http://math.ucr.edu/home/baez/qg-spring2004/}
{$\langle$\texttt{http://math.ucr.edu/home/baez/qg-spring2004/}$\rangle$}.

\vskip 1em \noindent
On this page you will find links to a pedagogical introduction to 
properties, structure and stuff by Toby Bartels, and also
to a long online conversation in which $(-1)$-categories and
$(-2)$-categories were discovered.   See also:

\vskip 1em \noindent
Simon Byrne, {\sl On Groupoids and Stuff}, honors thesis, 
Macquarie University, 2005, available at 
\href{http://www.maths.mq.edu.au/~street/ByrneHons.pdf}
{$\langle$\texttt{http://www.maths.mq.edu.au/$\sim$street/ByrneHons.pdf}$\rangle$}
\hfill \break
and 
\href{http://math.ucr.edu/home/baez/qg-spring2004/ByrneHons.pdf}
{$\langle$\texttt{http://math.ucr.edu/home/baez/qg-spring2004/ByrneHons.pdf}$\rangle$}.

\vskip 1em \noindent
Jeffrey Morton, Categorified algebra and quantum mechanics, 
available as \hfill \break 
\href{http://www.arxiv.org/abs/math.QA/0601458}{math.QA/0601458}.

\vskip 1em \noindent
\textbf{3. Cohomology: The Layer-Cake Philosophy.}  In topology, 
it's most common to generalize the basic principle of 
Galois theory from covering spaces to fiber bundles along these
lines:
\begin{quote}
\textbf{Principal $G$-bundles over a base space $B$
are classified by maps from $B$ to the classifying space $BG$.}
\end{quote}
For example, if $B$ is a CW complex and $G$ is a
topological group, then isomorphism classes of principal 
$G$-bundles over $M$ are in one-to-one correspondence with
homotopy classes of maps from $B$ to $BG$.  
Good references on this theory include:

\vskip 1em \noindent
John Milnor and James Stasheff, {\sl Characteristic classes},
Ann.\ Math.\ Studies {\bf 76}, Princeton U.\ Press, Princeton, 1974.
\vskip 1em \noindent
Dale Husemoller, {\sl Fibre Bundles}, Springer, Berlin, 1993.
\vskip 1em

Another approach, more in line with higher category
theory, goes roughly as follows:

\begin{quote}
\textbf{Fibrations over a pointed connected base space $B$ with fiber 
$F$ are classified
by homomorphisms sending based loops in $B$ to automorphisms of $F$.}
\end{quote}

\noindent
Stasheff proved a version of this which is reviewed here:

\vskip 1em \noindent
James Stasheff, H-spaces and classifying spaces, I-IV, 
{\sl AMS Proc.\ Symp.\ Pure Math.\ }{\bf 22} (1971), 247--272.
\vskip 1em

\noindent
He treats the space $\Omega B$ of based loops in $B$ as an $A_\infty$ space,
i.e.\ a space with a product that is associative up to a homotopy
that satisfies the pentagon identity up to a homotopy that satisfies a 
further identity up to a homotopy... ad infinitum.
He classifies fibrations over $B$ with fiber $F$ in terms of 
$A_\infty$-morphisms from $\Omega B$ into the topological monoid 
$\Aut(F)$ consisting of homotopy equivalences of $F$.

Another version was proved here:

\vskip 1em \noindent
J.\ Peter May, {\sl Classifying Spaces and Fibrations}, AMS Memoirs
{\bf 155}, American Mathematical Society, Providence, 1975.
\vskip 1em

\noindent
Moore loops in $B$ form a topological monoid $\Omega_{\rm M} B$.
May defines a {\bf transport} to be a 
homomorphism of topological monoids from $\Omega_{\rm M} B$ 
to $\Aut(F)$.
After replacing $F$ by a suitable homotopy-equivalent space, he defines
an equivalence relation on transports such that the equivalence 
classes are in natural one-to-one correspondence with the equivalence 
classes of fibrations over $B$ with fiber $F$.

By iterating the usual classification of principal $G$-bundles over $B$
in terms of maps $B \to BG$, we obtain the theory of Postnikov towers.
A good exposition of this can be found at the end of Chapter 4 of 
Hatcher's book {\sl Algebraic Topology}, already cited in the notes for 
Section \S\ref{sec:fundamental group}.
Unfortunately this treatment, like most expository accounts, limits itself
to `simple' spaces, namely those which $\pi_1$ acts trivially on the higher
homotopy groups.  For the general case see:

\vskip 1em \noindent
C.\ Alan Robinson, Moore--Postnikov systems for non-simple fibrations,
{\sl Ill.\ Jour.\ Math.\ }{\bf 16} (1972), 234--242.  

\vskip 1em \noindent
For a treatment of Postnikov towers based on simplicial sets rather
than topological spaces, try:

\vskip 1em \noindent
J.\ Peter May, {\sl Simplicial Objects in Algebraic Topology}, 
Van Nostrand, Princeton, 1968.  

\vskip 1em\noindent
Hatcher also discusses the cohomology groups of Eilenberg--Mac Lane 
spaces.  Elements of these are called {\bf cohomology operations}, 
and more information on them can be found here:

\vskip 1em \noindent
Norman E.\ Steenrod and David B.\ A. Epstein, {\sl Cohomology Operations},
Princeton U.\ Press, Princeton, 1962.

\vskip 1em \noindent
Robert E.\ Mosher and Martin C.\ Tangora, {\sl Cohomology 
Operations and Applications in Homotopy Theory}, Harper and Row,
New York, 1968.

\vskip 1em \noindent
Given $m > n > 1$,
elements of $H^m(K(G,n),A) = [K(G,n),K(A,m)]$
classify connected `simple spaces' with only $\pi_n$ and $\pi_{m-1}$
nontrivial.  We can think of these as $(m-1)$-groupoids with only
two nontrivial layers.  So, a lot of information about higher categories
can be dug out of cohomology operations.  For example, we have seen that
when $m = 3$ and $n = 1$, elements of $H^m(K(\Pi,n),A)$ classify 
possible associators for 2-groups.  When $m = 4$ and $n = 2$, they
classify possible associators and braidings for braided 2-groups.
When $m = 5$ and $n = 3$, they classify the same thing for symmetric
2-groups.  The pattern becomes evident upon consulting the periodic
table.  An understanding of this theory was what led Breen to notice
a flaw in Kapranov and Voevodksy's original definition of braided
monoidal 2-category.

\vskip 1em \noindent
Here is the paper by Street on cohomology with coefficients in 
an $\infty$-category:

\vskip 1em \noindent
Ross Street, Categorical and combinatorial aspects of descent theory,
available at 
\href{http://www.arxiv.org/abs/math.CT/0303175}{math.CT/0303175}.

\vskip 1em \noindent
The idea of cohomology with coefficients in an $\infty$-category
seems to have originated here:

\vskip 1em \noindent
John E.\ Roberts, Mathematical aspects of local cohomology, in
{\sl Alg\`ebres d'Op\'erateurs et Leurs Applications 
en Physique Math\'ematique,} CNRS, Paris, 1979, pp.\ 321--332.

\vskip 1em \noindent
\textbf{4. A Low-Dimensional Example.}  This section will make
more sense if one is comfortable with the cohomology of 
groups.  To get started, try:

\vskip 1em \noindent
Joseph J.\ Rotman, {\sl An Introduction to Homological Algebra}, 
Academic Press, New York, 1979.

\vskip 1em \noindent
or this more advanced book with the same title:

\vskip 1em \noindent
Charles A.\ Weibel, {\sl An Introduction to Homological Algebra},
Cambridge U.\ Press, Cambridge, 1995.

\vskip 1em \noindent
For more detail, we recommend:

\vskip 1em \noindent
Kenneth S.\ Brown, {\sl Cohomology of Groups}, 
Graduate Texts in Mathematics {\bf 182}, Springer, Berlin, 1982.

\vskip 1em
The classification of 2-groups up to equivalence using group cohomology
was worked out by a student of Grothendieck whom everybody 
calls `Madame Sinh':

\vskip 1em \noindent
Hoang X.\ Sinh, {\sl Gr-categories}, Universit\'e Paris VII
doctoral thesis, 1975.

\vskip 1em \noindent
She called them {\bf gr-categories} instead of 2-groups, and this
terminology remains common in the French literature.
Her thesis, while very influential, was never published.  
Later, Joyal and Street described the whole 2-category of
2-groups using group cohomology here:

\vskip 1em \noindent
Andr\'e Joyal and Ross Street, Braided monoidal categories,
Macquarie Mathematics Report No.\ 860081, November 1986.
Also available at \hfill \break
\href{http://rutherglen.ics.mq.edu.au/~street/JS86.pdf}
{$\langle$\texttt{http://rutherglen.ics.mq.edu.au/$\sim$street/JS86.pdf}$\rangle$}.

\vskip 1em \noindent
Joyal and Street call them {\bf categorical groups} instead
of 2-groups.  Like Sinh's thesis, this paper was never published 
--- the published paper with a similar title leaves out
the classification of 2-groups and moves directly to the 
classification of braided 2-groups.  These are also nice examples 
of the general `layer-cake philosophy' we are discussing
here.  Since braided 2-groups are morally the same as connected pointed
homotopy types with only $\pi_2$ and $\pi_3$ nontrivial, 
their classification involves $H^4(K(G,2),A)$ (for $G$ abelian) instead
of the cohomology group we are considering here, $H^3(K(G,1),A) =
H^3(G,A)$.   For details, see:

\vskip 1em \noindent
Andr\'e Joyal and Ross Street, Braided tensor categories,
{\sl Adv.\ Math.\ }{\bf 102} (1993), 20--78.

\vskip 1em \noindent
Finally, since it was hard to find a clear treatment of the classification 
of 2-groups in the published literature, an account was included here:

\vskip 1em \noindent
John C.\ Baez and Aaron Lauda, Higher-dimensional algebra V: 
2-Groups, {\sl Th.\ Appl.\ Cat.\ }{\bf 12} (2004), 423--491.
Also available as 
\href{http://www.arxiv.org/abs/math.QA/0307200}{math.QA/0307200}.

\vskip 1em \noindent
along with some more history of the subject.  
The analogous classification of Lie 2-algebras using
Lie algebra cohomology appears here:

\vskip 1em \noindent
John C.\ Baez and Alissa S.\ Crans, Higher-dimensional algebra VI: 
Lie 2-Algebras, {\sl Th.\ Appl.\ Cat.\ }{\bf 12} (2004), 492--528.
Also available as 
\href{http://www.arxiv.org/abs/math.QA/0307263}{math.QA/0307263}.

\vskip 1em \noindent
This paper also shows that an element of $H^{n+1}_{\rho}(\mathfrak{g},
\mathfrak{a})$ gives a Lie $n$-algebra with $\mathfrak{g}$ as the Lie
algebra of objects and $\mathfrak{a}$ as the abelian Lie algebra of
$(n-1)$-morphisms.  

\vskip 1em \noindent
\textbf{5.  Appendix: Posets, Fibers, and $n$-Topoi.}
As mentioned in the notes to Section \S\ref{sec:Grothendieck},
fibrations as functors between categories were introduced by
Grothendieck in SGA1.  They were then extensively developed by Jean
B\'enabou, and his handwritten notes of a 1980 course {\sl Des
Categories Fibr\'ees} have been very influential.  Unfortunately these
remain hard to get, so we suggest:

\vskip 1em \noindent
Thomas Streicher, Fibred categories \`a la Jean B\'enabou, April 2005,
available as \href{http://www.mathematik.tu-darmstadt.de/~streicher/FIBR/FibLec.pdf.gz}{$\langle$\texttt{http://www.mathematik.tu-darmstadt.de/$\sim$streicher/FIBR/FibLec.pdf.gz}$\rangle$}.

\vskip 1em \noindent
Street's fully weakened definition of fibrations in general weak
2-categories can be found here:

\vskip 1em \noindent
Ross Street, Fibrations in bicategories, {\sl Cah.\ Top.\ Geom.\
Diff.\ Cat.\ }{\bf 21} (1980), 111--160.  Errata, {\sl Cah.\ Top.\
Geom.\ Diff.\ Cat.\ }{\bf 28} (1987), 53--56.

\vskip 1em \noindent There are several good introductions to topos
theory; here are a couple:

\vskip 1em \noindent
Saunders Mac Lane and Ieke Moerdijk, {\sl Sheaves in Geometry and 
Logic: a First Introduction to Topos Theory}, Springer, New York, 1992.

\vskip 1em \noindent
Colin McLarty, {\sl Elementary Categories, Elementary Toposes},
Oxford U.\ Press, Oxford, 1992.

\vskip 1em \noindent
The serious student will eventually want to spend time with the {\sl
  Elephant}:

\vskip 1em \noindent
Peter Johnstone, {\sl Sketches of an Elephant: a Topos Theory Compendium}, 
Oxford U. Press, Oxford.  Volume 1, comprising Part A: Toposes as 
categories, and Part B: 2-Categorical aspects of topos theory, 2002.  
Volume 2, comprising Part C: 
Toposes as spaces, and Part D: Toposes as theories, 
2002.  Volume 3, comprising Part E: Homotopy and 
cohomology, and Part F: Toposes as mathematical universes, in preparation.

\vskip 1em \noindent
The beginning of part C is a good introduction to
locales and Heyting algebras.  We eagerly await part E for
more illumination on one of the main themes of this paper, 
namely the foundations of cohomology theory.  

The idea of defining sheaves as categories enriched over a 
certain bicategory is due to Walters:

\vskip 1em \noindent
Robert F.\ C.\ Walters, Sheaves on sites as Cauchy-complete categories,
{\sl J.\ Pure Appl.\ Algebra} {\bf 24} (1982), 95--102
\vskip 1em

Street has some papers on cosmoi:

\vskip 1em \noindent
Ross Street and Robert F.\ C.\ Walters, Yoneda structures on 2-categories,
{\sl J.\ Algebra} {\bf 50} (1978), 350--379.

\vskip 1em \noindent
Ross Street, Elementary cosmoi, I,
{\sl Category Seminar}, 
Lecture Notes in Math. {\bf 420}, Springer, Berlin, 1974,
pp.\ 134--180. 

\vskip 1em \noindent
Ross Street, Cosmoi of internal categories,
{\sl Trans.\ Amer.\ Math.\ Soc.\ }{\bf 258} (1980), 271--318.
\vskip 1em

\noindent
and Weber's 2-topos paper can now be found here:

\vskip 1em \noindent
Mark Weber, Yoneda structures from 2-toposes, {\sl Applied Categorical
Structures} {\bf 15} (2007), 259--323.  Slightly different
version available as 
\href{http://www.arxiv.org/abs/math.CT/0606393}{math.CT/0606393}.
\vskip 1em

The theory of $\infty$-topoi --- or what we prefer to call
$(\infty,1)$-topoi, to emphasize the room left for further
expansion --- is new and still developing.  A book just came out
on this subject: 

\vskip 1em\noindent
Jacob Lurie, Higher topos theory, available as 
\href{http://www.arxiv.org/abs/math.CT/0608040}{math.CT/0608040}.
\vskip 1em

\noindent
For a good introduction to $\infty$-topoi which takes a 
slightly nontraditional approach to topoi, see:

\vskip 1em \noindent
Charles Rezk, Toposes and homotopy toposes, available at \hfill \break
\href{http://www.math.uiuc.edu/~rezk/homotopy-topos-sketch.dvi}
{$\langle$\texttt{http://www.math.uiuc.edu/$\sim$rezk/homotopy-topos-sketch.dvi}$\rangle$}.
\vskip 1em

\noindent
This paper is an overview of $\infty$-topoi using model 
categories and Segal categories: 

\vskip 1em \noindent
Bertrand Toen, Higher and derived stacks: a global overview, 
available as \hfill \break 
\href{http://www.arxiv.org/abs/math.AG/0604504}{math.AG/0604504}.
\vskip 1em \noindent

The correspondence between properties of small functors and properties
of geometric morphisms is, to our knowledge, not written down all
together anywhere.  Johnstone summarized it in his talk at the Mac
Lane memorial conference in Chicago in 2006.  One can extract this
information from {\sl Sketches of an Elephant} if one looks at the
examples in the sections on various types of geometric morphism in
part C, and always think `modulo splitting idempotents'.

Classifying topoi are explained somewhat in Mac Lane by Moerdijk, and more
in part D of the {\sl Elephant}.  We also recommend having a look at the 
version of classifying topoi in part B of sketches, which uses 2-categorical
limits to construct them.  This makes the connection with $n$-stuff a
little clearer.

Another good introduction to classifying topoi, and their relationship
to topology, is:

\vskip 1em \noindent
Steven Vickers, Locales and toposes as spaces, 
available at \hfill \break
\href{http://www.cs.bham.ac.uk/~sjv/LocTopSpaces.pdf}
{$\langle$\texttt{http://www.cs.bham.ac.uk/$\sim$sjv/LocTopSpaces.pdf}$\rangle$}.
\vskip 1em

\end{document}